\documentclass[11pt,reqno]{amsart}
\usepackage[T1]{fontenc}
\usepackage{amsmath,amssymb,mathtools}
\usepackage[a4paper,textwidth=158mm,textheight=235mm,centering]{geometry}
\usepackage[colorlinks=true,linkcolor=blue,citecolor=blue,urlcolor=blue]{hyperref}
\numberwithin{equation}{section}
\allowdisplaybreaks[2]
\renewcommand{\paragraph}[1]{%
  \par\addvspace{.6\baselineskip}\noindent
  {\normalfont\itshape #1}\enspace\ignorespaces}
\newtheorem{theorem}{Theorem}[section]
\newtheorem{proposition}[theorem]{Proposition}
\newtheorem{lemma}[theorem]{Lemma}
\newtheorem{corollary}[theorem]{Corollary}
\theoremstyle{definition}

\theoremstyle{remark}
\newtheorem{remark}[theorem]{Remark}
\newcommand{\R}{\mathbb R}
\newcommand{\C}{\mathbb C}
\newcommand{\ii}{\mathrm i}
\newcommand{\e}{\mathrm e}
\newcommand{\Rea}{\operatorname{Re}}
\newcommand{\Ima}{\operatorname{Im}}
\newcommand{\diver}{\operatorname{div}}
\newcommand{\rad}{\mathrm{rad}}
\newcommand{\dd}{\,\mathrm d}
\newcommand{\pair}[2]{\langle #1,#2\rangle}
\newcommand{\norm}[1]{\lVert #1\rVert}
\newcommand{\abs}[1]{\lvert #1\rvert}
\newcommand{\U}{\vec U}
\newcommand{\E}{\mathcal E}

\newcommand{\Dtwo}{D_2}
\newcommand{\eps}{\varepsilon}
\newcommand{\epsn}{\varepsilon_n}
\newcommand{\epsv}{\varepsilon_v}
\newcommand{\veps}{\vec\varepsilon}
\title[Asymptotic stability of Zakharov blow-up]{Asymptotic stability of self-similar blow-up\newline for the two-dimensional Zakharov system}
\author[Ghoul]{Tej-eddine Ghoul}
\author[Liu]{Jie Liu}
\author[Masmoudi]{Nader Masmoudi}
\thanks{All authors are affiliated with New York University Abu Dhabi,
Abu Dhabi, United Arab Emirates.
Email addresses: \href{mailto:teg6@nyu.edu}{teg6@nyu.edu} (T. Ghoul),
\href{mailto:jl15817@nyu.edu}{jl15817@nyu.edu} (J. Liu),
\href{mailto:nm30@nyu.edu}{nm30@nyu.edu} (N. Masmoudi).}
\date{5 November 2026}

\begin{document}
\begin{abstract}
We prove asymptotic stability of the small-parameter Glangetas--Merle
self-similar blow-up solutions of the two-dimensional Zakharov system
under nonradial perturbations in the energy space. The proof
establishes dissipation on the wave time scale by canceling the fast
density coupling. We impose nonlinear orthogonality to the ground
state of the density-dependent Schr\"odinger operator for modulation
and use successive energy corrections after subtracting the
instantaneous density response from the Schr\"odinger error. The
resulting dissipation yields power-law decay of the rescaled errors
in the Schr\"odinger gradient, weighted local mass, and wave energy
norms.
\end{abstract}
\maketitle

\section{Introduction}

\subsection{Background and self-similar blow-up solutions}

We consider the Zakharov system in two space dimensions,
\begin{equation}\label{eq:physical}
 \begin{cases}
 \ii\partial_t\psi+\Delta\psi=n\psi,\\
 \partial_t n=-\diver v,\\
 \partial_t v=-c_0^2\nabla(n+\abs\psi^2),
 \end{cases}
 \qquad (t,x)\in[0,T)\times\R^2,\quad c_0>0.
\end{equation}
Here $\psi$ is the complex envelope of the high-frequency electric
field of the Langmuir waves, $n$ is the real low-frequency density
perturbation about the equilibrium plasma density, and $v$ is the
normalized ion velocity field in the continuity equation.
The constant $c_0$ is the fixed ion-acoustic speed.
The density perturbation $n$ acts as a potential for the Langmuir
envelope, while the ponderomotive force associated with the intensity
$|\psi|^2$ drives the ion-acoustic response. This feedback can concentrate
Langmuir waves in low-density regions and underlies the collapse mechanism
described by Zakharov \cite{Zakharov}.

The energy space is
\begin{equation}\label{eq:full-physical-space}
 H^1(\R^2;\C)\times L^2(\R^2;\R)
                           \times L^2(\R^2;\R^2).
\end{equation}
The energy solutions conserve the mass $\norm\psi_2^2$
and the Hamiltonian
\begin{equation}\label{eq:ham-physical}
 \int\abs{\nabla\psi}^2+\int n\abs\psi^2
       +\frac12\int n^2+\frac1{2c_0^2}\int\abs v^2.
\end{equation}
Unless otherwise specified, spatial integrals are over $\R^2$.

\paragraph{Self-similar collapse and virial estimates.}
Self-similar wave collapse was studied by Zakharov and Shur
\cite{ZS}. In numerical simulations with radial and anisotropic initial data,
Landman, Papanicolaou, Sulem, Sulem and Wang \cite{LPSSW}
observed convergence toward isotropic self-similar collapse.
Glangetas and Merle \cite{GM1} subsequently constructed a branch
of finite-energy radial self-similar solutions in dimension two.
Their second paper \cite{GM2} established concentration properties
of blow-up solutions and instability results for standing waves.
Merle's virial and rate estimates describe further restrictions
imposed by the wave equation. In \cite{MerleVirial}, localized virial identities imply that radial
negative-energy solutions either blow up in finite time or have
unbounded norms along a sequence of times tending to infinity. In \cite{MerleRate}, Merle proved that
every finite-time blow-up solution in the two-dimensional energy
class satisfies a lower bound of order $(T-t)^{-1}$ for
$\norm{\nabla\psi(t)}_2$ and $\norm{n(t)}_2$. The self-similar
solutions realize this order. On the two-dimensional torus, Kishimoto
and Maeda \cite{KM} constructed finite-time blow-up by localizing the
Euclidean profile and correcting the resulting approximate solution.
Wang \cite{WangPeriodic} constructs a different periodic blow-up
mechanism below the Schr\"odinger regularity $H^{1/2}$, in which the
wave component blows up while the Schr\"odinger norm remains bounded.
In the Euclidean setting, Cher, Czubak and Sulem
\cite[Remark~2.3, arXiv version]{CherCzubakSulem}
identified dynamical stability as an open problem for both
two-dimensional self-similar solutions and three-dimensional
asymptotically self-similar solutions.
Theorem~\ref{thm:main} resolves the two-dimensional question
on the small-parameter Glangetas--Merle branch.
In three dimensions, Yin \cite{Yin3D} constructed smooth radial
blow-up solutions with an asymptotically self-similar core at spatial
scale $(T-t)^{2/3}$ and left their stability open in Remark~1.3.

\paragraph{Four-dimensional blow-up and the two time scales.}
For the energy-critical four-dimensional system, Candy, Herr and
Nakanishi \cite[Theorem~1.1]{CHN4D} proved global well-posedness for
data with energy below the ground-state energy and wave $L^2$ norm
at most its ground-state value. Candy and Nakanishi \cite{CN4D}
established scattering in the potential well. Krieger and Schmid \cite{KS1,KS2} constructed
type II blow-up by correcting an approximate solution built from the
four-dimensional cubic ground state, under numerical nondegeneracy
assumptions. Their proof uses separate Schr\"odinger and wave clocks,
whose rates are the inverse square and inverse first power of the
concentrating length. Integration by parts in Schr\"odinger time
relates the principal nonresonant Schr\"odinger error to the density
through an elliptic inverse. Feeding this response back into the wave equation yields
an effective wave operator with a spatially nonlocal term. The
intermediate wave temporal frequencies are handled by a Fredholm
inversion \cite[Sections~2, 11 and~12]{KS2}.

The same separation of time scales and the elliptic response of the
Schr\"odinger component are relevant to our two-dimensional problem.
Here the wave evolution is retained in Schr\"odinger time, with its
small scale coefficient. We subtract the instantaneous response
using the reduced inverse at the full rescaled density and impose
orthogonality to the ground state of that same operator. The
resulting cancellation and the corrections to the perturbation energy
yield a forward dissipation estimate on the wave time scale, which
controls an open neighborhood of the exact self-similar branch.

\paragraph{The mass-critical NLS comparison.}
Eliminating the wave variable formally in the limit $c_0\to\infty$
gives $n=-|\psi|^2$ and the focusing cubic NLS
\[
 \ii\partial_t\psi+\Delta\psi+|\psi|^2\psi=0
                       \qquad\hbox{on }\R^2.
\]
Its positive radial ground state $Q$ solves
$-\Delta Q+Q-Q^3=0$. Its uniqueness is due to \cite{Kwong}, and
its mass is the threshold in the sharp interpolation inequality
\cite{Weinstein}. Merle
\cite{MerleMinimal} classified minimal-mass blow-up solutions of
the mass-critical NLS; they arise from the ground state through
the pseudoconformal symmetry. Bourgain and Wang \cite{BW}
constructed solutions with a pseudoconformal core and a regular
remainder. These special regimes are distinct from the stable
log--log regime, whose construction began with Perelman's
one-dimensional work \cite{Perelman}.

Merle and Rapha\"el \cite{MRSharp,MRProfile,MRUpper,MRLower}
and Rapha\"el \cite{RaphaelLoglog} established the stable log--log
regime and convergence of its rescaled core to the ground state
for a class of blow-up data with mass slightly above the ground-state
mass. Merle and Rapha\"el \cite{MRQuantization} also proved a
decomposition into a concentrating core and a remainder with a
strong $L^2$ limit for a class of mass-critical blow-up solutions.
Fibich, Merle and Rapha\"el \cite{FMR} verified the spectral property
used by the log--log theory in dimension two. In the log--log regime,
\[
 \norm{\nabla\psi(t)}_2
       \asymp\left(\frac{\log|\log(T-t)|}{T-t}\right)^{1/2}.
\]
Merle and Rapha\"el explicitly discuss the comparison with the
Zakharov system in \cite[Section~1.1, comment~5]{MRSharp}.
For every fixed finite sound speed, Merle's Zakharov lower bound
excludes the NLS log--log rate. This is a structural instability
under a change of equation, even at large sound speed; it does
not concern perturbations of the data for a fixed Zakharov system.

For slightly mass-supercritical NLS, Merle, Rapha\"el and Szeftel
\cite{MRS} proved stability of a self-similar blow-up regime. Their
argument combines a local virial inequality with a second
monotonicity formula based on mass conservation and radiation
flux and establishes dynamical trapping. In the Zakharov problem,
the distinct Schr\"odinger and wave time scales must be retained
in the nonlinear error estimates.

\paragraph{Exact self-similar solutions.}
In our normalization, the Glangetas--Merle branch
$(P_a,N_a,V_a)$, for sufficiently small $a>0$, satisfies
\begin{align}
 &(-\Delta+1+N_a)P_a=0,\qquad
 V_a=-c_0a\,yN_a,\label{eq:profiles}\\
 &\diver V_a=-c_0a\Dtwo N_a,\qquad
 \nabla(N_a+P_a^2)=-\frac a{c_0}\Dtwo V_a,
 \qquad \Dtwo=2+y\cdot\nabla.\label{eq:wave-profiles}
\end{align}
The operator $\Dtwo$ acts componentwise on vectors, $P_a>0$, and
$P_a\to Q$ in $H^1$ as $a\downarrow0$. For a scale $\lambda>0$,
write $\beta=c_0a\lambda$ and define
\begin{equation}\label{eq:quadratic-phase-profile}
 Q_{a,\beta}(y)=P_a(y)\e^{-\ii\beta|y|^2/4}.
\end{equation}
The coefficient $\beta$ is determined by $a$ and $\lambda$.
Fix $a_0>0$, $\lambda_0>0$, $\gamma_0\in\R$ and $x_0\in\R^2$.
A superscript $\star$ denotes the fixed reference solution and its
parameters. For $0\le t<T_0=\lambda_0/(c_0a_0)$, set
\begin{equation}\label{eq:exact-orbit}
 \begin{aligned}
 \lambda^\star(t)&=\lambda_0-c_0a_0t,&
 \beta^\star(t)&=c_0a_0\lambda^\star(t),\\
 \gamma^\star(t)&=\gamma_0+
 \frac1{c_0a_0}\left(\frac1{\lambda^\star(t)}
                                      -\frac1{\lambda_0}\right).
 \end{aligned}
\end{equation}
The corresponding exact self-similar solution is
\begin{equation}\label{eq:exact-reference}
 \begin{aligned}
 \psi^\star(t,x)&=(\lambda^\star)^{-1}
       \e^{\ii\gamma^\star}
       Q_{a_0,\beta^\star}\left(\frac{x-x_0}{\lambda^\star}\right),\\
 n^\star(t,x)&=(\lambda^\star)^{-2}
       N_{a_0}\left(\frac{x-x_0}{\lambda^\star}\right),\\
 v^\star(t,x)&=(\lambda^\star)^{-2}
       V_{a_0}\left(\frac{x-x_0}{\lambda^\star}\right).
 \end{aligned}
\end{equation}
The identities \eqref{eq:profiles}--\eqref{eq:wave-profiles}
verify \eqref{eq:physical} directly. In particular,
$\partial_t\gamma^\star=(\lambda^\star)^{-2}$.
We study an open energy neighborhood of the initial datum of
\eqref{eq:exact-reference}, allowing all blow-up parameters to change.

\subsection{Spectral estimates and main theorem}

For complex functions we use the $L^2$ inner product
$\pair h g=\int h\overline g$, linear in its first argument, and
write $(h,g)_\R=\Rea\pair h g$ for the associated real pairing.
For real functions, $(f,g)=\int fg$; for real vectors,
$(v,w)=\int v\cdot w$. Angular brackets for complex vectors use
the same componentwise convention. When a source is a distribution,
these brackets mean its continuous dual pairing with the indicated
test, with the same conjugation convention. Thus real density and
velocity pairings use parentheses, while complex Schr\"odinger
projections use angular brackets.
Define
\begin{equation}\label{eq:local-norm}
 \norm h_{\mathrm{loc}}^2
       =\norm{\nabla h}_2^2+\int\e^{-|y|}|h(y)|^2\dd y.
\end{equation}
It measures the global gradient and mass near the concentrating
profile, and is weaker than the full $H^1$ norm.
The spaces $H^1_\rad(\R^2;\C)$ and
$L^2_\rad(\R^2;\R)$ consist of scalar functions depending only
on $r=|y|$. For vector fields we use the convention
\[
 L^2_\rad(\R^2;\R^2)
 =\left\{w\in L^2(\R^2;\R^2):
          w(y)=w_r(|y|)\frac y{|y|}\ \hbox{a.e. for some }w_r\right\}.
\]
These
subspaces are used in the intermediate radial argument; the
main theorem allows general nonradial perturbations.

Write $Q'$ for the radial derivative of $Q$, and set
\[
 L_-=-\Delta+1-Q^2,\qquad A_Q=-2\Delta+2rQQ'.
\]
The nonnegativity and simple kernel $\C Q$ of $L_-$ are classical
\cite[Proposition~2.8(a)]{WeinsteinModulation}. Its domain is
$H^2(\R^2)$; the decaying potential gives essential spectrum
$[1,\infty)$, so zero is isolated and the inverse on $Q^\perp$
is bounded. We use $A_Q$ as a quadratic form on $H^1$.

\begin{proposition}[Spectral property at the limiting profile]\label{ass:spectral}
There is $\delta_*>0$ for which the following two inequalities hold.

\textup{(i) Schr\"odinger coercivity.} For radial
$\eps\in H^1(\R^2;\C)$ with $\pair\eps Q=0$,
\begin{equation}\label{eq:spectral-scalar}
 (A_Q\eps,\eps)_\R\ge\delta_*\norm\eps_{\mathrm{loc}}^2.
\end{equation}

\textup{(ii) Density coercivity.} For real radial
$\epsn\in L^2$ with $(\epsn,Q^2)=0$, let $z\in H^2\cap Q^\perp$
be the unique solution of $L_-z=Q\epsn$. Then
\begin{equation}\label{eq:spectral-density}
 \norm {\epsn}_2^2+(A_Qz,z)-2({\epsn},(Q-rQ')z)
                       \ge\delta_*\norm {\epsn}_2^2.
\end{equation}

\end{proposition}

The condition $(\epsn,Q^2)=0$ makes $Q\epsn$ orthogonal to the kernel
of $L_-$, so the equation defining $z$ is compatible. Both
inequalities concern the limiting ground state $Q$.
Proposition~\ref{prop:verify-spectral} in Appendix~\ref{app:certificate}
proves both estimates using interval bounds for explicit integrals.
The appendix also treats the first angular density mode under the
translation constraints, the higher angular modes, and the velocity
estimates.

\begin{theorem}[Nonradial stability and convergence to a self-similar profile]
\label{thm:main}
There exists $a_*>0$ such
that, for every $a_0\in(0,a_*)$, there is $\lambda_*(a_0)>0$
with the following property. Fix $0<\lambda_0<\lambda_*(a_0)$,
$\gamma_0\in\R$, $x_0\in\R^2$, and the exact solution
\eqref{eq:exact-reference}. For every sufficiently small $\eta>0$
there exists $\delta>0$ such that each initial datum
$(\psi_{\mathrm{in}},n_{\mathrm{in}},v_{\mathrm{in}})$ in
\eqref{eq:full-physical-space} satisfying
\[
 \norm{\psi_{\mathrm{in}}-\psi^\star(0)}_{H^1}
 +\norm{n_{\mathrm{in}}-n^\star(0)}_2
 +c_0^{-1}\norm{v_{\mathrm{in}}-v^\star(0)}_2<\delta
\]
generates an energy solution of \eqref{eq:physical} with finite
maximal forward existence time $T$ and the following properties.

\textup{(i) Decomposition and orbital bound.}
There are parameters $\lambda(t)>0$, $a(t)>0$,
$\gamma(t)\in\R$ and $x_c(t)\in\R^2$ such that, with
$\beta(t)=c_0a(t)\lambda(t)$, $\dd s/\dd t=\lambda^{-2}$
and $y=(x-x_c(t))/\lambda(t)$,
\begin{equation}\label{eq:main-error}
 \begin{aligned}
 \psi(t,x)&=\lambda^{-1}\e^{\ii\gamma}
                              (Q_{a,\beta}+\eps)(s,y),\\
 n(t,x)&=\lambda^{-2}(N_a+\epsn)(s,y),\\
 v(t,x)&=\lambda^{-2}(V_a+\epsv)(s,y),
 \end{aligned}
\end{equation}
and
\begin{equation}\label{eq:main-orbital-bound}
 \sup_{0\le t<T}
 \left(\norm{\eps}_{H^1}+\norm{\epsn}_2
       +c_0^{-1}\norm{\epsv}_2+|a(t)-a_0|\right)<\eta.
\end{equation}
\textup{(ii) Convergence of the profile, scale and center.}
There are $a_\infty\in(0,a_*)$ and $x_\infty\in\R^2$ such that
\begin{align}
 a(t)&\longrightarrow a_\infty,\label{eq:main-limits}\\
 \frac{\lambda(t)}{T-t}&\longrightarrow c_0a_\infty.
 \label{eq:main-speed}
\end{align}
There is $\sigma>0$, uniform on the chosen neighborhood, such that
\begin{align}
 &\norm{\eps}_{\mathrm{loc}}+\norm{\epsn}_2
   +c_0^{-1}\norm{\epsv}_2+|a-a_\infty|
                   =O((T-t)^\sigma),\label{eq:main-power-rate}\\
 &\lambda(t)=c_0a_\infty(T-t)+O((T-t)^{1+\sigma})
                                  \qquad(t\uparrow T).
 \label{eq:main-scale-remainder}
\end{align}

For the center,
\begin{equation}\label{eq:main-center-limit}
 x_c(t)\longrightarrow x_\infty,\qquad
 \frac{|x_c(t)-x_\infty|}{\lambda(t)}\longrightarrow0,
\end{equation}
and, for the same choice of $\sigma$,
\begin{equation}\label{eq:main-center-power}
 |x_c(t)-x_\infty|=O((T-t)^{1+\sigma}).
\end{equation}
The constants may depend on the reference parameters and the
initial datum.
\end{theorem}

The choice of the decay exponent and the corresponding estimates
are given in Proposition~\ref{prop:quantitative-rates}.

\begin{remark}[Topology of convergence]\label{rem:topology}
The bound \eqref{eq:main-orbital-bound} gives orbital stability in
the full rescaled energy space. The attraction in
\eqref{eq:main-power-rate} uses the weaker norm
\eqref{eq:local-norm} for the Schr\"odinger component. More precisely,
\begin{equation}\label{eq:intro-radiation}
 \lim_{t\uparrow T}\norm{\eps(t)}_2^2
       =\norm{\psi_{\mathrm{in}}}_2^2-\norm{P_{a_\infty}}_2^2\ge0,
\end{equation}
by mass conservation and the vanishing profile overlap; see
Corollary~\ref{cor:radiation-mass}. Global strong $L^2$ convergence
to zero holds exactly when this mass difference vanishes.
Arbitrarily small perturbations can have positive limiting
radiation mass, as shown by the construction following that
corollary. Positive radiation mass and a zero local limit preclude
strong global $L^2$ convergence. This distinction also occurs in the
mass-critical NLS theory \cite[Remark~3]{MRProfile}.
\end{remark}

\begin{remark}[A quantitative choice of decay exponent]\label{rem:explicit-rate}
After reducing the upper bound of the small positive branch, one may
take $\sigma=1/300$ in Theorem~\ref{thm:main}. More generally, for any
fixed $0<\sigma<1/50$, one can choose an upper branch bound and then
the initial scale and perturbation thresholds so that the theorem
holds with this exponent. The branch bound and exponent can be
chosen independently of $c_0>0$ and of $a_0$ in the resulting branch.
The initial thresholds and the constants in the rates may depend on
the reference parameters. The comparison giving these exponents is
proved in Section~\ref{sec:explicit-rate}.
\end{remark}

\subsection{Outline of the proof}
\leavevmode\par

\paragraph{The two time scales.}
We first work with radial perturbations. At length $\lambda$, the
physical time scales are $\lambda^2$ for Schr\"odinger evolution and
$\lambda/c_0$ for wave propagation. In Schr\"odinger time
$\dd s/\dd t=\lambda^{-2}$, the wave couplings in \eqref{eq:rescaled}
carry $\lambda$ and $c_0^2\lambda$, while the density coupling in the
Schr\"odinger equation has coefficient one. The leading dissipation is
therefore sought at order $\lambda$. To use it, we must first cancel
the fast density contributions in the modulation and energy estimates.

\paragraph{Nonlinear modulation.}
Write $(\eps,\epsn,\epsv)$ for the errors relative to
$(Q_{a,\beta},N_a,V_a)$, with $\beta=c_0a\lambda$, and let $\phi(N)$ be
the normalized ground state of $L_N=-\Delta+1+N$, where $N=N_a+\epsn$.
The condition $\pair\eps{\phi(N)}=0$ eliminates
$\pair{L_N\eps}{\phi(N)}$, including the fast product $\epsn\eps$.
The projection of the remaining forcing $Q_{a,\beta}\epsn$ combines
with the phase frequency, leaving a contribution of order $\lambda$
times the error in the amplitude equation. In differentiating the
condition, the dilation of $\phi(N)$ cancels that of $\eps$, so the
wave source $-\lambda\diver V$ supplies the factor $\lambda$ in the
remaining derivative of the test. Together with $(\epsn,P_a^2)=0$,
this gives the bounds \eqref{eq:rate-bound} for $a_s$ and the scale
defect $\lambda_s/\lambda+\beta$.

\paragraph{Energy corrections.}
The fixed-orbit estimate in Section~\ref{sec:linear-energy} motivates
the construction. The interaction energy cancels the exchange between
the Schr\"odinger and wave components, but contains the global mass
$\|\eps\|_2^2$. Subtracting this mass directly brings the fast density
term back into the derivative. We remove the instantaneous density
response first and subtract the mass of
\[
 z=\eps+R_N(Q_{a,\beta}\epsn),\qquad
 R_N=\bigl(L_N|_{\phi(N)^\perp}\bigr)^{-1}\Pi_N,
\]
where $\Pi_N$ is the orthogonal projection onto $\phi(N)^\perp$.
Applying $L_N$ to the correction removes the projected density forcing
from the equation for $z$. The component along $\phi(N)$ has zero
pairing with $z$. The wave equation and the modulation bounds
give the factor $\lambda$ in the derivative of the response.

The mass subtraction leaves mixed terms and work from $a_s$.
The leading part of $a_s$ is itself a pairing with $z$, so the
same quadratic correction can cancel both. It pairs $\Ima\eps$
with a localized test linear in the wave errors. Differentiating
this pairing with the fast Schr\"odinger equation supplies the
opposite terms
\eqref{eq:direct-correction-derivative}. The scale work is included in
the velocity dissipation form.

\paragraph{Dissipation and radial closure.}
The corrected energy $\E$ controls the Schr\"odinger gradient and
weighted local mass, as well as the wave errors in $L^2$.
Both the leading dissipation and the nonlinear remainder in
\eqref{eq:energy-virial} carry $\beta$, which is comparable to
$\lambda$ for $a$ in a fixed positive compact interval $I$.
Thus fixed smallness of the perturbation suffices for absorption,
even as $\lambda$ tends to zero. The resulting inequalities are
\begin{equation}\label{eq:intro-energy}
 \E_s\le-\kappa_I\lambda\E,\qquad
 |a_s|+\left|\frac{\lambda_s}{\lambda}+\beta\right|
                         \le K_I\lambda\sqrt\E .
\end{equation}
Integrating the inequality for $\sqrt\E$ gives
\begin{equation}\label{eq:intro-integral}
 \int_{s_0}^S\lambda\sqrt\E\,\dd s
                 \le\frac{2}{\kappa_I}\sqrt{\E(s_0)} .
\end{equation}
This bounds the total variation of $a$ and keeps the profile on
the chosen interval. Control of the profile mass and conservation of
the total Schr\"odinger mass then keep $\|\eps\|_2$ small. Together
with the scale equation, these bounds close the bootstrap and give
the blow-up and convergence rates in Section~\ref{sec:closure}.

\paragraph{Nonradial perturbations.}
The limiting first angular density dissipation form vanishes when
evaluated on a translation direction. The conditions $(\epsn,\partial_jN_a)=0$,
$j=1,2$, fix the center and exclude these directions. Its motion also
changes the leading velocity form. Section~\ref{sec:nonradial}
proves coercivity with this term included, recovering the same
energy and scale estimates and an integrable bound for the center
velocity. The radial closure then extends to general perturbations.

\section{Notation}\label{sec:notation}

The profile mass is $M=M(a)=\norm{P_a}_2^2$.
Choose the small branch bound $a_*>0$ as required by the profile
and spectral estimates. Then fix numbers
\[
 0<a_-<a_0<a_+<a_*,
 \qquad I=[a_-,a_+].
\]
Here $a_0$ is the parameter of the exact solution in the theorem;
$a_-$ and $a_+$ are fixed endpoints for the modulation parameter.
All parameter bounds below are taken on $I$.
Constants $C_I$ are uniform in $a\in I$ and sufficiently small
$\lambda$; they may depend on the fixed sound speed and weights,
and may change between estimates. Constants used simultaneously
in the bootstrap are fixed before the initial neighborhood is chosen.

\subsection{Rescaled fields and error spaces}

The radial argument uses the subspace
\begin{equation}\label{eq:physical-space}
 H^1_\rad(\R^2;\C)\times L^2_\rad(\R^2;\R)
                 \times L^2_\rad(\R^2;\R^2),
\end{equation}
with the scalar and vector conventions given in the introduction.
Section~\ref{sec:nonradial} removes these symmetry restrictions.

For a positive scale $\lambda$, a phase $\gamma$ and a center $x_c$,
use the rescaling
\begin{equation}\label{eq:rescaling}
 \begin{aligned}
 &\psi(t,x)=\lambda^{-1}\e^{\ii\gamma}\Psi(s,y),\quad
 n(t,x)=\lambda^{-2}N(s,y),\quad v(t,x)=\lambda^{-2}V(s,y),\\
 &y=\frac{x-x_c}{\lambda},\qquad \frac{\dd s}{\dd t}=\lambda^{-2}.
 \end{aligned}
\end{equation}
The full state is $\U=(\Psi,N,V)$. We use
$\Lambda=1+y\cdot\nabla$ for the Schr\"odinger dilation and
$\Dtwo=2+y\cdot\nabla$ for the density and velocity dilations.
The latter acts componentwise on vectors.

The profile family $Q_{a,\beta}$ has independent arguments $a$
and $\beta$. The modulation parameters are $a,\lambda,\gamma$
and, for nonradial data, $x_c$. Along the modulation we evaluate
the profile at $\beta=c_0a\lambda$. Any additional quadratic phase
is part of $\eps$.

With the local norm \eqref{eq:local-norm}, the spaces for the local estimates are
\begin{equation}\label{eq:X}
 \begin{aligned}
 X&=\{h\in H^1_{\mathrm{loc},\rad}(\R^2;\C):
                         \norm h_{\mathrm{loc}}<\infty\},\\
 Y&=X\times L^2_\rad(\R^2;\R)\times L^2_\rad(\R^2;\R^2).
 \end{aligned}
\end{equation}
The wave factors are real. Despite the word ``local'', the $Y$
norm controls the global field gradient, a weighted local field
mass and both global wave $L^2$ norms. The full energy norm also contains $\norm\eps_2$, which is
propagated separately by mass conservation.
For weighted Sobolev norms, smoothly extend the multiplier
$\e^{d|y|}$ from $|y|\ge1$ across the unit ball, keeping it positive.
We denote this fixed extension by $\e^{dr}$ and set
$\norm g_{H^j_d}=\norm{\e^{dr}g}_{H^j}$. In particular, the
multiplier outside the ball is $\e^{d|y|}$. In pointwise estimates
and radial integrals, $r$ retains its ordinary meaning $|y|$.
Fix $1/2<d<d'<1$. The stronger decay $d'$ absorbs polynomial
factors when an estimate is made with weight $d$. The reference,
error and norm are
\begin{equation}\label{eq:error}
\begin{aligned}
 \U_*&=(Q_{a,\beta},N_a,V_a),\qquad
 \veps=(\eps,\epsn,\epsv)=\U-\U_*,\\
 \norm{\veps}_Y^2&=\norm{\eps}_{\mathrm{loc}}^2
 +\norm{\epsn}_2^2+c_0^{-2}\norm{\epsv}_2^2.
\end{aligned}
\end{equation}
The reference state $\U_*$ is the comparison triple at the current
modulation parameters $(a,\lambda)$; it varies along a perturbed
solution. It is distinct from the fixed exact solution marked by
$\star$ in the introduction. The stationary profiles are
$(P_a,N_a,V_a)$, and $Q_{a,\beta}$ includes the explicit quadratic
phase. Neither should be confused with $\phi(N)$, the normalized
ground state of the linear operator with the full potential
$N$. The relation at the profile density is
$\phi(N_a)=P_a/\sqrt M$. Unless a star is present, ``reference''
below means the current comparison state $\U_*$.

The same definitions without radial restrictions are used in
Section~\ref{sec:nonradial}. We write $r=|y|$ and
$\langle y\rangle=(1+|y|^2)^{1/2}$; the notation $\langle r\rangle$
denotes the same weight.

\subsection{Energy solutions and tested identities}

We use the solution supplied by the local Cauchy theory
in $C_t(H^1\times L^2\times L^2)$, as specified in
Appendix~\ref{app:energy-approximation}. The equations give
\[
 \psi_t,v_t\in C_tH^{-1},\qquad n_t\in C_t\dot H^{-1}\subset C_tH^{-1}.
\]
The homogeneous estimate for $n_t$ follows from
$n_t=-\diver v$ and the gradient projection. For the
inhomogeneous tested bounds, let $\eta\in H^1(\R^2)$. Then
\[
 \begin{aligned}
 |\pair{n\psi}{\eta}|
   &\le\|n\|_2\|\psi\|_4\|\eta\|_4
       \le C\|n\|_2\|\psi\|_{H^1}\|\eta\|_{H^1},\\
 \|\diver v\|_{H^{-1}}&\le\|v\|_2,\\
 \|\nabla(n+|\psi|^2)\|_{H^{-1}}
   &\le\|n\|_2+\|\psi\|_4^2 .
 \end{aligned}
\]
In particular, neither $\nabla n$ nor $\nabla v$ is assumed to
belong to $L^2$. Dilation terms in rescaled coordinates are
interpreted by their distributional pairings.

We perform differentiated calculations first for smooth solutions.
Spatial integrations by parts use the cutoffs justified in
Appendix~\ref{app:energy-approximation}. The resulting canceled
expressions extend continuously to energy solutions. Derivatives of an error in a localized pairing
act on its test function. For example,
\[
 \begin{aligned}
 (\diver\epsv,\eta)&=-(\epsv,\nabla\eta),&
 (\Dtwo\epsn,\eta)&=-(\epsn,y\cdot\nabla\eta),\\
 \pair{-\Delta\eps}{p}
      &=\pair{\nabla\eps}{\nabla p},&
 \pair{\Lambda\eps}{p}
      &=\pair\eps p+\pair{\nabla\eps}{yp}.
 \end{aligned}
\]
The tests used below have the indicated weighted Sobolev
regularity. In particular the mixed nonlinear products are
paired with the localized $H^1$ tests established below. The final passage to energy solutions uses convergence
of smooth \emph{solutions} on compact time intervals, after the
identities have been integrated. That passage is proved in
Lemma~\ref{lem:energy-solutions}.

\subsection{Operators and differentiation}

Set $L_a=-\Delta+1+N_a$. Its nonnegativity follows directly from
the positive profile equation: for $h\in C_c^\infty$,
\[
 (L_ah,h)_\R=\int P_a^2\left|\nabla\left(\frac h{P_a}\right)\right|^2.
\]
This ground-state identity extends by approximation. Its kernel is
$\C P_a$, and the decaying potential $N_a$ gives essential spectrum
$[1,\infty)$, hence a gap on $P_a^\perp$. This is the same
argument as for $L_-$ in \cite[Proposition~2.8(a)]{WeinsteinModulation}.
The gap remains uniform on a fixed positive compact interval $I$
by continuity of the potentials and of the simple eigenprojection.
For real $N$ close to $N_a$ in $L^2$, let $\mu(N)$ and $\phi(N)$ be the
simple ground-state eigenvalue and the positive normalized ground state
of $-\Delta+N$. At the reference,
\begin{equation}\label{eq:phi0}
 \mu(N_a)=-1,\qquad \phi(N_a)=P_a/\sqrt M.
\end{equation}
For a real $L^2$ density direction $f$, the notation
$\phi'_N[f]$ means
$\left.\frac{\mathrm d}{\mathrm d\tau}\phi(N+\tau f)\right|_{\tau=0}$;
$N$ is the base density, and higher primes have the same convention.
The map is differentiated from $L^2$ into the weighted spaces
specified in Lemma~\ref{lem:spectral}.
For $V\in L^2(\R^2;\R^2)$ and a complex test $\eta\in H^1$,
its divergence is defined by
$\pair{\diver V}{\eta}=-\int V\cdot\nabla\overline\eta$.
For real tests this is equivalently
$(\diver V,\eta)=-(V,\nabla\eta)$.
Lemma~\ref{lem:spectral} defines $\phi'_N[\diver V]$ by solving
the differentiated eigenvalue equation with this $H^{-1}$ source.
Only this continuous extension is used; the eigenfunction map
itself is differentiated in an $L^2$ neighborhood.

Let $\Pi_a$ be the orthogonal projection onto $P_a^\perp$, and set
\begin{equation}\label{eq:RK}
 R_a=(L_a|_{P_a^\perp})^{-1}\Pi_a,\qquad
 K_\beta {\epsn}=-R_a(Q_{a,\beta}{\epsn}).
\end{equation}
A profile next to an operator denotes pointwise multiplication.
The density-response operator has full dependence
$K_{a,\beta}f=-R_a(P_a\e^{-\ii\beta r^2/4}f)$; we suppress $a$
and evaluate $\beta=c_0a\lambda$. At zero quadratic phase the
response to a density direction $f$ is $-R_a(P_af)$.

We also suppress $a$ in the phase-indexed operator $\mathcal A_\beta$
and the quadratic forms $E_\beta^{(2)}$, $\mathcal E_\beta$ and
$B_\beta$. In $\mathcal A_\beta$ the wave coefficient is written
as $\lambda=\beta/(c_0a)$ to retain its rescaling interpretation.
For $K_{a,\beta}$, each partial derivative holds the other
profile argument fixed. In the fixed-orbit calculation only $\beta$
varies in this operator, so $(K_\beta)_s=\beta_s\partial_\beta K_{a,\beta}$.
The conversion between the two coefficients is always
\[
 \beta=c_0a\lambda,\qquad
 c_0a_-\lambda\le\beta\le c_0a_+\lambda,\qquad
 \frac{\lambda}{\beta}=\frac1{c_0a}.
\]
A conversion from a $\lambda$ bound to a $\beta$ bound may
therefore cost $(c_0a_-)^{-1}$. Such constants are allowed only
after $I$ is fixed. In estimates used to choose the small branch
bound $a_*$, constants are uniform as $a\downarrow0$ and the
powers of $a$ are retained explicitly.

The zero-phase quadratic form $\mathcal E_0$ is the limit
of $\mathcal E_\beta$ as $\beta\downarrow0$ at fixed $a$ and
fixed error. This sets $Q_{a,\beta}=P_a$ and removes the correction
proportional to $\beta$. The nonlinear energy at the initial time
is written $\mathcal E(s_0)$.

We use the pairings fixed in the introduction. Arrows distinguish triples of field, density and velocity components
from their scalar components. Full-state
dualities are sums of the Schr\"odinger, density and velocity
component dualities; the subscript $\R$ takes the real part of
the complex Schr\"odinger contribution.
For a real quadratic form $B$ we use the symmetric bilinear
polarization
\[
 B(\vec u,\vec h)=\tfrac14\{B(\vec u+\vec h)-B(\vec u-\vec h)\},
 \qquad dB(\vec u)[\vec h]=2B(\vec u,\vec h).
\]
All differentials below are real Fr\'echet differentials:
the complex Schr\"odinger component is regarded as a real Banach
space. Thus, for a scalar or vector valued map $F$, a base state
$\U$ and real directions $\vec h_1,\ldots,\vec h_j$ in its stated domain,
\[
 d^jF(\U)[\vec h_1,\ldots,\vec h_j]
 =\left.\frac{\partial^j}{\partial\tau_1\cdots\partial\tau_j}
 F\left(\U+\sum_{\ell=1}^j\tau_\ell \vec h_\ell\right)
                                                  \right|_{\tau=0},
 \qquad \tau_\ell\in\R.
\]
For an error functional at fixed $(a,\lambda)$, the base variable
is $\veps\in Y$ instead. The partial derivatives of \eqref{eq:quadratic-phase-profile}
hold the other profile argument fixed:
\[
 \partial_a Q_{a,\beta}=\e^{-\ii\beta r^2/4}\partial_a P_a,\qquad
 \partial_\beta Q_{a,\beta}=-\tfrac{\ii r^2}{4}Q_{a,\beta}.
\]
The total $s$ derivative and the derivative of the phase coefficient are
\begin{equation}\label{eq:profile-time-derivative}
 \begin{aligned}
 (Q_{a,\beta})_s&=a_s\partial_a Q_{a,\beta}+\beta_s\partial_\beta Q_{a,\beta},\\
 \beta_s+\beta^2&=\beta\left(\frac{\lambda_s}{\lambda}+\beta\right)
                   +\frac{\beta}{a}a_s .
 \end{aligned}
\end{equation}
Spatial subscripts on
$\partial_j$ refer to differentiation in $y_j$. Parameter
subscripts on profiles, operators and quadratic forms
record dependence, rather than differentiation. Continuity of a
differential on the full energy space or on $Y$ is specified separately.

\section{Modulation equations and bootstrap conditions}\label{sec:coordinates}

We first consider radial perturbations with center zero.
The modulation conditions use the ground state associated with the
full rescaled density. Differentiating these conditions gives the
slow amplitude and scale estimates needed by the energy argument.
We then fix the bootstrap regime.
The center is added in Section~\ref{sec:nonradial}.

The profile mass is strictly increasing on the small positive branch:
$M'(a)>0$. Proposition~\ref{prop:slope} in
Appendix~\ref{app:mass-slope} gives its asymptotic coefficient and proof.
This supplies the nonzero branch derivative in the parameter map below.

\subsection{Choice of modulation parameters}

Fix $I\Subset(0,a_*)$ as in Section~\ref{sec:notation}, now with
$M'>0$ by Proposition~\ref{prop:slope}. The radial rescaled equations are
\begin{equation}\label{eq:rescaled}
 \begin{aligned}
 \Psi_s&=-\ii(-\Delta+N+\gamma_s)\Psi
                      +\frac{\lambda_s}{\lambda}\Lambda\Psi,\\
 N_s&=\frac{\lambda_s}{\lambda}\Dtwo N-\lambda\diver V,\\
 V_s&=\frac{\lambda_s}{\lambda}\Dtwo V
                     -c_0^2\lambda\nabla(N+\abs\Psi^2).
 \end{aligned}
\end{equation}
Here $\Psi,N,V$ are the total rescaled fields. Along the exact path
with constant $a$ and $\lambda_s/\lambda=-c_0a\lambda$, the wave
equations have the coefficient $\lambda$, while the leading
Schr\"odinger operator has coefficient one. These are the slow and
fast evolutions at a fixed spatial length in $y$, both measured in $s$.

The error equation contains the product $\epsn\eps$ with
coefficient one. Testing against $P_a$ removes the reference
operator $L_a$ but leaves this nonlinear density term. We use
the ground state associated with $N$ so that the whole
self-adjoint operator disappears from the differentiated constraint.
For given physical data and candidate parameters, use
$\eps=\Psi-Q_{a,\beta}$ and $\epsn=N-N_a$, and impose
\begin{equation}\label{eq:chart}
 \begin{cases}
 \Rea\pair\eps{\phi(N)}=0,\\
 \Ima\pair\eps{\phi(N)}=0,\\
 (\epsn,P_a^2)=0.
 \end{cases}
\end{equation}
Thus the Schr\"odinger error is orthogonal to $\phi(N)$.
These are three real equations for $(a,\gamma,\lambda)$; they select
coordinates on a neighborhood of the reference orbit. The first two
are nonlinear in the density through $\phi(N)$.

\begin{lemma}[Choice of the parameters]\label{lem:chart}
Fix $a_0\in I$, sufficiently small $\lambda_0>0$ and
$\gamma_0\in\R$. In the radial energy space
\eqref{eq:physical-space}, an open neighborhood of the exact initial
datum \eqref{eq:exact-reference}, with center zero, admits a unique
choice of $(a,\lambda,\gamma)$ near $(a_0,\lambda_0,\gamma_0)$
satisfying \eqref{eq:chart}. The phase is chosen as a real number
near $\gamma_0$, excluding the other representatives modulo $2\pi$.
The neighborhood contains every sufficiently small perturbation in
the physical $H^1\times L^2\times L^2$ norm.
\end{lemma}

\begin{proof}
\paragraph{Profile identities.}
Commuting $L_a=-\Delta+1+N_a$ with $\Lambda=1+y\cdot\nabla$
and differentiating $L_aP_a=0$ give
\[
 [L_a,\Lambda]=-2\Delta-y\cdot\nabla N_a,\qquad
 L_a\partial_aP_a=-P_a\partial_aN_a.
\]
Since $\Delta P_a=(1+N_a)P_a$, pairing the resulting equations
with $P_a$ yields
\begin{equation}\label{eq:pivots}
 L_a\Lambda P_a=-2P_a-P_a\Dtwo N_a,\quad
 (\Dtwo N_a,P_a^2)=-2M,\quad(\partial_aN_a,P_a^2)=0.
\end{equation}

\paragraph{The constraint map on physical data.}
The normalized eigenfunction obeys
\[
 \phi(\lambda^2n(\lambda\cdot))=\lambda\phi(n)(\lambda\cdot).
\]
This follows by scaling its eigenvalue equation; positivity fixes
its sign. Consequently,
\begin{equation}\label{eq:physical-chart-tests}
 \begin{aligned}
 \pair\eps{\phi(N)}
 &=\e^{-\ii\gamma}\pair\psi{\phi(n)}
    -\pair{\lambda^{-1}Q_{a,\beta}(\cdot/\lambda)}{\phi(n)},\\
 (\epsn,P_a^2)
 &=\int n(x)P_a(x/\lambda)^2\dd x-(N_a,P_a^2).
 \end{aligned}
\end{equation}
Use the real and imaginary parts of the first expression and the
second expression as the constraint map. The physical eigenpair
near the reference datum is simple and isolated. The proof of
Lemma~\ref{lem:spectral}, after the fixed reference scaling, gives
$C^1$ dependence of $\phi(n)$ on $n\in L^2$; see also
\cite[Chapter~VII]{Kato}. The profile
$\lambda^{-1}Q_{a,\beta}(\cdot/\lambda)$ and its parameter
derivatives belong to $H^1$ and vary continuously near the fixed
positive scale. Thus every derivative in the first line acts on
a smooth profile or on the bounded eigenfunction differential.
It does not spatially differentiate the rough physical data.

For the density expression, use $\tau=\log\lambda$. Parameter
derivatives act on the smooth tests, with
\[
 \begin{aligned}
 \|P_a(\cdot/\lambda)^2\|_2&=\lambda\|P_a^2\|_2,\\
 \|\partial_\tau[P_a(\cdot/\lambda)^2]\|_2
       &=\lambda\|y\cdot\nabla(P_a^2)\|_2,\\
 \|\partial_a[P_a(\cdot/\lambda)^2]\|_2
       &=2\lambda\|P_a\partial_aP_a\|_2.
 \end{aligned}
\]
The constraint map is therefore $C^1$ in data and parameters.
The velocity ranges over its full $L^2$ space.

\paragraph{The parameter Jacobian.}
At the exact datum, $\eps=0$ and $\phi(N_a)=P_a/\sqrt M$.
With physical data fixed, differentiate the constraint map in
$(a,\gamma,\tau)$. Since
$\partial_a(c_0a\lambda)=\beta/a$ at fixed $\lambda$ and
$\partial_\tau(c_0a\e^\tau)=\beta$ at fixed $a$, this gives
\[
\begin{aligned}
 \partial_a\pair\eps{\phi(N)}
   &=-\pair{\partial_a Q_{a,\beta}+(\beta/a)\partial_\beta Q_{a,\beta}}
            {\phi(N_a)},\\
 \partial_\gamma\pair\eps{\phi(N)}
   &=-\ii\pair{Q_{a,\beta}}{\phi(N_a)},\\
 \partial_\tau\pair\eps{\phi(N)}
   &=\pair{\Lambda Q_{a,\beta}-\beta\partial_\beta Q_{a,\beta}}{\phi(N_a)}=0.
\end{aligned}
\]
The terms differentiating the test vanish because $\eps=0$.
The last equality follows from the explicit quadratic phase and
integration by parts:
\[
 \pair{\Lambda Q_{a,\beta}}{\phi(N_a)}
  =-\pair{Q_{a,\beta}}{\Lambda\phi(N_a)}
  =\beta\pair{\partial_\beta Q_{a,\beta}}{\phi(N_a)}.
\]
The density row has derivatives $(0,0,-2M)$ by \eqref{eq:pivots}.
The real Jacobian in $(a,\gamma,\tau)$ is therefore
\[
 \begin{pmatrix}
 -M'/(2\sqrt M)&0&0\\
 0&-\sqrt M&0\\
 0&0&-2M
 \end{pmatrix}+O_I(\lambda).
\]
The two profile pairings in the $a$ and $\gamma$ entries are
$M'/(2\sqrt M)+O_I(\beta)$ and $\sqrt M+O_I(\beta)$, respectively.
Their corrections involve only polynomial moments of
exponentially localized profiles. The positive lower
bounds for $M$ and $M'$ give a uniformly bounded inverse at small
reference scale. Continuity and the implicit-function theorem
then give the asserted local choice of parameters on an open
neighborhood of the datum.
\end{proof}

Scaling the ground-state equation and preserving its positive
$L^2$ normalization gives
\begin{equation}\label{eq:covariance}
 \phi(N_\sigma)(y)=\sigma\phi(N)(\sigma y),\qquad
 N_\sigma(y)=\sigma^2N(\sigma y).
\end{equation}
We use this covariance to differentiate the moving test in
\eqref{eq:chart} without taking a strong derivative of the rough density.

\subsection{Slow modulation equations}
\label{sec:slow-modulation}

We now differentiate the chosen conditions. The main issue is
the coefficient of the field error: its fast self-adjoint part
must cancel before estimating the moving ground state.
The resulting amplitude and scale rates carry $\lambda$;
the phase instead follows the ground-state eigenvalue associated with $N$.

\begin{lemma}[Modulation identities and estimates]\label{lem:modulation}
Write $\phi=\phi(N)$.
For smooth solutions satisfying \eqref{eq:chart},
\begin{equation}\label{eq:orthogonality-evolution}
\begin{aligned}
 a_s\pair{\partial_a Q_{a,\beta}+(\beta/a)\partial_\beta Q_{a,\beta}}\phi
 &-\left(\frac{\lambda_s}{\lambda}+\beta\right)
             \pair{\Lambda Q_{a,\beta}-\beta\partial_\beta Q_{a,\beta}}\phi\\
 &+\ii(\mu(N)+\gamma_s)\pair{Q_{a,\beta}}\phi
  =\ii\pair{L_aQ_{a,\beta}}\phi-\lambda\pair\eps{\phi'_N[\diver V]}.
\end{aligned}
\end{equation}
The density condition gives
\begin{equation}\label{eq:density-eq}
\begin{aligned}
 [-2M+(\Dtwo\epsn,P_a^2)]
            \left(\frac{\lambda_s}{\lambda}+\beta\right)
  ={}&\lambda(\diver\epsv,P_a^2)+\beta(\Dtwo\epsn,P_a^2)\\
    &-2a_s(\epsn,P_a\partial_aP_a).
\end{aligned}
\end{equation}
For sufficiently small scale and error, these equations imply
\begin{align}
 |a_s|+|\lambda_s/\lambda+\beta|
       &\le C_I\lambda\|\veps\|_Y,\label{eq:rate-bound}\\
 |\gamma_s-1|&\le C_I(\beta\|\veps\|_Y+\|\veps\|_Y^2).
 \label{eq:phase-bound}
\end{align}
The phase frequency also satisfies
\begin{equation}\label{eq:phase-equation}
 |\mu(N)+\gamma_s|\le C_I\lambda\|\veps\|_Y.
\end{equation}
For energy solutions the identities hold in the tested integrated
sense of Section~\ref{sec:notation}.
\end{lemma}

\begin{proof}
\paragraph{Differentiating the error orthogonality.}
The covariance \eqref{eq:covariance} and normalization give
\[
 \phi'_N[\Dtwo N]=\Lambda\phi,
 \qquad (\phi,\phi'_N[\diver V])=0.
\]
For rough densities, $(\Dtwo N)\phi$ means
$\diver(yN\phi)-N\,y\cdot\nabla\phi$. The weighted estimates
in Lemma~\ref{lem:spectral} justify the paired covariance and
the one-divergence direction by approximation. Thus
\[
 \phi_s=\frac{\lambda_s}{\lambda}\Lambda\phi
                         -\lambda\phi'_N[\diver V].
\]
Subtracting the current reference field from \eqref{eq:rescaled},
with $L_N=-\Delta+1+N$, gives
\[
\begin{aligned}
 \eps_s={}&-\ii(L_N+\gamma_s-1)\eps
      -\ii Q_{a,\beta}\epsn-\ii(\gamma_s-1)Q_{a,\beta}
      +\frac{\lambda_s}{\lambda}\Lambda\eps\\
 &+\left(\frac{\lambda_s}{\lambda}+\beta\right)\Lambda Q_{a,\beta}
       -a_s\partial_a Q_{a,\beta}
       -(\beta_s+\beta^2)\partial_\beta Q_{a,\beta}.
\end{aligned}
\]
Here we used
$-\ii L_aQ_{a,\beta}-\beta\Lambda Q_{a,\beta}=-\beta^2\partial_\beta Q_{a,\beta}$.
In $0=\partial_s\pair\eps\phi$, the self-adjoint field term
vanishes since
\[
 \pair{(L_N+\gamma_s-1)\eps}\phi
       =(\mu(N)+\gamma_s)\pair\eps\phi=0.
\]
The two dilation terms cancel by $\Lambda^*=-\Lambda$. The
remaining density projection has the exact form
\[
 \pair{Q_{a,\beta}\epsn}\phi
       =(1+\mu(N))\pair{Q_{a,\beta}}\phi-\pair{L_aQ_{a,\beta}}\phi.
\]
Using \eqref{eq:profile-time-derivative} to express
$\beta_s+\beta^2$ in terms of the two rates in the linear system
gives \eqref{eq:orthogonality-evolution}.
The factor $\pair{Q_{a,\beta}}\phi$ is nonzero near the reference state.
Dividing by it and taking real parts removes the real phase
frequency and yields the amplitude equation
\begin{equation}\label{eq:amplitude-eq}
\begin{aligned}
 a_s\Rea\frac{\pair{\partial_a Q_{a,\beta}+(\beta/a)\partial_\beta Q_{a,\beta}}\phi}{\pair{Q_{a,\beta}}\phi}
 &-\left(\frac{\lambda_s}{\lambda}+\beta\right)
       \Rea\frac{\pair{\Lambda Q_{a,\beta}-\beta\partial_\beta Q_{a,\beta}}\phi}
                       {\pair{Q_{a,\beta}}\phi}\\
 &=\Rea\frac{\ii\pair{L_aQ_{a,\beta}}\phi
                -\lambda\pair\eps{\phi'_N[\diver V]}}
                       {\pair{Q_{a,\beta}}\phi}.
\end{aligned}
\end{equation}

\paragraph{The density row and its coupling to the amplitude equation.}
The density error satisfies
\[
 (\epsn)_s=\left(\frac{\lambda_s}{\lambda}+\beta\right)\Dtwo N_a
       +\frac{\lambda_s}{\lambda}\Dtwo\epsn
       -\lambda\diver\epsv-a_s\partial_aN_a.
\]
Differentiating $(\epsn,P_a^2)=0$ and using \eqref{eq:pivots} gives
\[
\begin{aligned}
 0={}&-2M\left(\frac{\lambda_s}{\lambda}+\beta\right)
       +\frac{\lambda_s}{\lambda}(\Dtwo\epsn,P_a^2)
       -\lambda(\diver\epsv,P_a^2)\\
    &+2a_s(\epsn,P_a\partial_aP_a).
\end{aligned}
\]
Collecting the scale defect proves \eqref{eq:density-eq}.
All wave derivatives are transposed onto the smooth tests:
\[
 (\diver\epsv,P_a^2)=-(\epsv,\nabla P_a^2),\qquad
 (\Dtwo\epsn,P_a^2)=-(\epsn,y\cdot\nabla P_a^2).
\]

Equations \eqref{eq:amplitude-eq} and \eqref{eq:density-eq}
form a real system for $a_s$ and $\lambda_s/\lambda+\beta$.
At zero error its coefficient matrix is
\begin{equation}\label{eq:modulation-matrix}
 \begin{pmatrix}
 \displaystyle\Rea\frac{\pair{\partial_a Q_{a,\beta}+(\beta/a)\partial_\beta Q_{a,\beta}}{\phi(N_a)}}
                              {\pair{Q_{a,\beta}}{\phi(N_a)}}&0\\[2mm]
 0&-2M
 \end{pmatrix}
 =\begin{pmatrix}M'/(2M)&0\\0&-2M\end{pmatrix}+O_I(\beta^2).
\end{equation}
The zero off-diagonal entry is the same scale cancellation used
in the parameter Jacobian. The first phase corrections are purely
imaginary, giving the stated real $O_I(\beta^2)$ error.
For small density error the matrix changes by $O_I(\|\epsn\|_2)$,
so its inverse is bounded uniformly on $I$.

To estimate the forcing, use the localized eigenfunction bounds:
\[
\begin{aligned}
 |\pair{Q_{a,\beta}}\phi|&\ge c_I>0,\\
 |\pair{\Lambda Q_{a,\beta}-\beta\partial_\beta Q_{a,\beta}}\phi|
     &\le C_I\|\epsn\|_2,\\
 |\pair{L_aQ_{a,\beta}}\phi|
     &=|\pair{L_aQ_{a,\beta}}{\phi-\phi(N_a)}|
                            \le C_I\beta\|\epsn\|_2,\\
 |\pair\eps{\phi'_N[\diver V]}|
     &\le C_I\|\veps\|_Y(1+\|\veps\|_Y).
\end{aligned}
\]
The last bound follows by splitting $V=V_a+\epsv$ and using the
one-divergence estimate. Both forcing rows are therefore bounded
by $C_I\lambda\|\veps\|_Y$, and inversion proves
\eqref{eq:rate-bound}. In particular,
\begin{equation}\label{eq:beta-rate}
 |\beta_s+\beta^2|\le C_I\beta^2\|\veps\|_Y .
\end{equation}
This follows from the second line of
\eqref{eq:profile-time-derivative}, since $a$ stays in $I$.
The tested coefficients and forcing are continuous in the parameters and the $Y$ error, and the inverse
is uniformly bounded on the bootstrap regime. Thus the rates
have the continuous dependence needed for the approximation
argument in Lemma~\ref{lem:energy-solutions}.

\paragraph{The phase estimate.}
Return to \eqref{eq:orthogonality-evolution}, divide by
$\pair{Q_{a,\beta}}\phi$, and take imaginary parts. The two rates already
estimated and the preceding forcing bounds give
\eqref{eq:phase-equation}. The expansion
\[
 \mu(N_a+\epsn)=-1+\frac{(\epsn,P_a^2)}M
                                   +O_I(\|\epsn\|_2^2)
\]
and the density condition prove \eqref{eq:phase-bound}.
The same tested pairings and matrix inversion give continuous
dependence of the phase rate. These canceled expressions are
continuous at energy regularity, so the integrated approximation
argument applies.
\end{proof}

The factor $\lambda$ in \eqref{eq:rate-bound} will allow parameter
convergence even for a small positive decay exponent. The quadratic
term in \eqref{eq:phase-bound} need not have this factor; its role in
the energy calculation is explained in Remark~\ref{rem:phase-work}.

\subsection{Bootstrap assumptions and improvement}\label{sec:bootstrap-setup}

With the fixed endpoints $0<a_-<a_0<a_+<a_*$ chosen in
Section~\ref{sec:notation}, set
\[
 I=[a_-,a_+],\qquad
 \delta_a=\tfrac14\min\{a_0-a_-,a_+-a_0\}>0.
\]
The small bound $a_*$ is chosen first; the interval $I$ and its
interior point $a_0$ are then fixed. The initial energy margin will be chosen relative to this
fixed distance from the boundary.
The reference $a_0$ is fixed; $a(s_0)$ is selected by the initial
modulation. Choose common thresholds $0<\beta_*,\delta<1$ for the
radial and full-space estimates. These thresholds control the
rescaled errors and are fixed before the radius of the physical
initial neighborhood.

Let $[0,T_{\max})$ be the maximal forward interval of the physical
energy solution, and let $t(s_0)=0$. On any interval where the
modulation of Lemma~\ref{lem:chart} is defined and $\lambda>0$, use
$\mathrm ds/\mathrm dt=\lambda^{-2}$. The bootstrap assumptions are
\begin{equation}\label{eq:bootstrap-assumptions}
 \begin{gathered}
 a(s)\in(a_-,a_+),\qquad
 0<\beta(s)=c_0a(s)\lambda(s)<\beta_*,\qquad \lambda(s)>0,\\
 \norm{\veps(s)}_Y<\delta,\qquad \norm{\eps(s)}_2<\delta,
 \end{gathered}
\end{equation}
together with the existence of a continuous modulation satisfying
\eqref{eq:chart}.
For nonradial data the center is also chosen by
$(\epsn,\partial_jN_a)=0$, $j=1,2$; the same assumptions use the
norms without symmetry restrictions. The center ranges over $\R^2$ and the phase is a continuous real
lift. The bootstrap imposes only the displayed smallness conditions;
parameter monotonicity and convergence will be conclusions.

For initialization we require the strict bounds
\begin{equation}\label{eq:bootstrap-initialization}
 \begin{aligned}
 |a(s_0)-a_0|&<\delta_a,&
 c_0a_+\lambda(s_0)&<\beta_*/2,\\
 \norm{\veps(s_0)}_Y&<\delta/2,&
 \norm{\eps(s_0)}_2&<\delta/2.
 \end{aligned}
\end{equation}
First choose the reference scale sufficiently small and then shrink
its open physical energy neighborhood. Continuity of the initial
modulation gives these inequalities. Local existence and continuity
therefore give a nonempty interval on which
\eqref{eq:bootstrap-assumptions} holds.

Let $[s_0,s_*)$ be the maximal continuous interval starting at
$s_0$ on which the solution, its modulation and
\eqref{eq:bootstrap-assumptions} all exist. Its endpoint is the
supremum of the admissible right endpoints; it may be finite at
this stage.

\begin{proposition}[Bootstrap improvement]\label{prop:bootstrap}
There is a sufficiently
small branch bound $a_*>0$ such that, after fixing
$0<a_-<a_0<a_+<a_*$ and $\delta_a$ as above, thresholds
$\beta_*=\beta_*(I)>0$ and $\delta=\delta(I)>0$ can be chosen
with the following property. For every sufficiently late point
of the exact solution with parameter $a_0$, there is an open
physical energy neighborhood of that datum such that every
solution starting in this neighborhood satisfies, on $[s_0,s_*)$,
\begin{equation}\label{eq:bootstrap-improvement}
 \begin{aligned}
 |a(s)-a_0|&<2\delta_a,&0<\beta(s)&<\beta_*/2,\\
 \norm{\veps(s)}_Y&<\delta/2,&
 \norm{\eps(s)}_2&<\delta/\sqrt2.
 \end{aligned}
\end{equation}
Moreover $\lambda_s<0$ and $s_*=\infty$. The statement holds for
radial perturbations and, with translation modulation, for
nonradial perturbations.
\end{proposition}

The energy estimates below use only
\eqref{eq:bootstrap-assumptions}. Their proof does not assume
convergence or a decay rate. The stricter initial margins needed
to improve these assumptions are chosen in
Section~\ref{sec:closure}; the proof for general perturbations is completed in
Section~\ref{sec:nonradial}.

\section{Linear energy estimates along a self-similar solution}
\label{sec:linear-energy}

We prove a conditional linear energy estimate for radial perturbations
of a fixed self-similar solution, with $a\in I$ and its exact scale and phase.
The proof first uses the interaction energy to cancel the exchange,
then subtracts the density response before removing field mass.
A quadratic correction removes the remaining mixed terms.
The two orthogonality conditions below are imposed throughout
the interval of the estimate. The stability proof in Section~\ref{sec:energy}
derives the energy estimate directly along the nonlinear modulation.
It retains the full density in the Schr\"odinger operator and uses
the parameter equations to estimate the remainder.

\subsection{Linearized evolution and orthogonality conditions}\label{app:linear}

The reference parameters satisfy
\[
 a_s=0,\qquad \lambda_s=-c_0a\lambda^2,\qquad
 \gamma_s=1,\qquad \beta=c_0a\lambda,\qquad \beta_s=-\beta^2.
\]
The profile identities are
\[
 \begin{aligned}
 &L_aQ_{a,\beta}
 =\e^{-\ii\beta r^2/4}
       \left(\ii\beta\Lambda P_a+\frac{\beta^2r^2}{4}P_a\right),\\
 &(Q_{a,\beta})_s
 =-\beta^2\partial_\beta Q_{a,\beta}
   =-\ii L_aQ_{a,\beta}-\beta\Lambda Q_{a,\beta},\\
 &-\beta\Dtwo N_a-\lambda\diver V_a=0,\\
 &-\beta\Dtwo V_a-c_0^2\lambda\nabla(N_a+P_a^2)=0 .
 \end{aligned}
\]
Write the rescaled solution as
\[
 (\Psi,N,V)=(Q_{a,\beta}+\eps,N_a+\epsn,V_a+\epsv).
\]
Subtraction of the reference equations from \eqref{eq:rescaled}
gives the exact perturbation equation on this fixed parameter path:
\[
 \veps_s=\mathcal A_\beta\veps+
       ( -\ii\epsn\eps,\ 0,\ -c_0^2\lambda\nabla|\eps|^2).
\]
Dropping the quadratic products defines the real linear operator
\begin{equation}\label{eq:G}
 \mathcal A_\beta\veps=
 \begin{pmatrix}
 -\ii(L_a\eps+Q_{a,\beta}\epsn)-\beta\Lambda\eps\\
 -\beta\Dtwo\epsn-\lambda\diver\epsv\\
 -\beta\Dtwo\epsv-c_0^2\lambda\nabla
                    (\epsn+2\Rea(\overline{Q_{a,\beta}}\eps))
 \end{pmatrix},
 \qquad \veps_s=\mathcal A_\beta\veps .
\end{equation}
In the real components, \eqref{eq:G} has the block matrix form
\[
\partial_s
\begin{pmatrix}
 \Rea\eps\\ \Ima\eps\\ \epsn\\ \epsv
\end{pmatrix}
=
\begin{pmatrix}
 -\beta\Lambda & L_a & \Ima Q_{a,\beta} & 0\\
 -L_a & -\beta\Lambda & -\Rea Q_{a,\beta} & 0\\
 0 & 0 & -\beta\Dtwo & -\lambda\diver\\
 -2c_0^2\lambda\nabla[(\Rea Q_{a,\beta})\,\cdot]
   & -2c_0^2\lambda\nabla[(\Ima Q_{a,\beta})\,\cdot]
   & -c_0^2\lambda\nabla & -\beta\Dtwo
\end{pmatrix}
\begin{pmatrix}
 \Rea\eps\\ \Ima\eps\\ \epsn\\ \epsv
\end{pmatrix}.
\]
Profile entries act by multiplication. In the last row the
gradient acts on the product with the input scalar. The last
component is a vector, and $\Dtwo$ acts on it componentwise.

The tangent conditions to \eqref{eq:chart} at the reference state are
\begin{equation}\label{eq:old-slice}
 \pair\eps{P_a}=0,\qquad(\epsn,P_a^2)=0.
\end{equation}
Here the first equality is one complex condition, hence two real
conditions. These are ordinary $L^2$ pairings with the displayed
profiles. Along the fixed-orbit linear equation,
\[
 \begin{aligned}
 \frac{\mathrm d}{\mathrm ds}\pair\eps{P_a}
          &=-\ii\pair{Q_{a,\beta}\epsn}{P_a}
                          +\beta\pair\eps{\Lambda P_a},\\
 \frac{\mathrm d}{\mathrm ds}(\epsn,P_a^2)
          &=\beta(\epsn,y\cdot\nabla P_a^2)
                       +\lambda(\epsv,\nabla P_a^2).
 \end{aligned}
\]
These identities show that the unmodulated fixed-orbit evolution
does not automatically preserve the tangent conditions. They are
imposed throughout the interval in the following conditional estimate.
The nonlinear argument derives the parameter work directly from
\eqref{eq:chart}.

\begin{proposition}[Linear energy estimate under orthogonality]\label{prop:linear}
Choose $a_*>0$ sufficiently small. For every $I\Subset(0,a_*)$ there are
$\lambda_*,c_I,C_I,d_I>0$ and a real quadratic energy
$\mathcal E_\beta$ on $Y$ such that
\[
 c_I\|\veps\|_Y^2\le\mathcal E_\beta(\veps)
                   \le C_I\|\veps\|_Y^2
\]
for radial $\veps\in Y$ satisfying \eqref{eq:old-slice},
$a\in I$, and $0<\lambda<\lambda_*$.
Let $\veps(s)$ be a smooth radial finite-energy solution of
\eqref{eq:G}, with $a$ fixed and
$\lambda_s=-c_0a\lambda^2$. On every interval where
\eqref{eq:old-slice} holds, its energy satisfies
\begin{equation}\label{eq:linear-dissipation}
 \frac{\mathrm d}{\mathrm ds}\mathcal E_\beta(\veps)
       \le-d_I\beta\|\veps\|_Y^2+C_I\beta^2\|\veps\|_Y^2 .
\end{equation}
\end{proposition}

\begin{remark}[Energy regularity]\label{rem:linear-regularity}
The quadratic energy and the canceled expression for its
derivative extend continuously to the slice \eqref{eq:old-slice}
in $Y$ by smooth approximation of the triples. For the nonlinear
flow, Lemma~\ref{lem:energy-solutions} proves the integrated
identities by approximation of solutions on compact intervals.
\end{remark}

\subsection{Interaction energy and mass subtraction}
\label{sec:linear-response}

There are two separate tasks. The interaction energy cancels the
field--wave exchange, but contains global field mass. We remove
that mass only after subtracting the density response from the
field equation. The resulting derivative is slow; its mixed terms
will determine the final correction.

\paragraph{The interaction identity.}
We augment the quadratic perturbation energy by the
Schr\"odinger mass. Its interaction term cancels the coupled fast pairings, as the
product rule below shows. The form is
\[
 E_\beta^{(2)}=(L_a{\eps},{\eps})_\R
 +2({\epsn},\Rea(\overline{Q_{a,\beta}} {\eps}))
 +\tfrac12\norm{\epsn}_2^2+\tfrac1{2c_0^2}\norm {\epsv}_2^2.
\]
The $L_a$ term contains $\|\eps\|_2^2$, which is not part
of the $Y$ norm. We first differentiate this energy to retain
its exact cancellations; the corrected mass subtraction below removes the global mass.
Set $A_a=-2\Delta-y\cdot\nabla N_a$, the scalar dilation
commutator $[L_a,\Lambda]$. From \eqref{eq:G},
\begin{align*}
 (E_\beta^{(2)})_s
 ={}&2\Rea\pair{\eps_s}{L_a\eps+Q_{a,\beta}\epsn}
       +((\epsn)_s,\epsn+2\Rea(\overline{Q_{a,\beta}}\eps))\\
 &+c_0^{-2}((\epsv)_s,\epsv)
       +2(\epsn,\Rea((\overline{Q_{a,\beta}})_s\eps))\\
 ={}&\underbrace{2\Rea\pair{-\ii(L_a\eps+Q_{a,\beta}\epsn)}
                               {L_a\eps+Q_{a,\beta}\epsn}}_{=0}\\
 &-\lambda\underbrace{\left[
   (\diver\epsv,\epsn+2\Rea(\overline{Q_{a,\beta}}\eps))
  +(\epsv,\nabla(\epsn+2\Rea(\overline{Q_{a,\beta}}\eps)))
                              \right]}_{=0}\\
 &-2\beta\Rea\pair{\Lambda\eps}{L_a\eps}
       -\beta\bigl[(\Dtwo\epsn,\epsn)
                         +c_0^{-2}(\Dtwo\epsv,\epsv)\bigr]\\
 &-2\beta\bigl[(\epsn,\Rea(\overline{Q_{a,\beta}}\Lambda\eps))
             +(\Dtwo\epsn,\Rea(\overline{Q_{a,\beta}}\eps))\bigr]\\
 &+2(\epsn,\Rea((\overline{Q_{a,\beta}})_s\eps)).
\end{align*}
Integration by parts and $\beta_s=-\beta^2$ give
\[
 \begin{aligned}
 2\Rea\pair{\Lambda\eps}{L_a\eps}&=(A_a\eps,\eps)_\R,\\
 (\Dtwo\epsn,\epsn)+c_0^{-2}(\Dtwo\epsv,\epsv)
                 &=\norm{\epsn}_2^2+c_0^{-2}\norm{\epsv}_2^2,\\
 (\epsn,\Rea(\overline{Q_{a,\beta}}\Lambda\eps))
 +(\Dtwo\epsn,\Rea(\overline{Q_{a,\beta}}\eps))
                 &=(\epsn,\Rea([(2-\Lambda)\overline{Q_{a,\beta}}]\eps)),\\
 (\overline{Q_{a,\beta}})_s
                 &=-\frac{\ii\beta^2r^2}{4}\overline{Q_{a,\beta}}.
 \end{aligned}
\]
This yields
\begin{align}
 (E_\beta^{(2)})_s
 &=-\beta D_{a,\beta}({\veps})
   +\frac{\beta^2}{2}({\epsn},r^2\Ima(\overline{Q_{a,\beta}} {\eps})),
 \label{eq:global-linear-identity}\\
 D_{a,\beta}({\veps})&=(A_a{\eps},{\eps})_\R
 +2({\epsn},\Rea([(2-\Lambda)\overline{Q_{a,\beta}}]{\eps}))
 +\norm{\epsn}_2^2+c_0^{-2}\norm {\epsv}_2^2.\nonumber
\end{align}
The form $D_{a,\beta}$ acts on the full triple. In $D_{a,0}$ the
phase is set to zero at fixed $a$, so $Q_{a,\beta}$ becomes $P_a$.
The mixed term has no sign and must be canceled before proving
dissipation.
Moreover,
\[
 \left|\frac{\beta^2}{2}
       (\epsn,r^2\Ima(\overline{Q_{a,\beta}}\eps))\right|
 \le C_I\beta^2\norm{\epsn}_2\norm\eps_{\mathrm{loc}}
 \le C_I\beta^2\norm\veps_Y^2.
\]

\paragraph{Subtraction of the corrected mass.}
To remove the global Schr\"odinger mass one might subtract
$\norm\eps_2^2$, but \eqref{eq:G} gives
\[
 \begin{aligned}
 (\norm\eps_2^2)_s
 &=2\Rea\pair{-\ii(L_a\eps+Q_{a,\beta}\epsn)-\beta\Lambda\eps}{\eps}\\
 &=2\Ima\pair{Q_{a,\beta}\epsn}{\eps}.
 \end{aligned}
\]
This reintroduces an order-one Schr\"odinger--density term.
We first remove the density forcing from the Schr\"odinger
equation, then subtract the mass of the corrected error. Recall from \eqref{eq:RK}
that
\[
 K_\beta f=-R_a(Q_{a,\beta}f),\qquad
 L_aK_\beta f=-\Pi_a(Q_{a,\beta}f).
\]
For $z=\eps-K_\beta\epsn$,
\[
 L_a\eps+Q_{a,\beta}\epsn
 =L_az+(1-\Pi_a)(Q_{a,\beta}\epsn)
 =L_az+\frac{\pair{Q_{a,\beta}\epsn}{P_a}}M P_a,
\]
and
\[
 \begin{aligned}
 z_s={}&\eps_s-(K_\beta)_s\epsn-K_\beta(\epsn)_s\\
 ={}&-\ii L_az-\frac{\ii}{M}\pair{Q_{a,\beta}\epsn}{P_a}P_a
          -\beta\Lambda(z+K_\beta\epsn)\\
 &+K_\beta(\beta\Dtwo\epsn+\lambda\diver\epsv)
          -(K_\beta)_s\epsn.
 \end{aligned}
\]
Thus
\begin{equation}\label{eq:good-unknown-linear}
\begin{aligned}
 z_s={}&-\ii L_az-\frac{\ii}{M}\pair{Q_{a,\beta}\epsn}{P_a}P_a
              -\beta\Lambda z\\
 &+\beta(K_\beta\Dtwo-\Lambda K_\beta)\epsn
       +\lambda K_\beta\diver\epsv-(K_\beta)_s\epsn.
\end{aligned}
\end{equation}
The fast density source is now only a ground-state term. Its
contribution to the mass derivative vanishes under the first
condition in \eqref{eq:old-slice}; all other terms have wave or
dilation coefficients.
The two dilation terms in $K_\beta\Dtwo-\Lambda K_\beta$
are combined before applying the operator to $L^2$ density data.
In the leading derivative we set the quadratic phase to zero in
the coefficients of the source and keep $z$ in the pairings.
Set $\chi_a=R_a\Lambda P_a$. The explicit formula follows from
$L_a(r^2P_a)=-4\Lambda P_a$ by projecting off $P_a$:
\begin{equation}\label{eq:chi}
 \chi_a=R_a\Lambda P_a=-\tfrac14(r^2-\frac{\norm{yP_a}_2^2}{M})P_a.
\end{equation}
The commutator and the product rule give
\begin{align}
 [L_a,\Lambda]&=A_a,\nonumber\\
 [\Lambda,R_a]&=R_aA_aR_a
       -\frac{\chi_a\otimes P_a+P_a\otimes\chi_a}{M},
       \label{eq:R-commutator}\\
 \Lambda R_a(P_af)-R_a(P_a\Dtwo f)
      &=[\Lambda,R_a](P_af)-R_a([(2-\Lambda)P_a]f).
       \label{eq:K-commutator}
\end{align}
The tensors are rank-one real duality operators, extended componentwise
to complex functions.
In the last line, $P_a\Dtwo f-\Lambda(P_af)=[(2-\Lambda)P_a]f$.
The right-hand side of \eqref{eq:K-commutator} extends to
$L^2$ densities. This is the density term used in the quadratic
correction below. Appendix~\ref{app:linear-details} gives the
projection calculation and the bounds for the phase-dependent terms.
On the tangent constraint space, the kernel pairing vanishes and
the remaining density coefficient is small:
\[
 \begin{aligned}
 \pair z{P_a}&=\pair\eps{P_a}=0,\\
 |\pair{Q_{a,\beta}\epsn}{P_a}|
 &=|\pair{(Q_{a,\beta}-P_a)\epsn}{P_a}|
                                  \le C_I\beta\norm{\epsn}_2.
 \end{aligned}
\]
Here $K_\beta\epsn\perp P_a$ and $(\epsn,P_a^2)=0$.
Equation~\eqref{eq:good-unknown-linear} yields
\[
 \begin{aligned}
 (\norm z_2^2)_s
 ={}&2\Rea\pair{z_s}z\\
 ={}&\frac2M\Ima\left\{
       \pair{Q_{a,\beta}\epsn}{P_a}\overline{\pair z{P_a}}\right\}
       +2\beta\Rea\pair{(K_\beta\Dtwo-\Lambda K_\beta)\epsn}z\\
 &+2\lambda\Rea\pair{K_\beta\diver\epsv}z
       -2\Rea\pair{(K_\beta)_s\epsn}z.
 \end{aligned}
\]
The commutator and phase bounds in
Appendix~\ref{app:linear-details} imply
\[
 \begin{aligned}
 &\frac2M\left|\pair{Q_{a,\beta}\epsn}{P_a}\pair z{P_a}\right|
      +2\left|\pair{(K_\beta)_s\epsn}z\right|
                                      \le C_I\beta^2\norm\veps_Y^2,\\
 &\norm{z-\eps-R_a(P_a\epsn)}_{H^2_d}\le C_I\beta\norm{\epsn}_2,\\
 &|D_{a,\beta}(\veps)-D_{a,0}(\veps)|\le C_I\beta\norm\veps_Y^2.
 \end{aligned}
\]
Consequently,
\[
 \begin{aligned}
 (\norm z_2^2)_s
 ={}&2\beta\Rea\pair{\Lambda R_a(P_a\epsn)-R_a(P_a\Dtwo\epsn)}{z}\\
 &-2\lambda\Rea\pair{R_a(P_a\diver\epsv)}{z}
       +O_I(\beta^2\norm\veps_Y^2).
 \end{aligned}
\]
The dilation difference is evaluated by \eqref{eq:K-commutator}.
For the velocity term, the product of $Q_{a,\beta}-P_a$ with
$\diver\epsv$ has $H^{-1}_{d'}$ norm at most $C_I\beta\|\epsv\|_2$
by the distributional product rule. Applying $R_a$ gives the
corresponding $H^1_d$ bound.

Subtracting this identity from \eqref{eq:global-linear-identity} gives
\begin{equation}\label{eq:mass-subtracted-identity}
\begin{aligned}
 (E_\beta^{(2)}-\norm z_2^2)_s
 ={}&-\beta D_{a,0}(\veps)
       -2\beta\Rea\pair{\Lambda R_a(P_a\epsn)-R_a(P_a\Dtwo\epsn)}{z}\\
 &+2\lambda\Rea\pair{R_a(P_a\diver\epsv)}{z}
       +O_I(\beta^2\norm\veps_Y^2).
\end{aligned}
\end{equation}
The same subtraction removes the global Schr\"odinger mass
from the energy itself:
\[
 \begin{aligned}
 (L_a\eps,\eps)_\R-\norm\eps_2^2
       &=\norm{\nabla\eps}_2^2+\int N_a|\eps|^2,\\
 \norm\eps_2^2-\norm{\eps-K_\beta\epsn}_2^2
       &=2\Rea\pair\eps{K_\beta\epsn}-\norm{K_\beta\epsn}_2^2.
 \end{aligned}
\]
All surviving Schr\"odinger factors are gradients or are paired
with localized coefficients, so the expanded expression defines
a form on $Y$. Its derivative still has the two mixed terms in
\eqref{eq:mass-subtracted-identity}; we cancel them next.

\subsection{Quadratic correction of the mixed terms}
\label{sec:linear-correction}

To separate the diagonal and mixed terms, expand the full dilation
form $D_{a,0}$ using $\eps=z-R_a(P_a\epsn)+O_{H^2_d}(\beta\|\epsn\|_2)$.
Write $[(2-\Lambda)P_a]$ for multiplication by $P_a-y\cdot\nabla P_a$.
The real linear coefficient of the mixed density term is
\begin{equation}\label{eq:Fa}
 \mathcal M_a=-A_aR_aP_a+[(2-\Lambda)P_a]
                    +[\Lambda,R_a]P_a-R_a[(2-\Lambda)P_a].
\end{equation}
The last two summands are the dilation difference in
\eqref{eq:mass-subtracted-identity}, by \eqref{eq:K-commutator}.
The other two come from the cross terms in $D_{a,0}$. This yields
\begin{equation}\label{eq:diagonal-expansion}
\begin{aligned}
 &D_{a,0}(\veps)+2\Rea\pair{\Lambda R_a(P_a\epsn)-R_a(P_a\Dtwo\epsn)}{z}\\
 &\qquad=(A_az,z)_\R+D_a(\epsn)+c_0^{-2}\|\epsv\|_2^2
             +2(\Rea z,\mathcal M_a\epsn)+O_I(\beta\|\veps\|_Y^2),
\end{aligned}
\end{equation}
where $D_a(f)$ acts only on the density
\[
\begin{aligned}
 D_a(f)={}&\|f\|_2^2+(A_aR_a(P_af),R_a(P_af))\\
          &-2(f,[(2-\Lambda)P_a]R_a(P_af)).
\end{aligned}
\]
The error in \eqref{eq:diagonal-expansion} follows from boundedness of
$D_{a,0}$ on $Y$ and the localized bound on the response difference.
The derivative therefore has the mixed term
\[
 -2\beta\bigl(\Rea z,\mathcal M_a\epsn
                     -(c_0a)^{-1}R_a(P_a\diver\epsv)\bigr).
\]
Choose
\[
 p=R_a\mathcal M_a\epsn-(c_0a)^{-1}R_a^2(P_a\diver\epsv),\qquad
 (p,P_a)=0,
\]
and use the fast field equation to differentiate the pairing:
\[
 (\Ima z,p)_s=-(\Rea z,L_ap)
     +(\Ima(z_s+\ii L_az),p)+(\Ima z,p_s).
\]
Thus the negative of $2\beta$ times this pairing differentiates
to the mixed term with the opposite sign. We define the energy by
\[
 \mathcal E_\beta=E_\beta^{(2)}-\norm z_2^2
                                      -2\beta(\Ima z,p).
\]
Expanding the mass difference,
\begin{equation}\label{eq:F}
\begin{aligned}
 \mathcal E_\beta({\veps})={}&\norm{\nabla {\eps}}_2^2+\int N_a\abs {\eps}^2
 +2({\epsn},\Rea(\overline{Q_{a,\beta}}{\eps}))
 +\tfrac12\norm{\epsn}_2^2+\tfrac1{2c_0^2}\norm {\epsv}_2^2\\
 &+2\Rea\pair {\eps}{K_\beta {\epsn}}-\norm{K_\beta {\epsn}}_2^2\\
 &-2\beta\bigl(\Ima z,
       R_a\mathcal M_a\epsn-(c_0a)^{-1}R_a^2(P_a\diver\epsv)\bigr).
\end{aligned}
\end{equation}
The phase enters through $Q_{a,\beta}$, the response $K_\beta$
and the explicit factor in the correction. Along the fixed orbit
these three occurrences vary with $\beta(s)$.
The first two lines are the expanded mass difference from
the preceding subtraction. The last line is a bounded quadratic correction on $Y$.
We now check its derivative, retaining the projection in
$L_aR_a=\Pi_a$.

For the correction term,
\[
 \begin{aligned}
 (-2\beta(\Ima z,p))_s
 ={}&-2\beta_s(\Ima z,p)-2\beta((\Ima z)_s,p)
                                      -2\beta(\Ima z,p_s)\\
 ={}&2\beta(\Rea z,L_ap)-2\beta_s(\Ima z,p)\\
 &-2\beta(\Ima(z_s+\ii L_az),p)-2\beta(\Ima z,p_s).
 \end{aligned}
\]
Since $L_ap=\Pi_a[\mathcal M_a\epsn-(c_0a)^{-1}R_a(P_a\diver\epsv)]$
and $z\perp P_a$ on \eqref{eq:old-slice}, the leading term is exactly
\[
 2\beta(\Rea z,L_ap)
 =2\beta\bigl(\Rea z,\mathcal M_a\epsn
                  -(c_0a)^{-1}R_a(P_a\diver\epsv)\bigr).
\]
The remaining terms involve $\beta_s$, $p_s$, or the slow part
of the corrected equation. Appendix~\ref{app:linear-details}
bounds each of them by $C_I\beta^2\|\veps\|_Y^2$.
Consequently the fixed-orbit derivative identity is
\begin{equation}\label{eq:unprojected-diagonal}
 (\mathcal E_\beta)_s=-\beta[(A_az,z)_\R+D_a({\epsn})
                         +c_0^{-2}\norm {\epsv}_2^2]
                         +O_I(\beta^2\norm{\veps}_Y^2).
\end{equation}
The mixed terms have canceled; only the diagonal forms remain.

\subsection{Coercivity and linear dissipation}
\label{sec:linear-positivity}

\begin{proof}[Proof of Proposition~\ref{prop:linear}]
The preceding construction gives the energy \eqref{eq:F}
and its derivative \eqref{eq:unprojected-diagonal}.
We prove the signs of their principal forms separately.
The perturbation estimates in Appendix~\ref{app:energy-coercivity}
transfer the estimates of Proposition~\ref{ass:spectral} to the small positive branch.
The scalar dissipation form controls $\|z\|_{\mathrm{loc}}^2$
because $z\perp P_a$. The density form $D_a(\epsn)$ controls
$\|\epsn\|_2^2$ on $(\epsn,P_a^2)=0$. Together with the velocity
square and the bounded response $\eps-z=K_\beta\epsn$, these
bounds give \eqref{eq:linear-dissipation}.

At zero phase, the exact completed square is
\begin{equation}\label{eq:completed-square}
\begin{aligned}
 \mathcal E_0(\veps)
 ={}&\pair{(L_a-1)(\eps+R_a(P_a\epsn))}{\eps+R_a(P_a\epsn)}_\R\\
    &+\tfrac12(\epsn,(1-2P_aR_aP_a)\epsn)
      +\tfrac1{2c_0^2}\|\epsv\|_2^2,
 \qquad (\epsn,P_a^2)=0 .
\end{aligned}
\end{equation}
Here $L_aR_a(P_a\epsn)=P_a\epsn$ on the density constraint.
The phase-dependent terms in \eqref{eq:F} change this form
by at most $C_I\beta\|\veps\|_Y^2$. Replacing
$\eps+R_a(P_a\epsn)$ in its scalar square by $z$ has the same cost
by the localized phase estimate and boundedness of the scalar form
on $X$. Apply the Schr\"odinger energy gap to $z\perp P_a$ and the density
energy gap to $(\epsn,P_a^2)=0$ to obtain
\[
 \mathcal E_\beta(\veps)
       \ge c_I\|\veps\|_Y^2-C_I\beta\|\veps\|_Y^2.
\]
Taking $\lambda_*(I)$ sufficiently small proves the positive
lower bound. The upper bound follows from the bounded quadratic
forms. Approximation in the same norms proves the extension
stated in Remark~\ref{rem:linear-regularity}.
\end{proof}

\section{Nonlinear energy estimates and dissipation}\label{sec:energy}

We now carry out the construction for the modulated nonlinear
solution. The operator $L_N$ retains the nonlinear density coupling,
and the full cubic interaction is included before differentiation.
These choices allow exact cancellation of the terms with no
wave coefficient. We then subtract the mass of the corrected error
and construct the quadratic correction, keeping the parameter work
explicit. The purpose of each step is to obtain a positive energy
whose nonlinear remainder retains $\lambda$.

All estimates hold on the bootstrap interval of
Section~\ref{sec:bootstrap-setup}, with $a\in I$ and
$\beta=c_0a\lambda$. In particular,
\[
 |a_s|+\left|\frac{\lambda_s}{\lambda}+\beta\right|
       \le C_I\lambda\|\veps\|_Y.
\]

\begin{proposition}[Energy and dissipation estimates]\label{prop:remainder}
Let $a$ lie in a fixed
sufficiently small positive compact interval $I$, and let the perturbation satisfy \eqref{eq:chart}. For small
$\beta$, $\norm{\veps}_Y$ and $\norm{\eps}_2$, there is a real energy $\E$, constructed below, such that
\begin{equation}\label{eq:core-coercivity}
 c_I\norm{\veps}_Y^2\le\E\le C_I\norm{\veps}_Y^2.
\end{equation}
Along the modulated solution, for constants $d_I,C_I>0$ depending on $I$,
\begin{equation}\label{eq:energy-virial}
 \E_s\le-d_I \beta\norm{\veps}_Y^2
 +C_I \beta(\beta+\norm{\veps}_Y+\norm {\eps}_2)\norm{\veps}_Y^2.
\end{equation}
Consequently, after fixing the smallness thresholds,
\begin{equation}\label{eq:closed-energy}
 \E_s\le-\kappa_I\lambda\E,\qquad
 \abs{a_s}+\abs{\frac{\lambda_s}{\lambda}+\beta}\le K_I\lambda\sqrt\E,
 \qquad \kappa_I>0.
\end{equation}
\end{proposition}

\subsection{Density response and the nonlinear energy}
\label{sec:nonlinear-cancellation}

All equations in this section concern the actual modulated solution
of \eqref{eq:rescaled}, subject to \eqref{eq:chart}. Write
$\phi=\phi(N)$. The parameters $a$ and $\lambda$ vary
with $s$, and $\beta=c_0a\lambda$. We first calculate for smooth
solutions and then pass to energy solutions.

\paragraph{The modulated error equations.}
Subtracting the current profiles, including their parameter
derivatives, gives
\begin{equation}\label{eq:actual-nonlinear-error}
\begin{aligned}
 \ii\eps_s={}&(L_N+\gamma_s-1)\eps+Q_{a,\beta}\epsn
       +\ii\frac{\lambda_s}{\lambda}\Lambda\eps+(\gamma_s-1)Q_{a,\beta}\\
 &+\ii\left(\frac{\lambda_s}{\lambda}+\beta\right)\Lambda Q_{a,\beta}
       -\ii a_s\partial_a Q_{a,\beta}
       -\ii(\beta_s+\beta^2)\partial_\beta Q_{a,\beta},\\
 (\epsn)_s={}&\frac{\lambda_s}{\lambda}\Dtwo\epsn-\lambda\diver\epsv
       +\left(\frac{\lambda_s}{\lambda}+\beta\right)\Dtwo N_a
       -a_s\partial_aN_a,\\
 (\epsv)_s={}&\frac{\lambda_s}{\lambda}\Dtwo\epsv
       -c_0^2\lambda\nabla(\epsn+2\Rea(\overline{Q_{a,\beta}}\eps)+|\eps|^2)\\
 &+\left(\frac{\lambda_s}{\lambda}+\beta\right)\Dtwo V_a
       -a_s\partial_aV_a,
 \qquad L_N=-\Delta+1+N=L_a+\epsn .
\end{aligned}
\end{equation}
The first row uses the chain rule \eqref{eq:profile-time-derivative}
and the reference identity
\[
 -\ii L_a Q_{a,\beta}-\beta\Lambda Q_{a,\beta}
       =-\beta^2\partial_\beta Q_{a,\beta}.
\]

\paragraph{The reduced inverse.}
Let $\Pi_N$ be the orthogonal projection onto $\phi(N)^\perp$.
We use the inverse of $L_N$ on this complement:
\begin{equation}\label{eq:current-inverse}
 \Pi_Ng=g-\pair g\phi\phi,\qquad
 R_N=\left(L_N\big|_{\phi^\perp}\right)^{-1}\Pi_N,
 \qquad K_Nf=-R_N(Q_{a,\beta}f).
\end{equation}
The real density argument $f$ belongs to $L^2$.
In $K_N$, the subscript records the density used in the inverse.
The profile factor also depends on the current $(a,\beta)$.
The eigenvalue of $L_N$ along $\phi$ is $1+\mu(N)$ and is
excluded from the inverse. At the same $(a,\beta)$, setting $N=N_a$
gives $R_{N_a}=R_a$ and $K_{N_a}=K_\beta$.

The localized bounds used below are $R_N:L^2_{d'}\to H^2_d$
and $R_N:H^{-1}_{d'}\to H^1_d$, together with the corresponding
bounds for one velocity-divergence direction. Their precise
statement and proof are Lemma~\ref{lem:current-inverse} in
Appendix~\ref{app:current-inverse}.

\paragraph{The corrected error and the energy.}
Define
\begin{equation}\label{eq:nonlinear-good-error}
 z=\eps-K_N\epsn,\qquad K_N\epsn=-R_N(Q_{a,\beta}\epsn).
\end{equation}
The density response is localized, with
$\|K_N\epsn\|_{H^2_d}\le C_I\|\epsn\|_2$.

Use the same quadratic interaction energy $E_\beta^{(2)}$
as in Section~\ref{sec:linear-energy}, now evaluated at the
current $(a,\beta)$, and add the cubic term
\[
 H_3=\int\epsn|\eps|^2.
\]
Subtract the mass of $z$ to define
\begin{equation}\label{eq:nonlinear-mass-energy}
 E=E_\beta^{(2)}+H_3-\|z\|_2^2.
\end{equation}
Expanding the square defines it without a bare global
Schr\"odinger mass:
\begin{equation}\label{eq:relative-action}
 \begin{aligned}
 E={}&\|\nabla\eps\|_2^2+\int N_a|\eps|^2
       +2(\epsn,\Rea(\overline{Q_{a,\beta}}\eps))
       +\tfrac12\|\epsn\|_2^2+\tfrac1{2c_0^2}\|\epsv\|_2^2\\
 &+2\Rea\pair\eps{K_N\epsn}-\|K_N\epsn\|_2^2+H_3.
 \end{aligned}
\end{equation}
Every term starts at degree two in the error. All terms except
$H_3$ are smooth on $Y$. The global mass enters only through
\begin{equation}\label{eq:H3}
 |H_3|\le\|\epsn\|_2\|\eps\|_4^2
       \le C\|\eps\|_2\|\epsn\|_2\|\nabla\eps\|_2
       \le C\|\eps\|_2\|\veps\|_Y^2.
\end{equation}
For later estimates, integration of the inverse derivative
along $N_a+\theta\epsn$, $0\le\theta\le1$, gives
\[
 \|K_N\epsn-K_\beta\epsn\|_{H^2_d}
       \le C_I\|\epsn\|_2^2.
\]
Thus the quadratic part of $E-H_3$ is the mass-subtracted
energy of Section~\ref{sec:linear-energy}. To differentiate
$E$, we need the equation for $z$. We first establish the phase
cancellation that makes its source of order $\lambda$.

\paragraph{The phase coefficient and the ground-state source.}
The phase estimate \eqref{eq:phase-equation} gives
\begin{equation}\label{eq:combined-phase-bound}
 |\mu(N)+\gamma_s|\le C_I\lambda\|\veps\|_Y.
\end{equation}
Since $K_N\epsn\perp\phi$ and $\pair\eps\phi=0$, we also have
$\pair z\phi=0$. Consequently the ground-state source in the
equation below has exactly zero work in the derivative of
$\|z\|_2^2$.

\paragraph{Evolution of the corrected error.}
Substitution into \eqref{eq:actual-nonlinear-error} gives
\begin{equation}\label{eq:nonlinear-good-equation}
 \ii z_s=(L_N+\gamma_s-1)z
       +\ii\frac{\lambda_s}{\lambda}\Lambda z+\mathcal F,
\end{equation}
where
\begin{equation}\label{eq:good-source}
 \begin{aligned}
 \mathcal F={}&(\mu(N)+\gamma_s)\pair{Q_{a,\beta}}\phi\phi
        -\pair{L_aQ_{a,\beta}}\phi\phi+(\gamma_s-1)R_NL_aQ_{a,\beta}\\
 &+\ii\left(\frac{\lambda_s}{\lambda}\Lambda(K_N\epsn)-\partial_s(K_N\epsn)\right)
       +\ii\left(\frac{\lambda_s}{\lambda}+\beta\right)\Lambda Q_{a,\beta}\\
 &-\ii a_s\partial_a Q_{a,\beta}
       -\ii(\beta_s+\beta^2)\partial_\beta Q_{a,\beta}.
 \end{aligned}
\end{equation}
For the first row use the exact identities
\[
 \begin{aligned}
 Q_{a,\beta}+K_N\epsn&=\pair{Q_{a,\beta}}\phi\phi+R_NL_aQ_{a,\beta},\\
 (1-\Pi_N)(Q_{a,\beta}\epsn)
   &=\bigl[(1+\mu(N))\pair{Q_{a,\beta}}\phi-\pair{L_aQ_{a,\beta}}\phi\bigr]\phi.
 \end{aligned}
\]
\begin{lemma}[Source estimate for the corrected field]
\label{lem:slow-source}
For a smooth modulated solution on the bootstrap regime, the
source in \eqref{eq:good-source} satisfies
\begin{equation}\label{eq:good-source-bound}
 \|\mathcal F\|_{H^1_d}\le C_I\lambda\|\veps\|_Y.
\end{equation}
\end{lemma}

\begin{proof}
\label{app:nonlinear-source-details}
Use the actual error equations \eqref{eq:actual-nonlinear-error},
the inverse \eqref{eq:current-inverse} and the corrected field
\eqref{eq:nonlinear-good-error}, under the bootstrap assumptions.
\paragraph{The time derivative of the inverse.}
Dilation covariance gives a useful formula without derivatives
of the rough density. If $U_\sigma g(y)=\sigma g(\sigma y)$
and $N_\sigma(y)=\sigma^2N(\sigma y)$, then
\[
 R_{N_\sigma}
 =U_\sigma
   \left((\sigma^2(L_N-1)+1)|_{\phi^\perp}\right)^{-1}
       \Pi_NU_\sigma^{-1}.
\]
Differentiation at $\sigma=1$ gives
\begin{equation}\label{eq:inverse-dilation}
 R'_N[\Dtwo N]=[\Lambda,R_N]-2R_N+2R_N^2,
 \qquad
 (R_N)_s=\frac{\lambda_s}{\lambda}R'_N[\Dtwo N]
                              -\lambda R'_N[\diver V].
\end{equation}
Differentiate $K_N\epsn=-R_N(Q_{a,\beta}\epsn)$ and insert the density
equation from \eqref{eq:actual-nonlinear-error}. The product rule gives
\[
 Q_{a,\beta}\Dtwo\epsn-\Lambda(Q_{a,\beta}\epsn)
       =[(2-\Lambda)Q_{a,\beta}]\epsn .
\]
The result is
\begin{equation}\label{eq:good-inverse-time}
 \begin{aligned}
 \frac{\lambda_s}{\lambda}\Lambda(K_N\epsn)-\partial_s(K_N\epsn)
 ={}&R_N\bigl[((Q_{a,\beta})_s-(\lambda_s/\lambda)\Lambda Q_{a,\beta})\epsn\bigr]
       +2\frac{\lambda_s}{\lambda}R_N^2(Q_{a,\beta}\epsn)\\
 &-\lambda R'_N[\diver V](Q_{a,\beta}\epsn)
       -\lambda R_N(Q_{a,\beta}\diver\epsv)\\
 &+\left(\frac{\lambda_s}{\lambda}+\beta\right)R_N(Q_{a,\beta}\Dtwo N_a)
       -a_sR_N(Q_{a,\beta}\partial_aN_a).
 \end{aligned}
\end{equation}
The potentially differentiated density has canceled. The
identities are first obtained for smooth densities; their
combination in \eqref{eq:good-inverse-time} extends by the
$H^{-1}$ tests of Lemma~\ref{lem:current-inverse}. In particular
the argument does not require an $L^2$ time derivative of $N$.
Moreover
\[
 \begin{aligned}
 (Q_{a,\beta})_s-\frac{\lambda_s}{\lambda}\Lambda Q_{a,\beta}
 ={}&-\ii L_aQ_{a,\beta}
      -\left(\frac{\lambda_s}{\lambda}+\beta\right)\Lambda Q_{a,\beta}\\
 &+a_s\partial_a Q_{a,\beta}
      +(\beta_s+\beta^2)\partial_\beta Q_{a,\beta}.
 \end{aligned}
\]
The localized inverse estimates of Lemma~\ref{lem:current-inverse}
give, term by term,
\[
\begin{aligned}
 &\left\|\frac{\lambda_s}{\lambda}\Lambda(K_N\epsn)
                         -\partial_s(K_N\epsn)\right\|_{H^1_d}\\
 &\quad\le C_I\left[
       \left(\beta+\left|\frac{\lambda_s}{\lambda}\right|\right)
                                                    \|\epsn\|_2
       +\lambda(1+\|\epsn\|_2)\|\epsv\|_2
       +\left|\frac{\lambda_s}{\lambda}+\beta\right|+|a_s|\right]\\
 &\quad\le C_I\lambda\|\veps\|_Y.
\end{aligned}
\]
In the second line, the part of $\diver V$ from $V_a$ has been
included in the $\beta\|\epsn\|_2$ term; the part from $\epsv$
gives $\lambda\|\epsn\|_2\|\epsv\|_2$. The last inequality uses
\eqref{eq:rate-bound}, \eqref{eq:beta-rate} and $\beta=c_0a\lambda$.
The same bounds control the source terms in \eqref{eq:good-source}
involving $\lambda_s/\lambda+\beta$, $a_s$ and $\beta_s+\beta^2$. For the remaining terms,
$L_aP_a=0$ and self-adjointness give
\[
\begin{aligned}
 \pair{L_aQ_{a,\beta}}\phi&=\pair{L_aQ_{a,\beta}}{\phi-\phi(N_a)},\\
 |\pair{L_aQ_{a,\beta}}\phi|&\le C_I\beta\|\epsn\|_2,\\
 \|(\gamma_s-1)R_NL_aQ_{a,\beta}\|_{H^1_d}
 &\le C_I\beta(\beta\|\veps\|_Y+\|\veps\|_Y^2)
       \le C_I\lambda\|\veps\|_Y.
\end{aligned}
\]
Together with \eqref{eq:combined-phase-bound}, these bounds prove
\eqref{eq:good-source-bound}.
\end{proof}

The self-adjoint part of the equation for $z$ preserves its
mass. Only this localized source contributes to the mass
derivative, with the required factor $\lambda$.

\subsection{Nonlinear interaction energy identity}
\label{sec:interaction-balance}
\label{app:quadratic}

We differentiate the interaction energy $E_\beta^{(2)}+H_3$
and the corrected mass $\|z\|_2^2$ separately. The fast terms
cancel exactly. We then collect the terms of order
$\beta\|\veps\|_Y^2$ to determine the quadratic correction.

\paragraph{Differentiation of the interaction energy.}
For the actual solution, the product rule gives
\[
\begin{aligned}
 \partial_s(E_\beta^{(2)}+H_3)
 ={}&2\Rea\pair{\eps_s}{L_N\eps+Q_{a,\beta}\epsn}
       +((\epsn)_s,\epsn+2\Rea(\overline{Q_{a,\beta}}\eps)+|\eps|^2)\\
 &+c_0^{-2}((\epsv)_s,\epsv)
       +a_s(\partial_aN_a,|\eps|^2)
       +2(\epsn,\Rea(\overline{(Q_{a,\beta})_s}\eps)).
\end{aligned}
\]
Substituting \eqref{eq:actual-nonlinear-error} produces the
exact cancellations
\[
\begin{aligned}
 &2\Rea\pair{-\ii(L_N\eps+Q_{a,\beta}\epsn)}{L_N\eps+Q_{a,\beta}\epsn}=0,\\
 &-\lambda\bigl[(\diver\epsv,
           \epsn+2\Rea(\overline{Q_{a,\beta}}\eps)+|\eps|^2)\\
 &\hspace{34mm}+(\epsv,\nabla
           (\epsn+2\Rea(\overline{Q_{a,\beta}}\eps)+|\eps|^2))\bigr]=0.
\end{aligned}
\]
The phase term also has a localized small coefficient. Indeed
$L_N\eps+Q_{a,\beta}\epsn=L_N\Psi-L_aQ_{a,\beta}$, so self-adjointness gives
\begin{equation}\label{eq:phase-work}
 2\Rea\pair{-\ii(\gamma_s-1)\Psi}{L_N\eps+Q_{a,\beta}\epsn}
       =-2(\gamma_s-1)\Ima\pair\eps{L_aQ_{a,\beta}}.
\end{equation}
Here $\Ima\pair{Q_{a,\beta}}{L_aQ_{a,\beta}}=0$ and $L_aQ_{a,\beta}=O_{H^1_d}(\beta)$.

\begin{remark}[The phase error and the slow energy estimate]
\label{rem:phase-work}
The term without $\lambda$ in \eqref{eq:phase-bound} comes from
the eigenvalue shift $\mu(N)+1=O_I(\|\epsn\|_2^2)$ under the
density constraint. It changes the rotation rate of the field;
it is not by itself an energy growth term. In the equation for $z$,
the real coefficient $\gamma_s-1$ is part of a self-adjoint
operator and does no mass work. The source along $\phi(N)$ also
does no mass work, since $z\perp\phi(N)$.

The remaining phase terms in \eqref{eq:phase-work} and
\eqref{eq:good-source} contain $L_aQ_{a,\beta}$ or $R_NL_aQ_{a,\beta}$.
Both are localized and $O_I(\beta)$, so
\[
\begin{aligned}
 &|\gamma_s-1|
   \bigl(|\pair\eps{L_aQ_{a,\beta}}|+|\pair{R_NL_aQ_{a,\beta}}z|\bigr)\\
 &\qquad\le C_I\beta(\beta+\|\veps\|_Y)\|\veps\|_Y^2.
\end{aligned}
\]
The quadratic correction \eqref{eq:B}, constructed below, has an explicit factor
$\beta$; its phase work has the same bound, using the
orthogonality of its test to $P_a$, as verified below.
Thus the quadratic term in the bound for $\gamma_s-1$ contributes
to the energy estimate with a factor of $\beta$. No integrability of
$\gamma_s-1$ is required: the phase is modulated throughout,
and the theorem does not assert a limit for $\gamma(s)-s$.
\end{remark}

\paragraph{Differentiation of the corrected mass.}
Both $L_N+\gamma_s-1$ and the dilation have zero real mass work.
Thus \eqref{eq:nonlinear-good-equation} gives
\[
 \partial_s\|z\|_2^2
 =2\Rea\pair{-\ii(L_N+\gamma_s-1)z
              +(\lambda_s/\lambda)\Lambda z-\ii\mathcal F}{z}
 =2\Ima\pair{\mathcal F}z.
\]
This calculation initially uses smooth decaying functions.
The canceled right side is a continuous local pairing by
\eqref{eq:good-source-bound}; its integrated extension does
not impose a spatial moment on the energy solution.

For the cubic dilation terms, integration by parts gives
\[
 (\Dtwo\epsn,|\eps|^2)
       +2\Rea\pair{\Lambda\eps}{\epsn\eps}=2H_3.
\]
The quadratic dilation terms give the same $D_{a,\beta}(\veps)$
as in the fixed-orbit calculation, defined in
\eqref{eq:global-linear-identity}. The mixed dilation terms combine
by the product rule and $(\Dtwo g,g)=\|g\|_2^2$ for each wave
component.

Combining the preceding identities yields
\begin{equation}\label{eq:action-law}
\begin{aligned}
 E_s-2(\lambda_s/\lambda)H_3
 ={}&\frac{\lambda_s}{\lambda}D_{a,\beta}(\veps)
       +2(\epsn,\Rea(\overline{(Q_{a,\beta})_s}\eps))
       -2(\gamma_s-1)\Ima\pair\eps{L_aQ_{a,\beta}}\\
 &+2\Rea\left\langle
       \left(\frac{\lambda_s}{\lambda}+\beta\right)\Lambda Q_{a,\beta}
           -a_s\partial_a Q_{a,\beta},\,L_N\eps+Q_{a,\beta}\epsn\right\rangle\\
 &-2(\beta_s+\beta^2)\Rea\pair{\partial_\beta Q_{a,\beta}}
                                    {L_N\eps+Q_{a,\beta}\epsn}\\
 &+\left(\frac{\lambda_s}{\lambda}+\beta\right)
       (\Dtwo N_a,\epsn+2\Rea(\overline{Q_{a,\beta}}\eps)+|\eps|^2)\\
 &-a_s(\partial_aN_a,\epsn+2\Rea(\overline{Q_{a,\beta}}\eps))\\
 &+c_0^{-2}\left(
       \left(\frac{\lambda_s}{\lambda}+\beta\right)\Dtwo V_a
                   -a_s\partial_aV_a,\epsv\right)
       -2\Ima\pair{\mathcal F}z .
\end{aligned}
\end{equation}
The $a_s(\partial_aN_a,|\eps|^2)$ from the product rule has
canceled the same term in the density forcing. Each pairing on
the right-hand side of \eqref{eq:action-law} extends continuously
to $Y$ with the parameter rates fixed.
Laplacians are transferred to the localized profile tests;
the polynomial density tests multiply $|\eps|^2$ through
\eqref{eq:weighted-local}.

The identity has only wave-scale terms, but the mixed and
parameter pairings still have no sign. We next separate them
from the diagonal forms before choosing the quadratic correction.

\subsection{Principal energy balance and modulation terms}
\label{sec:principal-work}

\paragraph{Leading modulation equations.}
The amplitude equation must be kept to its leading order,
because its work is as large as the intended dissipation.
The following expansion also shows where the density response
enters that equation.
The density response satisfies $K_N\epsn\perp\phi$. Thus
\eqref{eq:chart} gives the exact orthogonality and its reference
moment bounds
\begin{equation}\label{eq:nonlinear-moment-bound}
 \pair z\phi=\pair\eps\phi=0,\qquad
 |\pair z{P_a}|+|\pair\eps{P_a}|
                         \le C_I\|\veps\|_Y^2.
\end{equation}
Indeed, pair $z$ and $\eps$ with $P_a-\sqrt M\phi$ and use
$\|P_a-\sqrt M\phi\|_{H^2_d}\le C_I\|\epsn\|_2$,
$\|z\|_{\mathrm{loc}}\le C_I\|\veps\|_Y$ and the localized
pairing estimate.

The scale and branch outputs needed for the quadratic correction are
\begin{align}
 \left|\frac{\lambda_s}{\lambda}+\beta-\frac{\beta(\epsn,y\cdot\nabla P_a^2)
                         +\lambda(\epsv,\nabla P_a^2)}{2M}\right|
       &\le C_I\beta\|\veps\|_Y^2,\label{eq:m-output}\\
 \left|a_s-\frac{2\beta}{M'}(\Rea z,\Lambda P_a)\right|
       &\le C_I\beta(\beta+\|\veps\|_Y)\|\veps\|_Y.
       \label{eq:eta-linear}
\end{align}
For \eqref{eq:m-output}, expand the denominator in
\eqref{eq:density-eq}; each remaining coefficient is a linear
density pairing multiplying a rate of size $C_I\beta\|\veps\|_Y$.

To expand the branch equation \eqref{eq:amplitude-eq}, use the term
linear in $\beta$ in $L_aQ_{a,\beta}$. The eigenfunction expansion
and profile identities give
\[
\begin{aligned}
 \phi-\phi(N_a)&=-M^{-1/2}R_a(P_a\epsn)+O_{H^2_d}(\|\epsn\|_2^2),\\
 L_aQ_{a,\beta}&=\ii\beta\Lambda P_a+O_{H^1_d}(\beta^2),\\
 \phi'_{N_a}[\diver V_a]&=-c_0a\Lambda P_a/\sqrt M.
\end{aligned}
\]
Use also $\pair{L_aQ_{a,\beta}}{\phi(N_a)}=0$,
$\pair{Q_{a,\beta}}\phi=\sqrt M+O_I(\beta+\|\epsn\|_2)$, and the
local norm bound on $\eps$.
The two forcing terms in \eqref{eq:amplitude-eq} then satisfy
\[
\begin{aligned}
 \Rea\frac{\ii\pair{L_aQ_{a,\beta}}\phi}{\pair{Q_{a,\beta}}\phi}
   &=\frac\beta M(R_a(P_a\epsn),\Lambda P_a)
          +O_I(\beta(\beta+\|\veps\|_Y)\|\veps\|_Y),\\
 -\lambda\Rea\frac{\pair\eps{\phi'_N[\diver V]}}{\pair{Q_{a,\beta}}\phi}
   &=\frac\beta M(\Rea\eps,\Lambda P_a)
          +O_I(\beta(\beta+\|\veps\|_Y)\|\veps\|_Y).
\end{aligned}
\]
The wave perturbation in the second line contributes
$C_I\lambda\|\eps\|_{\mathrm{loc}}\|\epsv\|_2$;
the difference of eigenfunction derivatives adds one density factor.
The scale-defect term in \eqref{eq:amplitude-eq} is
$O_I(\beta\|\veps\|_Y^2)$, since its profile pairing vanishes
at zero density error. The branch coefficient satisfies
\[
 \Rea\frac{\pair{\partial_a Q_{a,\beta}+(\beta/a)\partial_\beta Q_{a,\beta}}\phi}{\pair{Q_{a,\beta}}\phi}
       =\frac{M'}{2M}+O_I(\beta^2+\|\epsn\|_2).
\]
The two forcing terms combine into the projection of
$\Rea(\eps+R_a(P_a\epsn))$. We use the comparison
\begin{equation}\label{eq:response-comparison}
 \|z-\eps-R_a(P_a\epsn)\|_{H^2_d}
 =\|R_N(Q_{a,\beta}\epsn)-R_a(P_a\epsn)\|_{H^2_d}
 \le C_I(\beta+\|\epsn\|_2)\|\epsn\|_2.
\end{equation}
Indeed, the difference of the inverses costs one density factor
and the difference $Q_{a,\beta}-P_a$ costs $\beta$. After division by the
branch coefficient, this gives the leading term
$2\beta(\Rea z,\Lambda P_a)/M'$ in \eqref{eq:eta-linear}.

\paragraph{The diagonal and mixed forms.}
The nonlinear derivative uses the same $\mathcal M_a$ and $D_a$
as the fixed-orbit linear calculation. Write $\chi_a=R_a\Lambda P_a$.
The branch variation is the real
linear wave functional
\begin{equation}\label{eq:branch-work}
 \mathcal J_a(f,w)
       =(f,\partial_a(N_a+P_a^2))+c_0^{-2}(w,\partial_aV_a).
\end{equation}
It is bounded on $L^2\times L^2$. It appears because
differentiating the profile in $a$ changes both wave components.

The remaining work has three sources: the diagonal dilation
forms, the mixed field--wave pairing, and the motion of the
reference profiles. The scale work must ultimately stay in the
velocity form; the mixed and branch work will be canceled.

\begin{lemma}[Principal nonlinear energy balance]\label{lem:principal-balance}
For a smooth modulated solution on the bootstrap regime,
the energy \eqref{eq:nonlinear-mass-energy}
has the following principal balance:
\begin{equation}\label{eq:direct-principal}
\begin{aligned}
 \bigg|E_s-2(\lambda_s/\lambda)H_3
 &+\beta\bigl[(A_az,z)_\R+D_a(\epsn)+c_0^{-2}\|\epsv\|_2^2\bigr]\\
 &+2\beta\bigl(\Rea z,
       \mathcal M_a\epsn-(c_0a)^{-1}R_a(P_a\diver \epsv)\bigr)\\
 &-\left(\frac{\lambda_s}{\lambda}+\beta\right)\bigl[(\epsn,\Dtwo(N_a+P_a^2))
                         +c_0^{-2}(\epsv,\Dtwo V_a)\bigr]
       +a_s\mathcal J_a(\epsn,\epsv)\bigg|\\
 &\hspace{12mm}\le C_I\beta(\beta+\|\veps\|_Y)\|\veps\|_Y^2.
\end{aligned}
\end{equation}
\end{lemma}

\begin{proof}
\label{app:nonlinear-balance-details}
We first collect the corrected-mass and dilation terms, then
the parameter terms.

\paragraph{The density response and dilation.}
The derivative of the density response introduces two terms
containing $R_a^2$. They cancel before estimation, leaving the same commutator
that occurred in the linear calculation.
For the corrected mass, the profile identities are
\[
 L_aQ_{a,\beta}=\ii\beta\Lambda P_a+O_{H^1_d}(\beta^2),\qquad
 \diver V_a=-c_0a\Dtwo N_a.
\]
The leading terms linear in $\epsn$ combine exactly as follows:
\begin{equation}\label{eq:leading-response-cancellation}
\begin{aligned}
 &\beta R_a[(\Lambda P_a)\epsn]
       +\beta R'_{N_a}[\Dtwo N_a](P_a\epsn)
       -2\beta R_a^2(P_a\epsn)\\
 &\quad=\beta R_a[(\Lambda P_a)\epsn]
       +\beta([\Lambda,R_a]-2R_a)(P_a\epsn)\\
 &\quad=\beta\bigl[\Lambda R_a(P_a\epsn)
                              -R_a(P_a\Dtwo\epsn)\bigr]\\
 &\quad=\beta\bigl\{[\Lambda,R_a](P_a\epsn)
                           -R_a([(2-\Lambda)P_a]\epsn)\bigr\}.
\end{aligned}
\end{equation}
The first equality uses \eqref{eq:inverse-dilation} and cancels
the two $R_a^2$ terms. The remaining equalities follow from
the product rule in \eqref{eq:K-commutator}.
The divergence and parameter terms then give
\[
\begin{aligned}
 &\frac{\lambda_s}{\lambda}\Lambda(K_N\epsn)
                   -\partial_s(K_N\epsn)\\
 &\quad=\beta\bigl\{[\Lambda,R_a](P_a\epsn)
                    -R_a([(2-\Lambda)P_a]\epsn)\bigr\}
       -\lambda R_a(P_a\diver\epsv)\\
 &\qquad-\left(\frac{\lambda_s}{\lambda}+\beta\right)\Lambda P_a
       -a_sR_a(P_a\partial_aN_a)
       +O_{H^1_d}\bigl(\beta(\beta+\|\veps\|_Y)\|\veps\|_Y\bigr).
\end{aligned}
\]
The last two terms use the profile identities
\[
 -R_a(P_a\Dtwo N_a)=\Lambda P_a,\qquad
 -R_a(P_a\partial_aN_a)=\partial_aP_a-\frac{M'}{2M}P_a.
\]
Combining these terms with the explicit profile forcing in
\eqref{eq:good-source} gives the two cancellations
\[
\begin{aligned}
 &-\ii\left(\frac{\lambda_s}{\lambda}+\beta\right)\Lambda P_a
       +\ii\left(\frac{\lambda_s}{\lambda}+\beta\right)\Lambda P_a=0,\\
 &-\ii a_s(R_a(P_a\partial_aN_a)+\partial_aP_a)
                         =-\ii a_s\frac{M'}{2M}P_a.
\end{aligned}
\]
Thus the scale forcing cancels, and the branch forcing is along
the reference ground state. The mass derivative pairs the
remaining term with $z$; by \eqref{eq:nonlinear-moment-bound},
\[
 \left|2\Ima\pair{-\ii a_s\tfrac{M'}{2M}P_a}z\right|
       \le C_I\beta(\beta+\|\veps\|_Y)\|\veps\|_Y^2.
\]
Replacing $Q_{a,\beta}$ by $P_a$ and $\partial_a Q_{a,\beta}$ by
$\partial_aP_a$ costs the same allowed remainder. By
\eqref{eq:beta-rate}, the source
$-\ii(\beta_s+\beta^2)\partial_\beta Q_{a,\beta}$ has localized norm
at most $C_I\beta^2\|\veps\|_Y$ and hence has work at most
$C_I\beta^2\|\veps\|_Y^2$ against $z$.

The first two terms of \eqref{eq:good-source} point along $\phi$.
Their mass work vanishes exactly by \eqref{eq:nonlinear-moment-bound}:
\[
 \Ima\pair{
   [(\mu(N)+\gamma_s)\pair{Q_{a,\beta}}\phi-\pair{L_aQ_{a,\beta}}\phi]\phi}{z}=0.
\]
For the third term and the phase work in
\eqref{eq:phase-work}, the localized product estimates give
\[
 \begin{aligned}
 |\pair{R_NL_aQ_{a,\beta}}z|+|\pair\eps{L_aQ_{a,\beta}}|
          &\le C_I\beta\|\veps\|_Y,\\
 |\gamma_s-1|
    \bigl(|\pair{R_NL_aQ_{a,\beta}}z|+|\pair\eps{L_aQ_{a,\beta}}|\bigr)
          &\le C_I\beta(\beta+\|\veps\|_Y)\|\veps\|_Y^2.
 \end{aligned}
\]
Here we used $|\gamma_s-1|\le
C_I(\beta\|\veps\|_Y+\|\veps\|_Y^2)$ and
$R_NL_aQ_{a,\beta}=O_{H^2_d}(\beta)$.

Split the dilation term into $-\beta D_{a,\beta}$ and
$(\lambda_s/\lambda+\beta)D_{a,\beta}$, keeping the latter
with the scale terms below. In the first term,
$|D_{a,\beta}(\veps)-D_{a,0}(\veps)|\le C_I\beta\|\veps\|_Y^2$
by the mixed term in \eqref{eq:global-linear-identity} and the
localized bound on $Q_{a,\beta}-P_a$. This replacement costs
$C_I\beta^2\|\veps\|_Y^2$. In $D_{a,0}$, use
\eqref{eq:response-comparison} to write the scalar and mixed
terms in $z$, as in \eqref{eq:diagonal-expansion}, with the
density term $D_a(\epsn)$ unchanged. After multiplication by
$\beta$, the response replacement costs at most
$C_I\beta(\beta+\|\veps\|_Y)\|\veps\|_Y^2$.
\paragraph{The parameter terms.}
We now collect the scale and branch contributions in
\eqref{eq:action-law}. For the coefficient of
$\lambda_s/\lambda+\beta$, first set $\beta=0$ in the profile
factors, keeping $a$ unchanged. Its linear part in the errors is
\[
\begin{aligned}
 &2\Rea\pair{\Lambda P_a}{L_a\eps+P_a\epsn}
       +(\Dtwo N_a,\epsn+2P_a\Rea\eps)
       +c_0^{-2}(\Dtwo V_a,\epsv)\\
 &\quad=2(\Rea\eps,L_a\Lambda P_a+P_a\Dtwo N_a)
       +(\epsn,2P_a\Lambda P_a+\Dtwo N_a)
       +c_0^{-2}(\epsv,\Dtwo V_a)\\
 &\quad=-4(\Rea\eps,P_a)
       +(\epsn,\Dtwo(N_a+P_a^2))
       +c_0^{-2}(\epsv,\Dtwo V_a).
\end{aligned}
\]
The last equality uses
$L_a\Lambda P_a=-2P_a-P_a\Dtwo N_a$ and
$2P_a\Lambda P_a=\Dtwo(P_a^2)$.
The scale coefficient also contains quadratic error terms,
including $D_{a,\beta}(\veps)$ and the terms with $\epsn\eps$
or $|\eps|^2$. Their size is at most $C_I\|\veps\|_Y^2$.
Replacing $Q_{a,\beta}$ and $\Lambda Q_{a,\beta}$
by $P_a$ and $\Lambda P_a$ adds at most
$C_I\beta\|\veps\|_Y$. Since
$\left|\frac{\lambda_s}{\lambda}+\beta\right|\le C_I\beta\|\veps\|_Y$, these errors and the scalar
pairing $-4\left(\frac{\lambda_s}{\lambda}+\beta\right)(\Rea\eps,P_a)$ are bounded by the right side
of \eqref{eq:direct-principal}.

There are two explicit phase-derivative terms in
\eqref{eq:action-law}. By \eqref{eq:beta-rate},
$|\beta_s|\le C_I\beta^2$. Using the localization of
$\partial_\beta Q_{a,\beta}$ and transposing $L_a$ onto this
profile derivative gives
\[
\begin{aligned}
 &2|\beta_s|\left|(\epsn,\Rea(\overline{\partial_\beta Q_{a,\beta}}\eps))\right|
       \le C_I\beta^2\|\epsn\|_2\|\eps\|_{\mathrm{loc}},\\
 &2|\beta_s+\beta^2|
       \left|\pair{\partial_\beta Q_{a,\beta}}{L_N\eps+Q_{a,\beta}\epsn}\right|
       \le C_I\beta^2\|\veps\|_Y^2 .
\end{aligned}
\]
In the second line the localized pairing is at most $C_I\|\veps\|_Y$,
by the weighted product estimates for the $\epsn\eps$ term.
The terms proportional to $a_s$ combine exactly as follows:
\[
\begin{aligned}
 &2(\epsn,\Rea(\overline{\partial_a Q_{a,\beta}}\eps))
       -2\Rea\pair{\partial_a Q_{a,\beta}}{L_N\eps+Q_{a,\beta}\epsn}\\
 &\qquad-(\partial_aN_a,\epsn+2\Rea(\overline{Q_{a,\beta}}\eps))
       -c_0^{-2}(\partial_aV_a,\epsv)\\
 &\quad=-2\Rea\pair{L_a\partial_a Q_{a,\beta}}{\eps}
       -(\epsn,\partial_a(P_a^2))
       -(\partial_aN_a,\epsn+2\Rea(\overline{Q_{a,\beta}}\eps))
       -c_0^{-2}(\partial_aV_a,\epsv)\\
 &\quad=-\mathcal J_a(\epsn,\epsv)
       -2\Rea\pair{\partial_a(L_aQ_{a,\beta})}\eps .
\end{aligned}
\]
Here $L_N=L_a+\epsn$ cancels the two terms containing
$\epsn\eps$, while
$L_a\partial_a Q_{a,\beta}=\partial_a(L_aQ_{a,\beta})-(\partial_aN_a)Q_{a,\beta}$
cancels the remaining Schr\"odinger--density pairing.
The last scalar term is small because
\[
 L_aQ_{a,\beta}=\e^{-\ii\beta r^2/4}
       \left(\ii\beta\Lambda P_a+\tfrac14\beta^2r^2P_a\right),
 \qquad \|\partial_a(L_aQ_{a,\beta})\|_{H^1_d}\le C_I\beta,
\]
The derivative here is at fixed $\beta$. In particular,
\[
 |a_s\pair{\partial_a(L_aQ_{a,\beta})}\eps|
       \le C_I\beta^2\|\veps\|_Y^2.
\]
Together with the corrected-mass calculation, these identities
prove \eqref{eq:direct-principal}.

All replacements in this calculation preserve their scale
factor. In fact
\[
 R_N-R_a=\int_0^1R'_{N_a+\theta\epsn}[\epsn]\,\dd\theta.
\]
The inverse bounds in Lemma~\ref{lem:current-inverse} show
that such a replacement adds one $L^2$ density factor, also
when the source contains one velocity divergence. Replacing a
smooth profile by its zero-phase value adds $\beta$.
The other error products use the localized estimates
\[
 |\pair{\epsn\eps}p|
       \le C_I\|\epsn\|_2\|\eps\|_{\mathrm{loc}}\|p\|_{H^1_d},
 \qquad
 \int\langle r\rangle^{-3}|\eps|^2
       \le C\|\eps\|_{\mathrm{loc}}^2.
\]
Together with the exact phase cancellation, these bounds account
for the full right side of \eqref{eq:direct-principal}.
\end{proof}

\subsection{Quadratic energy correction}\label{sec:quadratic-correction}
In \eqref{eq:direct-principal}, the mixed Schr\"odinger--wave
pairing and $a_s\mathcal J_a$ have no sign. To cancel both,
choose the real localized test $p\perp P_a$ solving
\begin{equation}\label{eq:direct-test}
\begin{aligned}
 L_ap={}&\Pi_a[\mathcal M_a\epsn-(c_0a)^{-1}R_a(P_a\diver \epsv)]
                 +\frac{\Lambda P_a}{M'}\mathcal J_a(\epsn,\epsv),\\
 p={}&R_a\mathcal M_a\epsn-(c_0a)^{-1}R_a^2(P_a\diver\epsv)
                 +\frac{\chi_a}{M'}\mathcal J_a(\epsn,\epsv).
\end{aligned}
\end{equation}
The right-hand side of the first line is orthogonal to $P_a$
because $(\Lambda P_a,P_a)=0$; hence the second line is its
unique solution in $P_a^\perp$. Set
\begin{equation}\label{eq:B}
 B_\beta=-2\beta(\Ima\eps,p),\qquad
 \E=E+B_\beta.
\end{equation}
The first two summands in $p$ give the fixed-orbit correction test.
The term with $\chi_a$ cancels $a_s\mathcal J_a$ through
\eqref{eq:eta-linear}. All three are linear in the wave errors.

The evolution of this test is slow:
only wave equations and parameter rates enter its derivative.
We record the estimate before using the fast field equation.

\begin{lemma}[Evolution of the correction test]\label{lem:correction-test-time}
For a smooth modulated solution on the bootstrap regime,
the test \eqref{eq:direct-test} satisfies
\begin{equation}\label{eq:direct-test-time}
 \|p\|_{H^1_d}\le C_I\|\veps\|_Y,
 \qquad \|p_s\|_{H^1_d}\le C_I\lambda
                              (\|\veps\|_Y+\|\veps\|_Y^2).
\end{equation}
\end{lemma}
Its distributional source estimates are proved in
Appendix~\ref{app:correction-details}.

\paragraph{Cancellation of the mixed and branch terms.}
To differentiate $B_\beta$, use the imaginary part of the
Schr\"odinger error equation:
\[
\begin{aligned}
 \Ima\eps_s={}&-L_a\Rea\eps-P_a\epsn
       -\epsn\Rea\eps-\epsn(\Rea Q_{a,\beta}-P_a)
       -(\gamma_s-1)\Rea\Psi\\
 &+\frac{\lambda_s}{\lambda}\Lambda\Ima\eps
       +\left(\frac{\lambda_s}{\lambda}+\beta\right)\Ima(\Lambda Q_{a,\beta})\\
 &-a_s\Ima(\partial_a Q_{a,\beta})
       -(\beta_s+\beta^2)\Ima(\partial_\beta Q_{a,\beta}).
\end{aligned}
\]
Pair this equation with $p$. From $\Rea Q_{a,\beta}=P_a\cos(\beta r^2/4)$
and the exponential profile bounds,
\[
 \begin{gathered}
 \|\Rea Q_{a,\beta}-P_a\|_{H^1_d}\le C_I\beta^2,\\
 \|\Ima(\Lambda Q_{a,\beta})\|_{H^1_d}
       +\|\Ima(\partial_a Q_{a,\beta})\|_{H^1_d}\le C_I\beta.
 \end{gathered}
\]
The two density pairings satisfy
\[
\begin{aligned}
 \beta|(\epsn\Rea\eps,p)|
       &\le C_I\beta\|\veps\|_Y^3,\\
 \beta|(\epsn(\Rea Q_{a,\beta}-P_a),p)|
       &\le C_I\beta^3\|\veps\|_Y^2.
\end{aligned}
\]
For the phase contribution, $(p,P_a)=0$ gives
\[
\begin{aligned}
 (\Rea\Psi,p)&=(\Rea\eps+\Rea Q_{a,\beta}-P_a,p),\\
 \beta|(\gamma_s-1)(\Rea\Psi,p)|
 &\le C_I\beta(\beta\|\veps\|_Y+\|\veps\|_Y^2)
                    (\|\veps\|_Y+\beta^2)\|\veps\|_Y\\
 &\le C_I\beta(\beta+\|\veps\|_Y)\|\veps\|_Y^2.
\end{aligned}
\]
The last step uses the bootstrap bounds $\beta,\|\veps\|_Y<1$.
The remaining modulation terms are bounded after pairing with
$p$. The dilation of $\eps$ is tested by
$\pair{\Lambda\eps}p=\pair\eps p+\pair{\nabla\eps}{yp}$.
The scale and branch terms use the two $O_{H^1_d}(\beta)$ bounds
above and \eqref{eq:rate-bound}. For the term with $\beta_s+\beta^2$,
use \eqref{eq:beta-rate} and the uniform localized bound on
$\partial_\beta Q_{a,\beta}$. Including the prefactor $\beta$,
their total is at most $C_I\beta^2\|\veps\|_Y^2$.
Finally \eqref{eq:beta-rate} and \eqref{eq:direct-test-time} give
\[
 |\beta_s(\Ima\eps,p)|+\beta|(\Ima\eps,p_s)|
       \le C_I\beta^2(1+\|\veps\|_Y)\|\veps\|_Y^2.
\]

The product rule therefore gives
\[
\begin{aligned}
 (B_\beta)_s
 &=-2\beta_s(\Ima\eps,p)
       -2\beta(\Ima\eps_s,p)-2\beta(\Ima\eps,p_s)\\
 &=2\beta(\Rea\eps,L_ap)+2\beta(P_a\epsn,p)
       +O_I\bigl(\beta(\beta+\|\veps\|_Y)\|\veps\|_Y^2\bigr).
\end{aligned}
\]
For the main terms, use the adjoint equation
\eqref{eq:direct-test} and the self-adjointness of $R_a$:
\[
 (P_a\epsn,R_ag)=(R_a(P_a\epsn),g),\qquad
 (P_a\epsn,\chi_a)=(R_a(P_a\epsn),\Lambda P_a).
\]
The first identity holds for the localized real sources $g$ in
\eqref{eq:direct-test}, including its distributional divergence
source. Replacing $\Pi_a$ by the identity in that equation
costs at most
\[
 C_I\beta|\pair\eps{P_a}|\|\veps\|_Y
       \le C_I\beta(\beta+\|\veps\|_Y)\|\veps\|_Y^2.
\]
The adjoint identities combine the field and density pairings into
the pairing with $\Rea(\eps+R_a(P_a\epsn))$. By
\eqref{eq:response-comparison} this can be written with $\Rea z$.
The branch coefficient $2\beta(\Rea z,\Lambda P_a)/M'$ is
$a_s$ up to the error in \eqref{eq:eta-linear}. Multiplication by
$\mathcal J_a(\epsn,\epsv)$ gives the same allowed remainder.
Thus
\begin{equation}\label{eq:direct-correction-derivative}
\begin{aligned}
 \bigg|(B_\beta)_s
 &-2\beta(\Rea z,\mathcal M_a\epsn-(c_0a)^{-1}R_a(P_a\diver \epsv))
       -a_s\mathcal J_a(\epsn,\epsv)\bigg|\\
 &\hspace{12mm}\le C_I\beta(\beta+\|\veps\|_Y)\|\veps\|_Y^2.
\end{aligned}
\end{equation}
Thus $B_\beta$ cancels the mixed pairing and
$a_s\mathcal J_a$. The scale-rate term in
\eqref{eq:direct-principal} remains; it will be included in
the velocity dissipation form.

\subsection{Coercivity and nonlinear dissipation}\label{sec:nonlinear}\label{app:energy}

\begin{lemma}[Coercivity]\label{lem:coercivity}
Under \eqref{eq:chart}, for sufficiently small
$\beta$, $\|\veps\|_Y$ and $\|\eps\|_2$,
$c_I\|\veps\|_Y^2\le\E\le C_I\|\veps\|_Y^2$.
\end{lemma}
\begin{proof}
Use the exact square \eqref{eq:completed-square} at zero phase.
In \eqref{eq:relative-action} the replacement of the phase and
response coefficients is bounded by \eqref{eq:response-comparison}.
Writing the scalar square with $z$ gives
\begin{equation}\label{eq:energy-comparison}
\begin{aligned}
 E={}&((L_a-1)z,z)_\R
       +\tfrac12(\epsn,(1-2P_aR_aP_a)\epsn)
       +\tfrac1{2c_0^2}\|\epsv\|_2^2+H_3\\
 &+O_I((\beta+\|\veps\|_Y)\|\veps\|_Y^2).
\end{aligned}
\end{equation}
Lemma~\ref{lem:limiting-energy} and its positive-$a$ perturbation
argument in Appendix~\ref{app:energy-coercivity} give
\[
\begin{aligned}
 ((L_a-1)z,z)_\R
    &\ge c_I\|z\|_{\mathrm{loc}}^2-C_I|\pair z{P_a}|^2,\\
 (\epsn,(1-2P_aR_aP_a)\epsn)&\ge c_I\|\epsn\|_2^2
                     \qquad((\epsn,P_a^2)=0).
\end{aligned}
\]
For completeness, the scalar penalty follows from the exact
decomposition
\[
\begin{aligned}
 z_\perp&=z-M^{-1}\pair z{P_a}P_a,\qquad \pair{z_\perp}{P_a}=0,\\
 ((L_a-1)z,z)_\R
       &=((L_a-1)z_\perp,z_\perp)_\R
                                 -M^{-1}|\pair z{P_a}|^2,\\
 \|z\|_{\mathrm{loc}}^2
       &\le2\|z_\perp\|_{\mathrm{loc}}^2+C_I|\pair z{P_a}|^2.
\end{aligned}
\]
The cross term vanishes because $(L_a-1)P_a=-P_a$.
Apply the constrained scalar gap to $z_\perp$.
Since $\eps=z+K_N\epsn$ and $K_N:L^2\to H^2_d$ is bounded,
\[
 \|\veps\|_Y^2
 \le C_I\bigl(\|z\|_{\mathrm{loc}}^2+\|\epsn\|_2^2
                                      +c_0^{-2}\|\epsv\|_2^2\bigr).
\]
Thus the principal form in \eqref{eq:energy-comparison} is bounded
below by $c_I\|\veps\|_Y^2-C_I|\pair z{P_a}|^2$.
Now use \eqref{eq:nonlinear-moment-bound},
$|B_\beta|\le C_I\beta\|\veps\|_Y^2$, and \eqref{eq:H3}.
Choosing the three thresholds small absorbs the errors into
half the lower bound. The upper bound follows from the same
continuous forms and $\|K_N\epsn\|_{H^2_d}\le C_I\|\epsn\|_2$.
\end{proof}

\begin{proof}[Proof of Proposition~\ref{prop:remainder}]
Adding \eqref{eq:direct-principal} and
\eqref{eq:direct-correction-derivative} gives
\begin{equation}\label{eq:master}
\begin{aligned}
 \bigg|\E_s-2(\lambda_s/\lambda)H_3
 &+\beta\bigl[(A_az,z)_\R+D_a(\epsn)+c_0^{-2}\|\epsv\|_2^2\bigr]\\
 &-\left(\frac{\lambda_s}{\lambda}+\beta\right)\bigl[(\epsn,\Dtwo(N_a+P_a^2))
                         +c_0^{-2}(\epsv,\Dtwo V_a)\bigr]\bigg|\\
 &\hspace{12mm}\le C_I\beta(\beta+\|\veps\|_Y)\|\veps\|_Y^2.
\end{aligned}
\end{equation}
The mixed and branch terms have canceled. We now choose the
small branch bound $a_*$ from the leading forms. In this part of
the proof, $C$ denotes a constant uniform for small positive $a$;
constants $C_I$ are reserved for higher-order errors after $I$
is fixed. We retain the scale-rate term in the velocity form.
The Schr\"odinger dissipation inequality is
\[
 (A_az,z)_\R\ge c_I\|z\|_{\mathrm{loc}}^2
                                      -C_I|\pair z{P_a}|^2,
\]
proved by removing the same explicit ground-state component.
The density inequality controls $D_a(\epsn)$ on
$(\epsn,P_a^2)=0$. The scale contribution changes the velocity
form and must be retained. Indeed
\[
 \|\Dtwo(N_a+P_a^2)\|_2\le Ca^2,\qquad
 \Dtwo V_a=-c_0a\Dtwo(yN_a).
\]
Substitute the leading scale quotient from \eqref{eq:m-output}.
The error is $O_I(\beta\|\veps\|_Y^2)$ and multiplies a
bounded linear wave pairing, so it contributes at most
$C_I\beta\|\veps\|_Y^3$. The remaining products containing
$\epsn$ are bounded by
\[
 C\beta\bigl(a^2\|\epsn\|_2^2
                         +a\|\epsn\|_2\|\epsv\|_2\bigr).
\]
Indeed the density part of the leading quotient divided by
$\beta$ has size $C\|\epsn\|_2$, its velocity part has size
$Ca^{-1}\|\epsv\|_2$, and the density and velocity tests
in the scale bracket have norms $O(a^2)$ and $O(a)$.
The velocity square in the negative bracket is exactly
\[
\begin{aligned}
 &c_0^{-2}\|\epsv\|_2^2
   -\frac{\lambda(\epsv,\nabla P_a^2)}{2M\beta}
                         c_0^{-2}(\epsv,\Dtwo V_a)\\
 &\qquad=c_0^{-2}\left[\|\epsv\|_2^2
     +\frac{(\epsv,\nabla P_a^2)(\epsv,\Dtwo(yN_a))}{2M}\right],
\end{aligned}
\]
where $\beta=c_0a\lambda$ and $V_a=-c_0a\,yN_a$ fix both
the sign and the denominator. Its limit as $a\downarrow0$ is
\begin{equation}\label{eq:velocity-block}
 c_0^{-2}\left[\|{\epsv}\|_2^2+
       \frac{({\epsv},\nabla Q^2)({\epsv},\Dtwo(-yQ^2))}{\int Q^4}\right].
\end{equation}
Here $Q$ is the limiting NLS ground state and ${\epsv}$ is an arbitrary
real radial velocity direction. Lemma~\ref{lem:velocity-analytic}
gives
\begin{equation}\label{eq:spectral-velocity}
 \|{\epsv}\|_2^2+
       \frac{({\epsv},\nabla Q^2)({\epsv},\Dtwo(-yQ^2))}{\int Q^4}
                         \ge\tfrac1{50}\|{\epsv}\|_2^2.
\end{equation}
The profile bounds give
\[
 \begin{gathered}
 \nabla P_a^2\to\nabla Q^2,\qquad
 \Dtwo(yN_a)\to\Dtwo(-yQ^2)\quad\hbox{in }L^2,\\
 2M(a)\to\int Q^4.
 \end{gathered}
\]
For the polynomially decreasing test, use local convergence on
each fixed ball and the uniform exterior estimate
\[
 \|\mathbf1_{r>R}\Dtwo(yN_a)\|_2\le C/R,\qquad R\ge1.
\]
The exponentially decreasing test is handled by dominated
convergence.
Thus the velocity forms converge in $L^2$ operator norm:
each difference is a sum of products with one convergent
$L^2$ test difference. For sufficiently small $a$, their gap
is at least $1/(100c_0^2)$. The mixed term is absorbed using
\[
 Ca\|\epsn\|_2\|\epsv\|_2
       \le\frac1{200c_0^2}\|\epsv\|_2^2+C'a^2\|\epsn\|_2^2.
\]
Choose $a_*$ small enough to absorb the remaining density
coefficient into its limiting gap, together with the scalar
and density form-convergence errors. These choices precede
the choice of $I\Subset(0,a_*)$. Finally
\[
 \beta|\pair z{P_a}|^2
       \le C_I\beta(\beta+\|\veps\|_Y)^2\|\veps\|_Y^2
\]
is included in the allowed remainder by
\eqref{eq:nonlinear-moment-bound}. Consequently
\[
 \E_s\le-d_I\beta\|\veps\|_Y^2
       +C_I\beta(\beta+\|\veps\|_Y+\|\eps\|_2)\|\veps\|_Y^2,
\]
where the global mass factor comes only from
$|2(\lambda_s/\lambda)H_3|\le C_I\beta\|\eps\|_2\|\veps\|_Y^2$.
This proves \eqref{eq:energy-virial} directly for the actual
nonlinear errors.

Fix the constants, then choose the scale and error thresholds
so that $C_I(\beta_*+2\delta)\le d_I/2$.
Coercivity gives
\[
 \E_s\le-\tfrac12d_I\beta\|\veps\|_Y^2
       \le-\frac{d_Ic_0a_-}{2C_I}\lambda\E.
\]
The modulation bound \eqref{eq:rate-bound} gives the remaining
estimate in \eqref{eq:closed-energy}.
\end{proof}

\subsection{Extension to energy solutions}

The construction above was calculated on smooth solutions.
For energy data, its canceled scalar identities are the ones
that pass to the limit; no strong time derivative of the wave
errors is needed.

\begin{lemma}[Passage to energy solutions]\label{lem:energy-solutions}
Let the energy solution supplied by the local Cauchy theory
satisfy the modulation and bootstrap conditions on $[s_0,s_*)$.
On every compact subinterval, the tested modulation equations and
the canceled nonlinear energy identity hold in integrated form. In particular, after
the smallness choices of Proposition~\ref{prop:remainder},
\[
 \E(v)-\E(u)\le-\kappa_I\int_u^v\lambda(s)\E(s)\dd s,
 \qquad s_0\le u\le v<s_*.
\]
The energy is locally absolutely continuous, the modulation
parameters are $C^1$, and the corresponding differential
inequalities hold almost everywhere. The full-space conclusions
hold after the transport estimates of
Section~\ref{sec:centered-nonlinear} have been proved.
\end{lemma}

The proof is given in Appendix~\ref{app:energy-approximation}.
It first approximates exact solutions on a common physical
interval, then compares their own rescaled clocks. Only after
the cancellations does it pass the tested identities to the limit.

\section{Radial bootstrap closure and blow-up asymptotics}\label{sec:closure}

The inequality \eqref{eq:closed-energy} gives two kinds of
control. Its square-root integral bounds the total variation
of $a$ and closes the bootstrap. Once the interval is global,
comparison with the scale equation gives convergence rates.
The initial modulated values are $\mathcal E(s_0)$ and $\lambda(s_0)$.
The reference scale is $\lambda_0=\lambda^\star(0)$.

\subsection{Closing the bootstrap}

\begin{proof}[Proof of Proposition~\ref{prop:bootstrap} for radial data]
Fix any $S<s_*$. All estimates below are uniform in $S$.
Lemma~\ref{lem:energy-solutions} permits their use in integrated
form for energy solutions.

\paragraph{Step 1. Energy decay and variation of the profile parameter.}
On $[s_0,S]$, the bootstrap assumptions and
Proposition~\ref{prop:remainder} give $\E_s\le-\kappa_I\lambda\E$.
In particular $0\le\E(s)\le\mathcal E(s_0)$. For $\eta>0$, the absolutely continuous function
$\sqrt{\E+\eta}$ satisfies
\[
 (\sqrt{\E+\eta})_s\le
          -\frac{\kappa_I}{2}\frac{\lambda\E}{\sqrt{\E+\eta}}
          \qquad\text{almost everywhere}.
\]
Integrate and let $\eta\downarrow0$. The integrand is dominated by
$\lambda\sqrt\E$ on the compact interval. This gives, including
when $\E$ vanishes,
\begin{equation}\label{eq:sqrt-integral}
 \int_{s_0}^S\lambda\sqrt\E\dd s
 \le\frac2{\kappa_I}\bigl(\sqrt{\mathcal E(s_0)}-\sqrt{\E(S)}\bigr).
\end{equation}
The modulation bound $|a_s|\le K_I\lambda\sqrt\E$ therefore gives
\[
 \int_{s_0}^S|a_s|\dd s\le\frac{2K_I}{\kappa_I}\sqrt{\mathcal E(s_0)},
 \qquad |a(s)-a(s_0)|\le\int_{s_0}^S|a_s|\dd s.
\]
Consequently
\begin{equation}\label{eq:a-drift}
 \sup_{s_0\le s\le S}|a(s)-a(s_0)|
 +\int_{s_0}^S|a_s|\dd s\le C_I\sqrt{\mathcal E(s_0)}.
\end{equation}
\paragraph{Step 2. Profile overlap and propagation of the global mass error.}
Put $M_{\mathrm{tot}}=\norm{\psi_{\mathrm{in}}}_2^2$.
The two-dimensional rescaling preserves Schr\"odinger mass, and
$\Psi=Q_{a,\beta}+\eps$. Thus mass conservation gives the exact
identity
\begin{equation}\label{eq:mass-expansion}
 \norm{\eps(s)}_2^2=M_{\mathrm{tot}}-M(a(s))
                      -2\Rea\pair{\eps(s)}{Q_{a,\beta}}.
\end{equation}
To estimate the overlap, recall the consequence of the nonlinear
orthogonality condition proved in
\eqref{eq:nonlinear-moment-bound}:
\[
 |\pair\eps{P_a}|\le C_I\|\veps\|_Y^2.
\]
The difference $Q_{a,\beta}-P_a$ is an exponentially localized
test with weighted $H^1$ norm at most $C_I\beta$. Therefore
\[
 \begin{aligned}
 \pair\eps{Q_{a,\beta}}
       &=\pair\eps{P_a}+\pair\eps{Q_{a,\beta}-P_a},\\
 |\pair\eps{Q_{a,\beta}-P_a}|
       &\le C_I\beta\|\eps\|_{\mathrm{loc}},
 \end{aligned}
\]
and hence
\begin{equation}\label{eq:overlap-bound}
 \left|\pair\eps{Q_{a,\beta}}\right|
 \le C_I\bigl(\beta\|\veps\|_Y+\|\veps\|_Y^2\bigr).
\end{equation}
Only the local error norm is used in this overlap estimate.

Subtracting \eqref{eq:mass-expansion} at $s_0$ from the identity at
$s$ gives
\[
 \begin{aligned}
 \norm{\eps(s)}_2^2={}&\norm{\eps(s_0)}_2^2
       +M(a(s_0))-M(a(s))\\
 &+2\Rea\pair{\eps(s_0)}{Q_{a(s_0),\beta(s_0)}}
       -2\Rea\pair{\eps(s)}{Q_{a(s),\beta(s)}}.
 \end{aligned}
\]
The three differences on the right have the explicit bounds
\[
 \begin{aligned}
 |M(a(s))-M(a(s_0))|
 &\le\sup_{a\in I}|M'(a)|\,|a(s)-a(s_0)|
                                      \le C_I\sqrt{\mathcal E(s_0)},\\
 |\pair{\eps(s)}{Q_{a(s),\beta(s)}}|
 &\le C_I\bigl(\beta_*\sqrt{\E(s)}+\E(s)\bigr)
   \le C_I(\beta_*\sqrt{\mathcal E(s_0)}+\mathcal E(s_0)),\\
 |\pair{\eps(s_0)}{Q_{a(s_0),\beta(s_0)}}|
 &\le C_I(\beta_*\sqrt{\mathcal E(s_0)}+\mathcal E(s_0)).
 \end{aligned}
\]
The mass variation is controlled by \eqref{eq:a-drift};
the two overlaps use \eqref{eq:overlap-bound} and coercivity.
Since $\beta_*<1$, for
$\mathcal E(s_0)\le1$ we obtain
\begin{equation}\label{eq:global-mass}
 \sup_{s_0\le s\le S}\norm{\eps(s)}_2^2
 \le\norm{\eps(s_0)}_2^2+C_I\sqrt{\mathcal E(s_0)}.
\end{equation}
Together with energy decay, this estimate propagates and
improves the global mass bound on the same bootstrap interval.

\paragraph{Step 3. Fixed constants and initial margins.}
Fix $C_\star\ge1$ dominating the constants in the drift bound
\eqref{eq:a-drift}, the mass bound \eqref{eq:global-mass} and
$\norm{\veps}_Y\le C_\star\sqrt\E$. Keep $C_\star,K_I,\kappa_I$
fixed from now on. The required initial margins are
\begin{equation}\label{eq:bootstrap-margins}
 \begin{aligned}
 C_\star\sqrt{\mathcal E(s_0)}
   &<\min\{\delta/2,\delta^2/4,\delta_a\},&
 \norm{\eps(s_0)}_2&<\delta/2,\\
 K_I\sqrt{\mathcal E(s_0)}&<c_0a_-/2,&
 c_0a_+\lambda(s_0)&<\beta_*/2,\\
 |a(s_0)-a_0|&<\delta_a.
 \end{aligned}
\end{equation}
In particular these bounds imply $\mathcal E(s_0)<1$.
First take $c_0a_+\lambda^\star(0)<\beta_*/4$; continuity
preserves the scale margin for nearby data. Then shrink the
physical energy neighborhood. Continuity of the modulation and
energy makes the remaining errors and $\mathcal E(s_0)$ small, including
$|a(s_0)-a_0|$. These choices are uniform in $S$.

Energy decay, \eqref{eq:a-drift} and \eqref{eq:global-mass} now give
\[
 \norm{\veps(s)}_Y<\delta/2,\qquad
 |a(s)-a_0|\le |a(s_0)-a_0|+|a(s)-a(s_0)|<2\delta_a,
\]
and
\[
 \norm{\eps(s)}_2^2
 <\delta^2/4+\delta^2/4=\delta^2/2.
\]
Thus the global mass error improves to $\delta/\sqrt2$. Also
$\operatorname{dist}(a(s),\partial I)>2\delta_a$.
Only the scale and reference-phase improvement remain to be checked.

\paragraph{Step 4. The scale equation and the reference quadratic phase.}
The closed modulation estimate gives
\[
 \left|\frac{\lambda_s}{\lambda^2}+c_0a\right|
 \le K_I\sqrt\E\le K_I\sqrt{\mathcal E(s_0)}<c_0a_-/2.
\]
Consequently
\begin{equation}\label{eq:lambda-ode}
 \lambda_s=-[c_0a+O_I(\sqrt\E)]\lambda^2,
 \qquad 0<d_-\le(1/\lambda)_s\le d_+<\infty,
\end{equation}
where the fixed choices
\[
 d_-=c_0a_-/2,\qquad d_+=c_0a_++c_0a_-/2
\]
are valid. In particular $\lambda_s<0$. Integrating the bounds for
$(1/\lambda)_s$ and then inverting positive quantities yields
\begin{equation}\label{eq:scale-comparison}
 \frac{\lambda(s_0)}{1+d_+\lambda(s_0)(s-s_0)}
 \le\lambda(s)\le
 \frac{\lambda(s_0)}{1+d_-\lambda(s_0)(s-s_0)}.
\end{equation}
Since $a(s)\le a_+$ and $\lambda(s)\le\lambda(s_0)$, we have
\[
 0<\beta(s)=c_0a(s)\lambda(s)
                  \le c_0a_+\lambda(s_0)<\beta_*/2.
\]
This completes the strict improvement of
\eqref{eq:bootstrap-assumptions}.

\paragraph{Step 5. Exclusion of a finite endpoint.}
Suppose $s_*<\infty$. We continue the scale, physical solution
and modulation in that order, then use the strict margins to
contradict the definition of $s_*$.

First, \eqref{eq:scale-comparison} gives a positive lower bound
for $\lambda$ on $[s_0,s_*)$. Moreover
\[
 t_*:=\lim_{s\uparrow s_*}t(s)
      =\int_{s_0}^{s_*}\lambda(s)^2\dd s
      \le\lambda(s_0)^2(s_*-s_0)<\infty.
\]
The modulation rates bound $a_s$, $(\log\lambda)_s$ and
$\gamma_s$ on this finite interval. They therefore have finite
integrals and give limits for $a$, $\lambda$ and the chosen phase
lift. The limiting scale is strictly positive and the limiting
profile parameter stays a positive distance from $\partial I$.

Second, the rescaled profile bounds, coercivity and
\eqref{eq:global-mass} give bounded full rescaled energy norms.
In physical coordinates the relevant scaling identities are
\[
 \norm\psi_2^2=M_{\mathrm{tot}},\qquad
 \norm{\nabla\psi}_2=\lambda^{-1}\norm{\nabla\Psi}_2,
 \qquad (\norm n_2,\norm v_2)=\lambda^{-1}(\norm N_2,\norm V_2).
\]
Writing
$\lambda_{\min}=\lambda(s_0)/(1+d_+\lambda(s_0)(s_*-s_0))>0$, these
identities give, for a constant depending on the fixed initial data,
\[
 \sup_{0\le t<t_*}
  (\|\psi(t)\|_{H^1}+\|n(t)\|_2+\|v(t)\|_2)
 \le\sqrt{M_{\mathrm{tot}}}
       +\lambda_{\min}^{-1}C_I(1+\sqrt{\mathcal E(s_0)})<\infty.
\]
The local lifespan is bounded below in terms of this norm
by the energy-space contraction argument in
\cite[Theorem~1.1 and Section~5]{SanwalLWP}; see
Appendix~\ref{app:energy-approximation} for the equivalence
of the two wave formulations. Choose a common lower lifespan
$\tau_*>0$ for the displayed bound and restart at a time
$t_0<t_*$ with $t_*-t_0<\tau_*/2$. Local uniqueness joins the
new solution to the old one on their overlap, and its interval
extends beyond $t_*$.

Third, strong continuity of the continued solution and the
parameter limits give limiting rescaled errors. Lemma~\ref{lem:chart}
keeps the modulation Jacobian uniformly invertible: the reference
pivots are separated from zero on $I$, and the error lies in the
chosen uniform neighborhood. The implicit-function theorem extends
the modulation through $t_*$, agreeing with the preceding one by
local uniqueness.

Finally, the limiting local error is at most $\delta/2$, the
limiting global mass error is at most $\delta/\sqrt2$, the reference
quadratic phase coefficient is below $\beta_*/2$ and the profile parameter has distance
at least $2\delta_a$ from the boundary. Each is strictly inside
its original bootstrap range. Continuity preserves those original
strict conditions for a short interval beyond $t_*$, and the
positive scale extends the rescaled time beyond $s_*$. This
contradicts its defining maximality. Therefore $s_*=\infty$ and
proves the radial part of Proposition~\ref{prop:bootstrap}.
\end{proof}

The cited Cauchy theory applies to our velocity formulation
through the scalar wave variable $n_t$. With
$\norm h_{\dot H^{-1}}=\norm{\nabla(-\Delta)^{-1}h}_2$, a radial
no-swirl velocity satisfies
\[
 v=-\nabla\Delta^{-1}n_t,\qquad
 \norm{n_t}_{\dot H^{-1}}=\norm v_2.
\]
The reduction to unit sound speed is
\[
 (\psi,n,v)(t,x)
 =(c_0\psi_{\mathrm{unit}},c_0^2n_{\mathrm{unit}},
 c_0^3v_{\mathrm{unit}})(c_0^2t,c_0x).
\]
The corresponding energy norms are equivalent.

The orbital bound follows directly from coercivity, mass propagation
and parameter drift. For $\mathcal E(s_0)\le1$ these give
\[
 \sup_{s\ge s_0}
 \left(\norm{\eps}_{H^1}+\norm{\epsn}_2
       +c_0^{-1}\norm{\epsv}_2+|a(s)-a_0|\right)
 \le C_I\bigl(\norm{\eps(s_0)}_2+\mathcal E(s_0)^{1/4}+|a(s_0)-a_0|\bigr).
\]
Shrinking the initial neighborhood makes the right-hand side
smaller than the prescribed $\eta$ in Theorem~\ref{thm:main}.

\subsection{Convergence of the parameters and error}\label{sec:quantitative-rates}

The scale comparison makes $\lambda^2$ integrable. Define
\begin{equation}\label{eq:T}
 T=t(s_0)+\int_{s_0}^\infty\lambda(s)^2\dd s<\infty.
\end{equation}
The rate argument uses only the closed energy and modulation
bounds on $[s_0,\infty)$. It therefore applies to the nonradial
solution as soon as those same bounds are established.

\begin{proposition}[Rates from the closed energy inequality]
\label{prop:quantitative-rates}
Let $I$ be the fixed positive profile interval, and let $d_+$ be an
upper bound for $-\lambda_s/\lambda^2$ in \eqref{eq:lambda-ode}.
Set
\begin{equation}\label{eq:rate-exponent}
 \sigma_I=\frac{\kappa_I}{2d_+}>0.
\end{equation}
Then $a(s)$ converges to a limit $a_\infty$ in the interior of $I$,
and
\begin{equation}\label{eq:energy-power-rate}
 \E(s)\le\mathcal E(s_0)\left(\frac{\lambda(s)}{\lambda(s_0)}\right)^{2\sigma_I}.
\end{equation}
For either modulation system,
\begin{equation}\label{eq:all-power-rates}
\begin{aligned}
 &\norm{\eps(s)}_{\mathrm{loc}}+\norm{\epsn(s)}_2
       +c_0^{-1}\norm{\epsv(s)}_2\\
 &\quad+|a(s)-a_\infty|
       +\left|\frac{\lambda(s)}{T-t(s)}-c_0a_\infty\right|
 \le C_I\sqrt{\mathcal E(s_0)}
                      \left(\frac{\lambda(s)}{\lambda(s_0)}\right)^{\sigma_I}.
\end{aligned}
\end{equation}
For solutions with translation modulation, $x_c(s)$ has a limit $x_\infty$ and
\begin{equation}\label{eq:center-power-rate}
 \frac{|x_c(s)-x_\infty|}{T-t(s)}
 \le C_I\sqrt{\mathcal E(s_0)}
                      \left(\frac{\lambda(s)}{\lambda(s_0)}\right)^{\sigma_I}.
\end{equation}
The constants and exponent are uniform on the chosen neighborhood.
\end{proposition}

\begin{proof}
\paragraph{Decay of the energy.}
Since $-\lambda_s/\lambda\le d_+\lambda$,
\eqref{eq:closed-energy}, or its full-space version, gives directly
\[
 (\E\lambda^{-2\sigma_I})_s
 =\lambda^{-2\sigma_I}(\E_s-2\sigma_I(\lambda_s/\lambda)\E)
 \le(-\kappa_I+2\sigma_Id_+)\lambda^{1-2\sigma_I}\E=0.
\]
This argument also covers zeros of $\E$ without dividing by it.
For energy solutions the product rule holds almost everywhere by
the absolute continuity in Lemma~\ref{lem:energy-solutions}.
Integration gives \eqref{eq:energy-power-rate}; since
$\lambda(s)\to0$, it also gives $\E(s)\to0$.

\paragraph{Limits of the parameters.}
Apply the square-root estimate \eqref{eq:sqrt-integral} from time
$s$ to infinity:
\begin{equation}\label{eq:tail-square-root}
 \int_s^\infty\lambda(\tau)\sqrt{\E(\tau)}\dd\tau
                                  \le\frac2{\kappa_I}\sqrt{\E(s)}.
\end{equation}
The parameter rate makes $a_s$ integrable, so $a$ has a limit
$a_\infty$. It lies in the interior of $I$, since
\[
 \operatorname{dist}(a(s),\partial I)\ge2\delta_a.
\]
Apply the same rate estimate on $[s,\infty)$ to obtain
\[
 |a(s)-a_\infty|\le\frac{2K_I}{\kappa_I}\sqrt{\E(s)}.
\]
Moreover
\[
 \lambda_t=-c_0a+\frac{\frac{\lambda_s}{\lambda}+\beta}{\lambda},\qquad
 |\lambda_t+c_0a_\infty|\le C_I\sqrt{\E(s(t))}.
\]
Integrate up to $T$, using $\lambda(T)=0$. Monotonicity of $\E$
then gives
\[
 |\lambda(t)-c_0a_\infty(T-t)|
       \le C_I(T-t)\sqrt{\E(s(t))}.
\]
Coercivity controls the three error norms by the same square root. This proves
\eqref{eq:all-power-rates}, and in particular
\begin{equation}\label{eq:speed-limit}
 \frac{\lambda(t)}{T-t}\longrightarrow c_0a_\infty.
\end{equation}

For solutions with translation modulation, $|x_{c,s}|=\lambda|\nu|
\le K_I\lambda^2\sqrt\E$ is integrable. It defines a limit
$x_\infty$ and gives
\begin{equation}\label{eq:center-convergence-quantitative}
 |x_c(s)-x_\infty|
 \le K_I\int_s^\infty\lambda(\tau)^2\sqrt{\E(\tau)}\dd\tau
 \le K_I(T-t(s))\sqrt{\E(s)}.
\end{equation}
Together with \eqref{eq:energy-power-rate}, this proves
\eqref{eq:center-power-rate}.

\paragraph{Conversion to physical time.}
The two scale bounds give
\begin{equation}\label{eq:remaining-time-comparison}
 \frac{\lambda(s)}{d_+}\le T-t(s)
                           \le\frac{\lambda(s)}{d_-}.
\end{equation}
Indeed,
\[
 \frac1{d_+}\int_s^\infty(-\lambda_\tau)\dd\tau
 \le\int_s^\infty\lambda(\tau)^2\dd\tau
 \le\frac1{d_-}\int_s^\infty(-\lambda_\tau)\dd\tau,
\]
and the outer integrals are $\lambda(s)/d_+$ and
$\lambda(s)/d_-$ because $\lambda(\tau)\to0$.
Thus the norm and parameter bounds are $O((T-t)^{\sigma_I})$,
the scale remainder is $O((T-t)^{1+\sigma_I})$, and so is the
center displacement. The scale comparison also gives the equivalent
error rate $O((1+s-s_0)^{-\sigma_I})$ in rescaled time. The constants
in these last $O$ bounds may depend on the fixed reference scale,
as permitted in Theorem~\ref{thm:main}.
\end{proof}

\begin{corollary}[Convergence of the radiation mass]\label{cor:radiation-mass}
Under the hypotheses of Proposition~\ref{prop:quantitative-rates},
\begin{equation}\label{eq:radiation-mass-limit}
 \lim_{s\to\infty}\norm{\eps(s)}_2^2
 =M_{\mathrm{tot}}-M(a_\infty)\ge0,
\end{equation}
and
\begin{equation}\label{eq:radiation-mass-rate}
 \left|\norm{\eps(s)}_2^2-
                 [M_{\mathrm{tot}}-M(a_\infty)]\right|
 \le C_I\sqrt{\mathcal E(s_0)}
                      \left(\frac{\lambda(s)}{\lambda(s_0)}\right)^{\sigma_I}.
\end{equation}
\end{corollary}

\begin{proof}
The overlap estimate \eqref{eq:overlap-bound} and coercivity give
\[
 |\Rea\pair\eps{Q_{a,\beta}}|
                         \le C_I(\beta\sqrt\E+\E).
\]
For the radiation mass, subtraction of its candidate limit from
\eqref{eq:mass-expansion} gives the exact equality
\[
 \|\eps\|_2^2-[M_{\mathrm{tot}}-M(a_\infty)]
 =M(a_\infty)-M(a)-2\Rea\pair\eps{Q_{a,\beta}}.
\]
Consequently
\[
 \begin{aligned}
 \left|\norm\eps_2^2-[M_{\mathrm{tot}}-M(a_\infty)]\right|
 &\le |M(a)-M(a_\infty)|+2|\Rea\pair\eps{Q_{a,\beta}}|\\
 &\le C_I\sqrt\E.
 \end{aligned}
\]
Here $M'$ is bounded on $I$, $|a-a_\infty|\le C_I\sqrt\E$,
and $\beta$ and $\E$ are small. Applying
\eqref{eq:energy-power-rate} proves \eqref{eq:radiation-mass-rate}.
The limit follows as $\E\to0$; its nonnegativity follows because
it is a limit of squared norms.
\end{proof}

\paragraph{Nonzero radiation mass arbitrarily near the exact solution.}\label{par:radiation-example}
We justify the assertion in Remark~\ref{rem:topology} without
assuming a scattering description of the radiation. Fix a reference
initial datum with parameters $a_0,\lambda_0,\gamma_0,x_0$, and
use the normalized ground state $\phi(N_{a_0})=P_{a_0}/\norm{P_{a_0}}_2$.
Choose a real radial $\chi\in C_c^\infty(\{1<|y|<2\})$ with
$\norm\chi_2=1$. For a fixed sufficiently small number $d>0$ set
\[
 g_R(y)=\frac dR\chi(y/R),\qquad
 h_R=g_R-(g_R,\phi(N_{a_0}))\phi(N_{a_0}).
\]
These symbols are used only in this construction. As $R\to\infty$,
\[
 (h_R,\phi(N_{a_0}))=0,\qquad \norm{h_R}_{\mathrm{loc}}\to0,
 \qquad \norm{h_R}_2\to d.
\]
Indeed, rescaling and the annular support give
\[
 \begin{gathered}
 \|g_R\|_2=d,\qquad \|\nabla g_R\|_2=\frac dR\|\nabla\chi\|_2,
 \qquad \int e^{-r}|g_R|^2\le d^2e^{-R},\\
 |(g_R,\phi(N_{a_0}))|\le d\|\mathbf1_{\{R<r<2R\}}\phi(N_{a_0})\|_2\to0.
 \end{gathered}
\]
The subtracted projection is exponentially localized, so its local
and gradient norms also tend to zero. Take the
rescaled initial state
\[
 (\Psi_R,N_R,V_R)
       =(Q_{a_0,c_0a_0\lambda_0}+h_R,N_{a_0},V_{a_0}).
\]
The projection and density conditions \eqref{eq:chart} hold
exactly at $a_0,\lambda_0,\gamma_0$; the translation conditions hold
by radiality. In physical variables the field perturbation has
$L^2$ norm tending to $d$ and gradient norm
$\lambda_0^{-1}\norm{\nabla h_R}_2\to0$. Thus $d$ can be
chosen so that these data lie in any prescribed stability
neighborhood for all sufficiently large $R$.

By coercivity the initial nonlinear energy of these data tends
to zero. The drift estimate
\eqref{eq:a-drift} consequently implies
$a_{\infty,R}\to a_0$. Since
$\pair{h_R}{Q_{a_0,c_0a_0\lambda_0}}\to0$, the conserved masses satisfy
$M_{\mathrm{tot},R}\to M(a_0)+d^2$. It follows that
\[
 M_{\mathrm{tot},R}-M(a_{\infty,R})\longrightarrow d^2>0.
\]
For large $R$ the limit in \eqref{eq:radiation-mass-limit} is
therefore strictly positive. On the other hand, the local
convergence forces any strong global $L^2$ limit of the rescaled
error to be zero. This proves the claimed failure of global
strong $L^2$ attraction on an open energy neighborhood.

\paragraph{Conclusion for radial perturbations.}
The bootstrap gives the full orbital bound, and
Proposition~\ref{prop:quantitative-rates} gives all the stated limits
and rates. Finally,
\[
 \norm{\nabla\Psi}_2\ge\norm{\nabla Q_{a,\beta}}_2
                         -\norm{\nabla\eps}_2
\]
has a positive lower bound after shrinking $\mathcal E(s_0)$. Indeed, since
$P_a$ is real,
\[
 \norm{\nabla Q_{a,\beta}}_2^2
 =\norm{\nabla P_a}_2^2+\frac{\beta^2}{4}\norm{yP_a}_2^2,
 \qquad \inf_{a\in I}\norm{\nabla P_a}_2>0.
\]
Coercivity makes $\norm{\nabla\eps}_2$ smaller than half this
positive infimum. Therefore
\[
 \norm{\nabla\psi(t)}_2
     =\lambda^{-1}\norm{\nabla\Psi}_2\longrightarrow\infty
                    \qquad(t\uparrow T).
\]
Thus $T$ is the maximal physical existence time.

\section{Stability under nonradial perturbations}\label{sec:nonradial}

We keep the radial profiles and the energy \eqref{eq:B} for general
perturbations in $X$ and $Y$. The two density orthogonality conditions
in the translation directions determine the center. Its motion
contributes to the principal velocity form. We derive the translation
rates, prove coercivity under the modulation constraints, and compute
the work from translations in the nonlinear energy balance.
The local operator estimates of Lemma~\ref{lem:anchoring} already hold
without symmetry.

\begin{proposition}[Schr\"odinger and velocity coercivity]
\label{prop:angular-scalar-velocity}
The energy estimate \eqref{eq:scalar-energy-gap} holds for all
$h\in H^1(\R^2;\C)$ with $\pair hQ=0$, without radial symmetry.
The Schr\"odinger coercivity estimate \eqref{eq:spectral-scalar}
also extends to that space,
with a possibly smaller positive constant. Both Schr\"odinger estimates
extend by density to $X$ without radial symmetry.
The velocity estimate \eqref{eq:spectral-velocity} holds for
every real vector field $w\in L^2(\R^2;\R^2)$.
\end{proposition}

The angular mean uses the radial Schr\"odinger estimates. On the
nonzero modes, each potential $Q^2$ and $-rQQ'$ costs at most
$\|Q\|_2^2/(4\pi)<1$ times the gradient square by
\eqref{eq:angular-potentials}. For the velocity form, its two
tests are radial, so the component orthogonal to radial no-swirl
fields contributes just its squared norm. The details are in
Appendix~\ref{app:angular-comparison}.

\subsection{Translation modulation and transport}

In coordinates $y=(x-x_c(t))/\lambda(t)$, the center contributes
transport $\nu\cdot\nabla$ to each rescaled equation, where
$\nu=x_{c,s}/\lambda\in\R^2$. We choose $x_c$ through two density
orthogonality conditions:
\begin{equation}\label{eq:center-conditions}
 (\epsn,\partial_jN_a)=0,\qquad j=1,2.
\end{equation}

\begin{proposition}[Translation modulation]
\label{prop:center-modulation}
Near a fixed sufficiently late point of the radial orbit, the
conditions \eqref{eq:chart} and \eqref{eq:center-conditions}
define $(a,\lambda,\gamma,x_c)$ for general energy
data. The error orthogonality acquires a localized translation term,
while its additional energy work is computed in
Section~\ref{sec:centered-nonlinear}. With all fields now unrestricted
in angle, the modulation bounds include
\begin{equation}\label{eq:center-rate}
 |a_s|+|\frac{\lambda_s}{\lambda}+\beta|+|\nu|
 \le C_I\lambda
       (\norm{\eps}_{\mathrm{loc}}+\norm{\epsn}_2
                                      +c_0^{-1}\norm{\epsv}_2).
\end{equation}
The phase bounds \eqref{eq:phase-bound} and \eqref{eq:phase-equation}
also hold with the full-space norm $\|\veps\|_Y$.
\end{proposition}

\begin{proof}
In physical coordinates the two added constraints are
\[
 \int n(x)\,\partial_jN_a\left(\frac{x-x_c}{\lambda}\right)\dd x=0,
 \qquad j=1,2,
\]
because $(N_a,\partial_jN_a)=0$. Parameter derivatives therefore act
on smooth profile tests, even for rough density data.
In the translation columns use $\lambda\partial_{x_{c,k}}$ at
fixed physical data. The corresponding density row has derivative
$-(N,\partial_k\partial_jN_a)$, so at the reference the block is
\[
 -(N_a,\partial_k\partial_jN_a)
   =(\partial_kN_a,\partial_jN_a)
   =\|\partial_1N_a\|_2^2\delta_{jk}.
\]
The diagonal coefficient has a positive lower bound on $I$, since
$N_a$ is nonconstant and depends continuously on $a$ in $H^1$.
The cross entries with the radial scale and profile directions
vanish by parity. In the complex constraint
\eqref{eq:physical-chart-tests}, the eigenfunction at the physical
density is fixed when taking parameter derivatives. Translating the
profile test gives $\pair{\partial_jQ_{a,\beta}}{\phi(N_a)}=0$
at the reference. The remaining pivots are those of
Lemma~\ref{lem:chart}. The physical tests are continuously
differentiable in the parameters, so the implicit-function
argument applies to energy data.

Translation covariance gives
$\phi'_N[\partial_jN]=\partial_j\phi(N)$. Thus the two error
transport terms in the derivative of \eqref{eq:chart} cancel:
\[
 \nu\cdot\bigl(\pair{\nabla\eps}{\phi(N)}
                       +\pair\eps{\nabla\phi(N)}\bigr)=0.
\]
The profile transport leaves $\sum_j\nu_j\pair{\partial_jQ_{a,\beta}}\phi$.
The identity holds first for smooth potentials; approximation and
Lemma~\ref{lem:spectral} extend the tested expression to energy
regularity. With $\phi=\phi(N)$, the complex modulation row is
\begin{equation}\label{eq:centered-error-projection}
\begin{aligned}
 a_s\pair{\partial_a Q_{a,\beta}+(\beta/a)\partial_\beta Q_{a,\beta}}\phi
 &-\left(\frac{\lambda_s}{\lambda}+\beta\right)
          \pair{\Lambda Q_{a,\beta}-\beta\partial_\beta Q_{a,\beta}}\phi
       -\sum_{j=1}^2\nu_j\pair{\partial_jQ_{a,\beta}}\phi\\
 &+\ii(\mu(N)+\gamma_s)\pair{Q_{a,\beta}}\phi
   =\ii\pair{L_aQ_{a,\beta}}\phi-\lambda\pair\eps{\phi'_N[\diver V]}.
\end{aligned}
\end{equation}
Radial parity and the eigenfunction bound give
\[
 |\pair{\partial_jQ_{a,\beta}}\phi|
   =|\pair{\partial_jQ_{a,\beta}}{\phi-\phi(N_a)}|
                         \le C_I\|\epsn\|_2.
\]
The eigenpair estimates use only two-dimensional Sobolev products
and the simple isolated eigenvalue, so they apply to nonradial
$L^2$ potentials near $N_a$ as well.

The exact density error equation is
\[
 \begin{aligned}
 (\epsn)_s={}&\left(\frac{\lambda_s}{\lambda}+\beta\right)\Dtwo N_a
       +\frac{\lambda_s}{\lambda}\Dtwo {\epsn}-\lambda\diver {\epsv}
       -a_s\partial_aN_a+\nu\cdot\nabla(N_a+{\epsn}).
 \end{aligned}
\]
The differentiated constraints are
\[
 \begin{aligned}
 0&=((\epsn)_s,P_a^2)+2a_s({\epsn},P_a\partial_aP_a),\\
 0&=((\epsn)_s,\partial_jN_a)+a_s({\epsn},\partial_a\partial_jN_a).
 \end{aligned}
\]
The profile equation and its derivatives give
\[
\begin{aligned}
 L_a\Lambda P_a&=-2P_a-P_a\Dtwo N_a,&
 L_a\partial_aP_a&=-P_a\partial_aN_a,\\
 (\Dtwo N_a,P_a^2)&=-2M,&
 (\partial_aN_a,P_a^2)&=0.
\end{aligned}
\]
The second line follows by pairing the first with $P_a$ and
using $L_aP_a=0$. Radial parity gives
\[
 (\nabla N_a,P_a^2)=(\Dtwo N_a,\partial_jN_a)
                         =(\partial_aN_a,\partial_jN_a)=0.
\]
For rough ${\epsn}$, integration by parts means
$(\partial_\ell {\epsn},\partial_jN_a)=-({\epsn},\partial_\ell\partial_jN_a)$.
Substitution yields the density row
\begin{equation}\label{eq:centered-density-row}
\begin{aligned}
 &[-2M+(\Dtwo {\epsn},P_a^2)](\frac{\lambda_s}{\lambda}+\beta)
       +2a_s({\epsn},P_a\partial_aP_a)+\nu\cdot(\nabla {\epsn},P_a^2)\\
 &\qquad=\lambda(\diver {\epsv},P_a^2)+\beta(\Dtwo {\epsn},P_a^2).
\end{aligned}
\end{equation}
The coefficient of $\nu$ vanishes at zero error. The translation rows are
\begin{equation}\label{eq:center-system}
\begin{aligned}
 &({\epsn},\partial_a\partial_jN_a)a_s
       +(\Dtwo {\epsn},\partial_jN_a)(\frac{\lambda_s}{\lambda}+\beta)\\
 &\quad+\sum_{\ell=1}^2
 [(\partial_\ell N_a,\partial_jN_a)
            -({\epsn},\partial_\ell\partial_jN_a)]\nu_\ell\\
 &\qquad=\lambda(\diver {\epsv},\partial_jN_a)
                         +\beta(\Dtwo {\epsn},\partial_jN_a).
\end{aligned}
\end{equation}
Taking the real part of \eqref{eq:centered-error-projection}
after division by $\pair{Q_{a,\beta}}\phi$, together with the density and translation rows,
gives a system for $(a_s,\lambda_s/\lambda+\beta,\nu)$.
At zero error its translation block is
$(\partial_\ell N_a,\partial_jN_a)
=\norm{\partial_1N_a}_2^2\delta_{j\ell}$; radial parity removes
the off-diagonal blocks. Its scalar blocks are the radial modulation
matrix above. All these inverses are bounded on $I$.

For the perturbed translation block,
\[
 |({\epsn},\partial_\ell\partial_jN_a)|
       +|({\epsn},\partial_a\partial_jN_a)|
       +|(\Dtwo {\epsn},\partial_jN_a)|\le C_I\norm {\epsn}_2,
\]
where $\Dtwo$ is transposed in the last pairing. The density-row
coefficient of $\nu$ has the same bound. Thus the full matrix
perturbation has norm at most $C_I\norm{\veps}_Y$.
The transposed velocity pairings are bounded by $C_I\norm {\epsv}_2$.
The forcing of the scalar row vanishes at zero error, since
$\pair{L_aQ_{a,\beta}}{\phi(N_a)}=0$ and $\eps=0$. Its translation coefficient is
$O_I(\|\epsn\|_2)$ by the bound above. Applying the zero-error
inverse and the same forcing estimates as in
Lemma~\ref{lem:modulation} gives
\[
 \begin{aligned}
 |a_s|+|\lambda_s/\lambda+\beta|+|\nu|
 &\le C_I\lambda\norm{\veps}_Y\\
 &\quad+C_I\norm{\veps}_Y
                   (|a_s|+|\lambda_s/\lambda+\beta|+|\nu|).
 \end{aligned}
\]
Choose $C_I\delta<1/2$ and absorb the last line. This proves
\eqref{eq:center-rate}, which also gives \eqref{eq:beta-rate}
for the modulation system with translations.
The tested coefficients and the bounded
matrix inverse also give continuous dependence of the rates,
as required by Lemma~\ref{lem:energy-solutions}. Divide
\eqref{eq:centered-error-projection} by $\pair{Q_{a,\beta}}\phi$
and take imaginary parts. The new term has size at most
$C_I|\nu|\|\epsn\|_2\le C_I\lambda\|\veps\|_Y^2$, so
$|\gamma_s+\mu(N)|\le C_I\lambda\|\veps\|_Y$ as in
\eqref{eq:phase-equation}. The expansion of $\mu(N)+1$ on
the density constraint then gives \eqref{eq:phase-bound}.
Appendix~\ref{app:energy-approximation} passes these tested
identities to energy solutions.
\end{proof}

To state the velocity form independently of the solution, let
$w\in L^2(\R^2;\R^2)$ be an arbitrary real test field.
After first taking the leading $\beta$ coefficient at fixed $a$ and
then letting $a\downarrow0$, the velocity form including translation terms is
\begin{equation}\label{eq:center-velocity-bracket}
\begin{aligned}
 &\norm {w}_2^2+
 \frac{({w},\nabla Q^2)({w},\Dtwo(-yQ^2))}{\int Q^4}\\
 &\qquad-\sum_{j=1}^2\frac{({w},\nabla\partial_j(Q^2))({w},\partial_j(yQ^2))}{\|\partial_j(Q^2)\|_2^2}.
\end{aligned}
\end{equation}
Indeed, the velocity part of the linear translation rate is
$\lambda(\diver {w},\partial_jN_a)/\norm{\partial_jN_a}_2^2$.
Its energy contribution pairs this rate with
$c_0^{-2}({w},\partial_jV_a)$. Since $V_a=-c_0a\,yN_a$
and $N_a\to-Q^2$, division by the negative dissipation
coefficient $-\beta$ gives the minus sign in
\eqref{eq:center-velocity-bracket}. The common factor
$c_0^{-2}$ has been omitted there. The corresponding density
test $\partial_j(N_a+P_a^2)$ has $L^2$ norm $O(a^2)$, by
$\nabla(N_a+P_a^2)=a^2\Dtwo(yN_a)$ and the uniform profile tails.

The two translation tests satisfy
\[
 \begin{aligned}
 (\nabla\partial_j(Q^2),\partial_j(yQ^2))
 & =-(\partial_j(Q^2),\partial_j\diver(yQ^2))\\
 &=-\bigl(\partial_j(Q^2),(3+y\cdot\nabla)\partial_j(Q^2)\bigr)
   =-2\|\partial_j(Q^2)\|_2^2.
 \end{aligned}
\]
The last equality is integration by parts in two dimensions.
Cross pairings for different indices vanish by parity. Each pair of
translation tests transforms in angular order one under simultaneous rotation
of its argument and vector values. The translation tests are orthogonal to
the two invariant radial vectors in the scale-modulation term. The rank-two
eigenvalue formula on each translation plane consequently shows that
strict positivity follows from the additional inequalities
\begin{equation}\label{eq:center-velocity-product}
 \norm{\nabla\partial_j(Q^2)}_2\norm{\partial_j(yQ^2)}_2<4\|\partial_j(Q^2)\|_2^2,\qquad j=1,2.
\end{equation}
The moment estimates establishing these inequalities are given
in Appendix~\ref{app:center-velocity-proof}.

\begin{proposition}[Velocity coercivity with translation terms]\label{prop:center-velocity-gap}
For the Townes profile, the bracket
\eqref{eq:center-velocity-bracket} is bounded below by
$\frac1{50}\norm {w}_2^2$ on all real $L^2$ vector fields.
\end{proposition}

The verification is given in Appendix~\ref{app:center-velocity-proof}.

\subsection{Density coercivity under translation constraints}\label{sec:density-coercivity}

On the full space $L^2(\R^2)$ define
\[
 R_Q=(L_-|_{Q^\perp})^{-1}\Pi_{Q^\perp},\qquad
 \mathcal B_Q=1-2QR_QQ .
\]
The inverse removes the one-dimensional kernel $\C Q$, and each
adjacent $Q$ denotes multiplication. The density energy is
$(\mathcal B_Qf,f)/2$ for a real density test $f\in L^2(\R^2)$.
Its radial coercivity is proved in
Lemma~\ref{lem:limiting-energy}. For the derivative, define the
bounded self-adjoint operator $\mathcal D_Q$ by
\begin{equation}\label{eq:full-density-form}
 (\mathcal D_Q {f},{f})=\norm {f}_2^2+(A_Qz,z)
                  -2({f},(Q-rQ')z),\qquad z=R_Q(Q{f}).
\end{equation}
This is the limiting pure-density form in the energy derivative;
$\mathcal B_Q$ instead occurs in the energy itself. Their positivity
statements serve different purposes. The radial restriction of
$\mathcal D_Q$ agrees with the operator in
Appendix~\ref{app:certificate}.

\begin{proposition}[Coercivity of the density energy under translation constraints]
\label{prop:centered-density-energy}
There is $c>0$ such that
\begin{equation}\label{eq:centered-density-energy}
 (\mathcal B_Q{f},{f})\ge c\norm {f}_2^2
\end{equation}
for real ${f}\in L^2(\R^2)$ satisfying
$({f},Q^2)=({f},Q\partial_1Q)=({f},Q\partial_2Q)=0$.
Together with the Schr\"odinger energy gap, this controls the density
part of the limiting nonlinear energy under the translation conditions.
\end{proposition}

The proof is given in Appendix~\ref{app:angular-density}.

The first angular mode requires the translation constraint.
The following identity explains its effect on the density
dissipation form.

\begin{lemma}[Translation identity for the density form]
\label{prop:density-translation}
For $j=1,2$,
\begin{equation}\label{eq:density-translation-null}
 \mathcal B_Q(Q\partial_jQ)=0,\qquad (\mathcal D_Q(Q\partial_jQ),(Q\partial_jQ))=0.
\end{equation}
\end{lemma}

\begin{proof}
Differentiating the ground-state equation gives
$L_-\partial_jQ=2Q^2\partial_jQ$. Hence
\[
 R_Q(Q^2\partial_jQ)=\tfrac12\partial_jQ,\qquad
 \mathcal B_Q(Q\partial_jQ)=0.
\]
The ground-state projection annihilates the first angular
subspace. The commutator $[L_-,\Lambda]=A_Q$ gives
$[\Lambda,R_Q]=R_QA_QR_Q$. Expanding
$\mathcal B_Q=1-2QR_QQ$ and using
$[\Lambda,Q]=rQ'$ therefore yields
\begin{equation}\label{eq:density-energy-commutator}
 \mathcal D_Q=\mathcal B_Q-\tfrac12[\Lambda,\mathcal B_Q].
\end{equation}
This identity is valid on the smooth exponentially decreasing
functions used here. Indeed $[\Lambda,\mathcal B_Q](Q\partial_jQ)=-\mathcal B_Q\Lambda (Q\partial_jQ)$, so
\[
\begin{aligned}
 (\mathcal D_Q(Q\partial_jQ),Q\partial_jQ)
 &=\tfrac12(\mathcal B_Q\Lambda(Q\partial_jQ),Q\partial_jQ)\\
 &=\tfrac12(\Lambda(Q\partial_jQ),\mathcal B_Q(Q\partial_jQ))=0.
\end{aligned}
\]
All functions in this calculation are smooth and exponentially
localized, so the commutator is defined by integration by parts.
\end{proof}

In each first angular amplitude, let $e$ be the normalization of
$QQ'$ in $L^2(r\dd r)$. Appendix~\ref{app:projected-density-criterion}
identifies $\mathcal D_Q-1$ on this mode with the integral operator
$K_1$ and gives its squared Hilbert--Schmidt norm $I_1$.
The translation identity gives $(K_1e,e)=-1$. The block on the
span of $e$ therefore contributes one to $I_1$, leaving at most
$I_1-1$ for the squared norm of the compression to $e^\perp$.
The interval bound $I_1<5/3$ in Appendix~\ref{sec:angular-certificate}
then gives a gap at least $1-\sqrt{I_1-1}>1/6$ under the
translation constraints. The next corollary combines this bound
with the radial and higher-mode estimates.

\begin{corollary}[Full-space coercivity of the limiting density dissipation form]
\label{cor:full-centered-density}
For every real ${f}\in L^2(\R^2)$ with
$({f},Q\partial_1Q)=({f},Q\partial_2Q)=0$,
\begin{equation}\label{eq:full-centered-density-gap}
 (\mathcal D_Q{f},{f})\ge\frac1{25}\norm {f}_2^2.
\end{equation}
Moreover, for all sufficiently small positive $a$, the moving
constraints $({f},\partial_jN_a)=0$, $j=1,2$, imply
\begin{equation}\label{eq:moving-centered-density-gap}
 (\mathcal D_Q{f},{f})\ge\frac1{50}\norm {f}_2^2.
\end{equation}
Both estimates concern the limiting operator $\mathcal D_Q$.
\end{corollary}

\begin{proof}
Decompose ${f}$ into its radial part and real sine and cosine modes.
The interval estimate in
Appendix~\ref{sec:radial-density-certificate} for the radial density
kernel gives
$\mathcal D_Q\ge1/20$ on the full radial $L^2$ subspace.
The preceding compression gives the $1/6$ gap in angular order one
under the translation constraints. All modes of order at least two
have gap $1/25$ by Proposition~\ref{prop:high-angular-density}.
Rotation invariance and
Parseval's identity give \eqref{eq:full-centered-density-gap}, first
on finite sums and then on $L^2$ by boundedness.

For the moving conditions, the profile convergence and the stationary identity
$\nabla(N_a+P_a^2)=a^2\Dtwo(yN_a)$ imply
\[
 \norm{-\tfrac12\partial_jN_a-Q\partial_jQ}_2\longrightarrow0.
\]
Indeed the right-hand side of that stationary identity tends to zero
in $L^2$ by the uniform profile tails, and
$\nabla(P_a^2)\to\nabla(Q^2)$ by local convergence and the uniform
exponential tails of $P_a$. Since the two nonzero functions
$Q\partial_jQ$ are orthogonal, define
\[
 {f}^{\,0}={f}-\sum_{j=1}^2
       \frac{({f},Q\partial_jQ)}{\norm{Q\partial_jQ}_2^2}Q\partial_jQ.
\]
If $({f},\partial_jN_a)=0$, then
\[
 |({f},Q\partial_jQ)|
 =|({f},Q\partial_jQ+\tfrac12\partial_jN_a)|=o_a(1)\|{f}\|_2,
 \qquad \|{f}-{f}^{\,0}\|_2=o_a(1)\|{f}\|_2.
\]
Boundedness of $\mathcal D_Q$ and \eqref{eq:full-centered-density-gap} give
\[
 \begin{aligned}
 (\mathcal D_Q{f},{f})
 &\ge\tfrac1{25}\|{f}^{\,0}\|_2^2
       -\|\mathcal D_Q\|\|{f}-{f}^{\,0}\|_2(\|{f}\|_2+\|{f}^{\,0}\|_2)\\
 &\ge\bigl(\tfrac1{25}-o_a(1)\bigr)\|{f}\|_2^2.
 \end{aligned}
\]
Taking $a$ sufficiently small proves
\eqref{eq:moving-centered-density-gap}.
\end{proof}

\subsection{Full-space coercivity of the energy and dissipation forms}
\label{sec:centered-linear}

The preceding estimates give separate Schr\"odinger, density and
velocity gaps. The motion of the center adds a wave pairing to
\eqref{eq:master}. We include it in the principal dissipation form
before absorbing the nonlinear errors. The translation directions
are the spatial derivatives of the reference state:
\begin{equation}\label{eq:center-tangent}
 \partial_j\U_*=(\partial_jQ_{a,\beta},\partial_jN_a,\partial_jV_a).
\end{equation}
For a real density test $f\in L^2(\R^2)$ and a real velocity test
$w\in L^2(\R^2;\R^2)$, the leading translation rate
$\nu^{\mathrm{lin}}(f,w)=(\nu_1^{\mathrm{lin}}(f,w),\nu_2^{\mathrm{lin}}(f,w))$
has components
\begin{equation}\label{eq:linear-center-output}
 \nu_j^{\mathrm{lin}}(f,w)
 =\frac{\beta(\Dtwo {f},\partial_jN_a)+\lambda(\diver {w},\partial_jN_a)}
             {\|\partial_jN_a\|_2^2}.
\end{equation}
Its size is at most $C_I\beta(\|f\|_2+c_0^{-1}\|w\|_2)$. Since $\lambda/\beta=(c_0a)^{-1}$,
the quotients $\nu_j^{\mathrm{lin}}(f,w)/\beta$ depend on $a$ at fixed
$c_0$ and are independent of $\beta$. In the forms below,
$h$ is a Schr\"odinger test in the corrected coordinate. For a solution
we take $h=z$. Define
\begin{equation}\label{eq:centered-leading-form}
\begin{aligned}
 \mathcal D_a^c(h,f,w)={}&(A_ah,h)_\R
                         +D_a({f})+c_0^{-2}\|{w}\|_2^2\\
 &-\frac{({f},y\cdot\nabla P_a^2)+(c_0a)^{-1}({w},\nabla P_a^2)}{2M}
       \bigl[({f},\Dtwo(N_a+P_a^2))\\
 &\hspace{57mm}+c_0^{-2}({w},\Dtwo V_a)\bigr]\\
 &-\sum_j\frac{\nu_j^{\mathrm{lin}}(f,w)}{\beta}
       [({f},\partial_j(N_a+P_a^2))+c_0^{-2}({w},\partial_jV_a)].
\end{aligned}
\end{equation}
\begin{proposition}[Full-space coercivity]\label{prop:centered-linear}
For sufficiently small $a_*>0$, let $w\in L^2(\R^2;\R^2)$ be arbitrary
and let $f\in L^2(\R^2;\R)$ satisfy
\[
 ({f},P_a^2)=0,\qquad ({f},\partial_jN_a)=0,\quad j=1,2.
\]
For every Schr\"odinger test $h\in X$, the energy form satisfies
\begin{equation}\label{eq:centered-energy-form}
\begin{aligned}
 &((L_a-1)h,h)_\R+\tfrac12({f},(1-2P_aR_aP_a){f})
                                  +\tfrac1{2c_0^2}\|{w}\|_2^2\\
 &\qquad\ge c_I\bigl(\|h\|_{\mathrm{loc}}^2+\|{f}\|_2^2
                   +c_0^{-2}\|{w}\|_2^2\bigr)-C_I|\pair h{P_a}|^2.
\end{aligned}
\end{equation}
For the same $(h,f,w)\in Y$, the dissipation form satisfies
\begin{equation}\label{eq:centered-quadratic-derivative}
 \mathcal D_a^c(h,f,w)\ge d_I\|(h,f,w)\|_Y^2
                         -C_I|\pair h{P_a}|^2.
\end{equation}
The constants are uniform for $a\in I\Subset(0,a_*)$.
\end{proposition}

\begin{proof}
\paragraph{Step 1. Extending the operator bounds to general angles.}
The Schr\"odinger and density operators have radial profile coefficients.
The translation terms use pairings with derivatives of the reference
profiles. Lemmas~\ref{lem:anchoring}
and~\ref{lem:RX} give the full-space localized pairings and reduced
inverse. The commutators
\eqref{eq:R-commutator}--\eqref{eq:K-commutator} are full-space
operator identities. In particular, multiplication by a profile
converts a differentiated wave error into a weighted
negative-Sobolev source; for example,
\[
 P_a\diver {w}=\diver(P_a{w})-\nabla P_a\cdot {w}.
\]
The reduced inverse restores two derivatives. The compositions
and adjoint estimates of Appendix~\ref{app:tests} therefore apply
to general angles. They make the displayed energy and dissipation
forms continuous on $Y$, also when the Schr\"odinger component has
infinite global mass. The same estimates extend
Lemma~\ref{lem:form-gradients} to the full space. The profile tests
in its wave components and their spatial gradients belong to
$L^2$, uniformly on $I$.

\paragraph{Step 2. The translation terms in the dissipation form.}
In \eqref{eq:centered-leading-form}, translation modulation contributes
the quadratic expression
\[
 -\sum_{j=1}^2\frac{\nu_j^{\mathrm{lin}}(f,w)}\beta
   \bigl[({f},\partial_j(N_a+P_a^2))
                    +c_0^{-2}({w},\partial_jV_a)\bigr].
\]
We estimate it together with the scale contribution. The
constants $C$ in the following small-$a$ estimates are uniform
as $a\downarrow0$; factors of $a$ are kept before choosing
$a_*$. The stationary identity and uniform profile bounds give
\[
 \|\nabla(N_a+P_a^2)\|_2\le Ca^2,\qquad
 \|\partial_jV_a\|_2\le Ca.
\]
The density and velocity parts of
$\nu_j^{\mathrm{lin}}(f,w)/\beta$ are bounded by
$C\|{f}\|_2$ and $Ca^{-1}\|{w}\|_2$, respectively.
All terms except the velocity square are therefore bounded by
\begin{equation}\label{eq:center-mixed-errors}
 C\{a^2\|{f}\|_2^2+a\|{f}\|_2\|{w}\|_2\}.
\end{equation}
For the velocity square, use the two profile tests in
\eqref{eq:center-velocity-bracket}. Integration by parts in
\eqref{eq:linear-center-output} and the profile limits give
\[
 \begin{aligned}
 \frac{\nu_j^{\mathrm{lin}}(0,w)}\beta
   &=\frac{({w},\nabla\partial_j(Q^2))}{c_0a\|\partial_j(Q^2)\|_2^2}+o(a^{-1})\|{w}\|_2,\\
 c_0^{-2}({w},\partial_jV_a)
   &=\frac a{c_0}({w},\partial_j(yQ^2))+o(a)\|{w}\|_2.
 \end{aligned}
\]
Their product tends to
$c_0^{-2}({w},\nabla\partial_j(Q^2))({w},\partial_j(yQ^2))/\|\partial_j(Q^2)\|_2^2$.
Indeed $\partial_jN_a\to-\partial_jQ^2$ and
$a^{-1}\partial_jV_a\to c_0\partial_j(yQ^2)$ in $L^2$,
and the denominator tends to $\|\partial_j(Q^2)\|_2^2>0$; the differentiated
test in the velocity output has the same convergence.
Combining this negative product with the scale contribution
and the bare velocity square gives exactly
$c_0^{-2}$ times \eqref{eq:center-velocity-bracket} in the
limit. Proposition~\ref{prop:center-velocity-gap} gives its
strict gap. The density terms in \eqref{eq:center-mixed-errors}
will be absorbed together with the other small-$a$ terms below.

\paragraph{Step 3. Transferring the gaps to positive $a$.}
The scalar forms converge in local form norm by splitting at a
large fixed radius. For $2<p<3$, the polynomial tail is controlled by
\[
 \int_{|y|>R}\langle y\rangle^{-3}|h|^2
 \le CR^{p-3}\int\langle y\rangle^{-p}|h|^2
 \le C_pR^{p-3}\norm h_{\mathrm{loc}}^2.
\]
On the fixed ball the form difference is
$o_a(1)\norm h_{\mathrm{loc}}^2$. Choose $R$ first and then $a$.
For the inverse and density blocks, the augmented-resolvent argument in
Appendix~\ref{app:energy-coercivity} uses the full-space kernel
$\C P_a$, the potential supremum norm and elliptic regularity;
none of these steps requires radial sources. It therefore gives
\[
 \|R_a-R_Q\|_{L^2(\R^2)\to H^2(\R^2)}\longrightarrow0.
\]
Thus $R_aP_a\to R_QQ$ in $L^2\to H^2$.
Applying $A_a=-2\Delta-y\cdot\nabla N_a$ to $R_aP_a$
and using $(2-\Lambda)P_a\to Q-rQ'$ in $L^\infty$ gives
convergence of $D_a$ to $\mathcal D_Q$ in $L^2$ operator norm.
The density energy operator $1-2P_aR_aP_a$ converges to
$\mathcal B_Q$ in the same norm. Since these are full-space
bounds, they are uniform in the angular mode.

The constraints converge as well. If $\pair h{P_a}=0$, then
\[
 |\pair hQ|=|\pair h{Q-P_a}|
                    \le o_a(1)\|h\|_{\mathrm{loc}}.
\]
Removing this small $Q$ component transfers both scalar gaps
to the condition $\pair h{P_a}=0$.

For the density constraints one has explicitly
\[
 \begin{aligned}
 |({f},Q^2)|&=|({f},Q^2-P_a^2)|=o_a(1)\norm {f}_2,\\
 |({f},Q\partial_jQ)|
 &=|({f},Q\partial_jQ+\tfrac12\partial_jN_a)|
                                      =o_a(1)\norm {f}_2.
 \end{aligned}
\]
The three limiting tests $Q^2,Q\partial_1Q,Q\partial_2Q$ are
nonzero and mutually orthogonal by parity. Subtracting their three
orthogonal projections gives ${f}^{\,0}$ satisfying the limiting
constraints and
$\norm{{f}-{f}^{\,0}}_2=o_a(1)\norm {f}_2$.
Since $\mathcal B_Q$ is bounded on $L^2$,
\[
 |(\mathcal B_Q{f},{f})-(\mathcal B_Q{f}^{\,0},{f}^{\,0})|
 \le\norm{\mathcal B_Q}_{\mathrm{op}}
       \norm{{f}-{f}^{\,0}}_2(\norm {f}_2+\norm{{f}^{\,0}}_2)
 =o_a(1)\norm {f}_2^2.
\]
Proposition~\ref{prop:centered-density-energy} therefore gives a
strict density energy gap for sufficiently small positive $a$. For dissipation use
Corollary~\ref{cor:full-centered-density}, including its version with moving constraints. The Schr\"odinger energy and dissipation gaps are given by
Proposition~\ref{prop:angular-scalar-velocity}.

For the dissipation estimate, first suppose
$\pair h{P_a}=0$. Let $c>0$ be the minimum of the scalar,
density and full-space velocity gaps. The leading form is bounded below by
\[
 \begin{aligned}
 &c\bigl(\|h\|_{\mathrm{loc}}^2+\|{f}\|_2^2
                                      +c_0^{-2}\|{w}\|_2^2\bigr)\\
 &\quad-o_a(1)\|(h,f,w)\|_Y^2
       -C\bigl(a^2\|{f}\|_2^2+a\|{f}\|_2\|{w}\|_2\bigr).
 \end{aligned}
\]
For any fixed $\eta>0$, Young's inequality gives
$Ca\|{f}\|_2\|{w}\|_2\le\eta c_0^{-2}\|{w}\|_2^2
+C_\eta a^2\|{f}\|_2^2$. Choose $\eta<c/4$, then $a_*$ so that
the $a^2$ and form-convergence errors consume less than another
quarter of the gaps. The leading form is therefore bounded below
by a fixed positive constant times $\|(h,f,w)\|_Y^2$.
The same positive branch bound $a_*$ works for the energy and
dissipation forms. Fix $I\Subset(0,a_*)$ after choosing this bound.

The preceding expressions are bounded quadratic forms on $Y$.
Smooth approximation with correction of the displayed pairings
therefore extends these algebraic estimates to $Y$.
For arbitrary $h$, subtract its $P_a$ component. The transferred
Schr\"odinger energy gap and continuity of the local form give
\[
 ((L_a-1)h,h)_\R\ge c_I\|h\|_{\mathrm{loc}}^2
                                      -C_I|\pair h{P_a}|^2.
\]
Together with the transferred density energy gap, this proves
\eqref{eq:centered-energy-form}. Applying the corresponding scalar
dissipation bound to $h$ leaves the density, scale and translation
terms unchanged. The lower bound just obtained for their sum
therefore proves \eqref{eq:centered-quadratic-derivative}.

\end{proof}

\subsection{Nonlinear energy estimate and bootstrap closure}
\label{sec:centered-nonlinear}

The radial construction provides the energy $\E=E+B_\beta$.
Center motion adds $\nu\cdot\nabla\U$ to the rescaled equations.
We compute its energy contribution in the next lemma: the leading
wave pairing is precisely the term retained in
\eqref{eq:centered-leading-form}. The full-space coercivity estimates give
the nonlinear dissipation estimate and bootstrap closure.

\begin{lemma}[Energy contribution of the translation term]
\label{lem:center-transport-work}
For a smooth solution in the nonradial modulation coordinates on the bootstrap interval, let
$(\E_s)_{\mathrm{center}}$ denote the contribution of
$\nu\cdot\nabla\U$ to the energy derivative, with the parameters
held fixed in the variation. Then
\[
\begin{aligned}
 \bigg|(\E_s)_{\mathrm{center}}
 &-\sum_j\nu_j^{\mathrm{lin}}(\epsn,\epsv)
       \bigl[(\epsn,\partial_j(N_a+P_a^2))
                         +c_0^{-2}(\epsv,\partial_jV_a)\bigr]\bigg|\\
 &\hspace{12mm}\le C_I\beta(\beta+\|\veps\|_Y)\|\veps\|_Y^2.
\end{aligned}
\]
The outputs \eqref{eq:m-output} and \eqref{eq:eta-linear}
remain valid with the same orders of remainder.
\end{lemma}

\begin{proof}
\label{app:center-transport-details}
Use the energy $\E=E+B_\beta$ on the
translation orthogonality conditions. The linear translation rates are those of
\eqref{eq:linear-center-output}.
\paragraph{Energy variation under translation.}
To compute the contribution of $\nu_j\partial_j\U$, hold
$a$ and $\lambda$ fixed and translate only the full state:
replace $\U$ by
$\U(y+h e_j)$, where $e_j$ is the $j$th coordinate vector.
Write $\mathcal T_j$ for differentiation at $h=0$.
The profiles in the error decomposition stay fixed in this
calculation; hence
\[
 \mathcal T_j\veps=\partial_j\U_*+\partial_j\veps.
\]
Translation covariance gives
$R'_N[\partial_jN]=[\partial_j,R_N]$. Since
$z=\eps+R_N(Q_{a,\beta}\epsn)$, the product rule yields
\[
\begin{aligned}
 \mathcal T_jz={}&\partial_j\eps+\partial_jQ_{a,\beta}
       +[\partial_j,R_N](Q_{a,\beta}\epsn)
       +R_N\bigl(Q_{a,\beta}(\partial_j\epsn+\partial_jN_a)\bigr)\\
 ={}&\partial_jz+b_j,\\
 b_j={}&\partial_jQ_{a,\beta}-R_N((\partial_jQ_{a,\beta})\epsn)
                                      +R_N(Q_{a,\beta}\partial_jN_a).
\end{aligned}
\]
All terms in $b_j$ are localized $H^2$ functions. For the
interaction energy, put $L_N=-\Delta+1+N$ and use
\[
\begin{aligned}
 L_N\eps+Q_{a,\beta}\epsn&=L_N\Psi-L_aQ_{a,\beta},\\
 \epsn+2\Rea(\overline{Q_{a,\beta}}\eps)+|\eps|^2
                       &=N+|\Psi|^2-(N_a+P_a^2).
\end{aligned}
\]
Differentiation at fixed $a,\lambda$ gives
\[
\begin{aligned}
 \mathcal T_j(E_\beta^{(2)}+H_3)
 ={}&2\Rea\pair{\partial_j\Psi}{L_N\eps+Q_{a,\beta}\epsn}
       +(\partial_jN,\epsn+2\Rea(\overline{Q_{a,\beta}}\eps)+|\eps|^2)\\
 &+c_0^{-2}(\partial_jV,\epsv).
\end{aligned}
\]
The full spatial derivatives cancel because
\[
\begin{aligned}
 2\Rea\pair{\partial_j\Psi}{L_N\Psi}
                         &=-(\partial_jN,|\Psi|^2),\\
 (\partial_jN,N)&=(\partial_jV,V)=0.
\end{aligned}
\]
Consequently
\[
\begin{aligned}
 \mathcal T_j(E_\beta^{(2)}+H_3)
 ={}&-2\Rea\pair{L_aQ_{a,\beta}}{\partial_j\Psi}
       -(N_a+P_a^2,\partial_jN)
       -c_0^{-2}(V_a,\partial_jV)\\
 ={}&2\Rea\pair\eps{\partial_j(L_aQ_{a,\beta})}
       +(\epsn,\partial_j(N_a+P_a^2))
       +c_0^{-2}(\epsv,\partial_jV_a).
\end{aligned}
\]
The reference constants cancel since
\[
 2\Rea\pair{L_aQ_{a,\beta}}{\partial_jQ_{a,\beta}}
       =-(\partial_jN_a,P_a^2),\qquad
 (N_a,\partial_jN_a)=(V_a,\partial_jV_a)=0.
\]
Since $\Rea\pair{\partial_jz}z=0$, the additional energy work is
\begin{equation}\label{eq:exact-center-transport-work}
\begin{aligned}
 (\E_s)_{\mathrm{center}}&=\sum_j\nu_jW_j(\veps),\\
 W_j(\veps)={}&2\Rea\pair\eps{\partial_j(L_aQ_{a,\beta})}
       +(\epsn,\partial_j(N_a+P_a^2))
       +c_0^{-2}(\epsv,\partial_jV_a)\\
 &-2\Rea\pair{b_j}z
       +2B_\beta(\veps,\partial_j\U_*+\partial_j\veps).
\end{aligned}
\end{equation}
The derivatives of the wave errors in the last term are
transferred to the gradients of the quadratic correction.
With the gradient representatives of
Lemma~\ref{lem:form-gradients} in Appendix~\ref{app:tests},
\[
\begin{aligned}
 2B_\beta(\veps,\partial_j\veps)
 ={}&\Rea\pair{\partial_j\eps}{d_\eps B_\beta(\veps)}
       -(\epsn,\partial_jd_{\epsn}B_\beta(\veps))\\
 &-(\epsv,\partial_jd_{\epsv}B_\beta(\veps)).
\end{aligned}
\]
Both the exponentially localized tests and the finite polynomial
profile tests have spatial gradients in $L^2$, with norms at most
$C_I\beta\|\veps\|_Y$. The Schr\"odinger term uses
$\|\partial_j\eps\|_2\le\|\eps\|_{\mathrm{loc}}$.
Since $\|\partial_j\U_*\|_Y\le C_I$, these estimates give
\begin{equation}\label{eq:B-transport-tests}
\begin{aligned}
 |2B_\beta(\veps,\partial_j\U_*)|
                              &\le C_I\beta\|\veps\|_Y,\\
 |2B_\beta(\veps,\partial_j\veps)|
                              &\le C_I\beta\|\veps\|_Y^2.
\end{aligned}
\end{equation}
The displayed integration by parts defines these terms continuously
on energy data, with $\epsn,\epsv$ only in $L^2$.

\paragraph{The leading translation term.}
At $\epsn=0$, $b_j=\partial_jQ_{a,\beta}+R_a(Q_{a,\beta}\partial_jN_a)$.
Differentiating $L_aP_a=0$ in $y_j$ gives
\[
 L_a\partial_jP_a=-P_a\partial_jN_a,\qquad
 (\partial_jP_a,P_a)=0,\qquad
 -R_a(P_a\partial_jN_a)=\partial_jP_a.
\]
Thus $\partial_jP_a+R_a(P_a\partial_jN_a)=0$.
Subtracting this zero from $b_j$ gives the exact expansion
\[
\begin{aligned}
 b_j={}&\partial_j(Q_{a,\beta}-P_a)
       +(R_N-R_a)(Q_{a,\beta}\partial_jN_a)
       +R_a((Q_{a,\beta}-P_a)\partial_jN_a)\\
 &-R_N((\partial_jQ_{a,\beta})\epsn).
\end{aligned}
\]
The phase differences are $O_{H^2_d}(\beta)$, and the inverse
difference contributes $C_I\|\epsn\|_2$ by
Lemma~\ref{lem:current-inverse}. Hence
\[
 \|b_j\|_{H^1_d}\le C_I(\beta+\|\epsn\|_2),\qquad
 \|\partial_j(L_aQ_{a,\beta})\|_{H^1_d}\le C_I\beta.
\]
Thus the exact work satisfies
\begin{equation}\label{eq:translation-relative-identities}
 \left|W_j(\veps)
       -(\epsn,\partial_j(N_a+P_a^2))
       -c_0^{-2}(\epsv,\partial_jV_a)\right|
       \le C_I(\beta+\|\veps\|_Y)\|\veps\|_Y.
\end{equation}
Differentiating the translation constraints and using
\eqref{eq:center-rate} gives
\begin{equation}\label{eq:center-output-remainder}
 |\nu-\nu^{\mathrm{lin}}(\epsn,\epsv)|
       \le C_I\beta\|\veps\|_Y^2.
\end{equation}
Indeed, subtracting the definition of $\nu_j^{\mathrm{lin}}$
in \eqref{eq:linear-center-output} from \eqref{eq:center-system}
gives the exact identity
\[
\begin{aligned}
 \|\partial_jN_a\|_2^2(\nu_j-\nu_j^{\mathrm{lin}}(\epsn,\epsv))
 ={}&-a_s(\epsn,\partial_a\partial_jN_a)
       -\left(\frac{\lambda_s}{\lambda}+\beta\right)
                                (\Dtwo\epsn,\partial_jN_a)\\
 &+\sum_{\ell=1}^2\nu_\ell
                         (\epsn,\partial_\ell\partial_jN_a).
\end{aligned}
\]
The denominator has a positive lower bound on $I$.
Transpose $\Dtwo$ to its smooth profile test and apply
\eqref{eq:center-rate}; the right-hand side is bounded by
$C_I\beta\|\veps\|_Y^2$. This proves
\eqref{eq:center-output-remainder}. Multiplying
\eqref{eq:translation-relative-identities} by $\nu_j$ and using
this remainder gives
\[
\begin{aligned}
 \bigg|\sum_j\nu_jW_j(\veps)
 &-\sum_j\nu_j^{\mathrm{lin}}(\epsn,\epsv)
       \bigl[(\epsn,\partial_j(N_a+P_a^2))
                         +c_0^{-2}(\epsv,\partial_jV_a)\bigr]\bigg|\\
 &\hspace{12mm}\le C_I\beta(\beta+\|\veps\|_Y)\|\veps\|_Y^2.
\end{aligned}
\]

The scale output follows directly from
\eqref{eq:centered-density-row}. The additional products beyond
the leading quotient in \eqref{eq:m-output} satisfy
\[
\begin{aligned}
 &\left|\left(\frac{\lambda_s}{\lambda}+\beta\right)
                                      (\Dtwo\epsn,P_a^2)\right|
       +2|a_s(\epsn,P_a\partial_aP_a)|\\
 &\qquad+|\nu\cdot(\nabla\epsn,P_a^2)|
                              \le C_I\beta\|\veps\|_Y^2.
\end{aligned}
\]
Here the spatial derivatives act on $P_a^2$, so each factor
paired with $\epsn$ is a bounded density test. This proves
\eqref{eq:m-output} with the same order of remainder.
The amplitude equation gains the real part of
$\sum_j\nu_j\pair{\partial_jQ_{a,\beta}}\phi/\pair{Q_{a,\beta}}\phi$ from
\eqref{eq:centered-error-projection}. Its size is at most
$C_I\beta\|\veps\|_Y^2$. Thus the expansion proving
\eqref{eq:eta-linear} is unchanged to its stated precision.
\end{proof}

\begin{proposition}[Nonlinear energy estimate for general perturbations]
\label{prop:centered-nonlinear}
On the modulation conditions \eqref{eq:chart} and
\eqref{eq:center-conditions}, for $a\in I$ and sufficiently small
$\beta$, $\norm{\veps}_Y$ and $\norm\eps_2$, one has
\begin{align}
 c_I\norm{\veps}_Y^2&\le\E\le C_I\norm{\veps}_Y^2,
 \label{eq:centered-nonlinear-coercivity}\\
 \E_s&\le-d_I \beta\norm{\veps}_Y^2
       +C_I \beta(\beta+\norm{\veps}_Y+\norm\eps_2)
                                            \norm{\veps}_Y^2.
 \label{eq:centered-nonlinear-energy}
\end{align}
In particular the smallness thresholds can be chosen independently
of the current $\lambda$ so that
\begin{equation}\label{eq:centered-closed-energy}
 \E_s\le-\kappa_I\lambda\E,\qquad
 |a_s|+|\frac{\lambda_s}{\lambda}+\beta|+|\nu|\le K_I\lambda\sqrt\E.
\end{equation}
\end{proposition}

\begin{proof}
The algebraic expansion \eqref{eq:energy-comparison} holds
without symmetry. Apply \eqref{eq:centered-energy-form} with
$h=z=\eps-K_N\epsn$ and use the inverse bound to compare
$\|z\|_{\mathrm{loc}}^2+\|\epsn\|_2^2+c_0^{-2}\|\epsv\|_2^2$
with $\|\veps\|_Y^2$. By \eqref{eq:nonlinear-moment-bound},
\[
 |\pair z{P_a}|^2
       \le C_I\|\veps\|_Y^4.
\]
The cubic interaction obeys \eqref{eq:H3} and
$|B_\beta|\le C_I\beta\|\veps\|_Y^2$. Hence
\[
 \E\ge\bigl[c_I-C_I(\beta+\|\veps\|_Y+\|\eps\|_2)\bigr]
                                                  \|\veps\|_Y^2.
\]
Choosing the thresholds small proves the lower bound in
\eqref{eq:centered-nonlinear-coercivity}; boundedness of the same
forms gives its upper bound.

\paragraph{The contribution of the translation term.}
Lemma~\ref{lem:center-transport-work} separates its leading wave
pairing from the higher-order transport terms. It also preserves
the amplitude and scale expansions used in the radial construction.
Thus the only change in the principal derivative is the translation
pairing already present in \eqref{eq:centered-leading-form}.

\paragraph{Dissipation.}
All terms computed in \eqref{eq:master} retain the same bounds
without symmetry. Adding the contribution of the translation term gives the
full derivative. Substitute \eqref{eq:m-output} and the leading
translation rates. The resulting principal term is
$-\beta\mathcal D_a^c(z,\epsn,\epsv)$. The norms of
$(z,\epsn,\epsv)$ and $\veps$ in $Y$ are equivalent by the bound on
$K_N\epsn$. By \eqref{eq:nonlinear-moment-bound},
$|\pair z{P_a}|\le C_I\|\veps\|_Y^2$.
The dissipation estimate \eqref{eq:centered-quadratic-derivative}
then gives, after decreasing $d_I$ if needed,
\[
 -\beta\mathcal D_a^c(z,\epsn,\epsv)
 \le-d_I\beta\|\veps\|_Y^2+C_I\beta\|\veps\|_Y^4.
\]
The last term fits the remainder since the scale and error
thresholds are less than one. The only global mass term is
estimated by \eqref{eq:H3} and
$|\lambda_s/\lambda|\le C_I\beta$:
\[
 |2(\lambda_s/\lambda)H_3|
       \le C_I\beta\|\eps\|_2\|\veps\|_Y^2.
\]
Combining these bounds proves
\begin{equation}\label{eq:centered-remainder}
 \E_s\le-d_I\beta\|\veps\|_Y^2
       +C_I\beta(\beta+\|\veps\|_Y+\|\eps\|_2)\|\veps\|_Y^2.
\end{equation}
Fix the constants in \eqref{eq:centered-nonlinear-coercivity}
and \eqref{eq:centered-remainder}. Choose the scale and error
thresholds so that
\[
 C_I(\beta_*+2\delta)\le d_I/2.
\]
On the bootstrap interval this gives
\[
 \E_s\le-\tfrac12d_I\beta\|\veps\|_Y^2
       \le-\frac{d_Ic_0a_-}{2C_I}\lambda\E.
\]
Coercivity and \eqref{eq:center-rate} give the parameter part
of \eqref{eq:centered-closed-energy}. The inequalities pass
to energy solutions by Lemma~\ref{lem:energy-solutions},
using the transport tests just proved. All constants are fixed
before the initial neighborhood in the bootstrap argument is
chosen.
\end{proof}

\begin{proof}[Proof of Theorem~\ref{thm:main}]
\paragraph{Step 1. Initial margins and bootstrap improvement.}
With the full-space coercivity estimates established, fix $I$ about $a_0$ and
then the thresholds and constants.
Use \eqref{eq:bootstrap-assumptions} on the full spaces with the
two translation conditions. Choose $C_\star$ as in Step 3 of the radial
proof, now covering the estimates for general perturbations. Then select the late
reference point and initial neighborhood satisfying
\eqref{eq:bootstrap-margins}.
The bounds \eqref{eq:centered-closed-energy} give
\[
 \int_{s_0}^S\lambda\sqrt\E\dd s\le C_I\sqrt{\mathcal E(s_0)},\qquad
 \int_{s_0}^S|a_s|\dd s\le C_I\sqrt{\mathcal E(s_0)}.
\]
The mass identity \eqref{eq:mass-expansion} is translation
invariant. The local pairing bound
\eqref{eq:nonlinear-moment-bound} follows from the normalized
ground state and the nonlinear constraints, without radial
symmetry. Thus \eqref{eq:overlap-bound} holds on the nonradial
decomposition, and the same two-time mass subtraction gives
\eqref{eq:global-mass}. This propagates the global $L^2$
smallness required in \eqref{eq:centered-remainder}. With $C_\star$ fixed for the constants in the full-space estimates, the margins
\eqref{eq:bootstrap-margins} give explicitly
\[
 \begin{aligned}
 \|\veps(s)\|_Y&\le C_\star\sqrt{\mathcal E(s_0)}<\delta/2,\\
 \|\eps(s)\|_2^2&\le\|\eps(s_0)\|_2^2+C_\star\sqrt{\mathcal E(s_0)}
                                      <\delta^2/2,\\
 |a(s)-a_0|&\le|a(s_0)-a_0|+C_\star\sqrt{\mathcal E(s_0)}<2\delta_a.
 \end{aligned}
\]
These improve the same conditions on the same interval; neither
parameter convergence nor a power decay estimate has been used.
The scale is strictly decreasing and
$\beta\le c_0a_+\lambda(s_0)<\beta_*/2$. The extra translation block stays
invertible: its reference Gram matrix has entries
$(\partial_jN_a,\partial_kN_a)$, which equal
$\norm{\partial_1N_a}_2^2\delta_{jk}$ by radial symmetry.
This positive number has a positive minimum for $a\in I$;
small error preserves that lower bound. These thresholds and inverse bounds are uniform under
translations of the center.

\paragraph{Step 2. Energy solutions and continuation.}
To apply the Cauchy theory of \cite{BC,GTV,SanwalLWP}, decompose the physical
velocity into gradient and divergence-free parts.
In Fourier variables the two orthogonal projections, for
$\xi\ne0$, are
\[
 \widehat{v_{\mathrm{grad}}}(\xi)
       =\frac{\xi\otimes\xi}{|\xi|^2}\widehat v(\xi),\qquad
 v_{\mathrm{df}}=v-v_{\mathrm{grad}}.
\]
The equations imply
\[
 \partial_tv_{\mathrm{df}}=0,\quad\diver v_{\mathrm{df}}=0,
 \quad v_{\mathrm{grad}}=-\nabla\Delta^{-1}n_t,\quad
 \|v_{\mathrm{grad}}\|_2=\|n_t\|_{\dot H^{-1}}.
\]
The divergence-free part equals its initial value. The gradient
part is recovered from the scalar wave formulation of the cited theory,
so its continuous dependence and bounded-energy continuation
criterion apply to our full velocity data.
The compact-interval approximation argument in
Lemma~\ref{lem:energy-solutions} applies after including translations:
the translation tests and correction gradients in the preceding proof
are continuous on the full energy space. Thus one passes the
full-space integrated inequalities to energy solutions before any
endpoint is considered.

If $s_*<\infty$, Step 4 of the radial proof gives a positive lower
scale bound and a finite limiting physical time. The center rate
$|x_{c,s}|\le K_I\lambda^2\sqrt{\mathcal E(s_0)}$ gives a limit for $x_c$
on this finite rescaled interval. The physical energy norms are
bounded by the same scaling identities used in Step 5. The local
theory therefore extends the solution; the uniformly invertible
translation block and the other modulation pivots extend the chart.
Strong continuity of translations and the strict improved margins
then contradict maximality of $s_*$. This proves $s_*=\infty$
in the nonradial case as well. The quantitative orbital bound at
the end of the radial closure applies without a symmetry restriction
and gives
\eqref{eq:main-orbital-bound}.

\paragraph{Step 3. Limits, rates and blow-up.}
The scale comparison \eqref{eq:lambda-ode}--\eqref{eq:scale-comparison}
and Proposition~\ref{prop:quantitative-rates}, with the nonradial
constants, give $T$, $a_\infty$, $x_\infty$ and all stated rates. In particular,
\eqref{eq:center-convergence-quantitative} and the scale comparison
give $|x_c-x_\infty|/\lambda\to0$. The radiation mass conclusion
follows from Corollary~\ref{cor:radiation-mass} on the same chart.

Translation does not change the
gradient norm; the last lower bound in Section~\ref{sec:closure}
therefore still proves genuine blow-up at $T$. This proves the
claimed open-neighborhood orbital stability and the stated local
attraction.
\end{proof}

\subsection{An explicit exponent on the small positive branch}
\label{sec:explicit-rate}

Under the sound-speed normalization at the end of
Section~\ref{sec:closure}, the unit-speed scale is $c_0\lambda$.
The time $s$, the parameter $a$ and $\beta=c_0a\lambda$ are unchanged,
and the rescaled velocity error becomes $\epsv/c_0$. Thus it suffices
to choose the branch bound and exponent at unit sound speed.

Use the corrected field $z$ in the comparison
\[
 B=\|\nabla z\|_2^2+\tfrac12\|\epsn\|_2^2+
                  \tfrac1{2c_0^2}\|\epsv\|_2^2 .
\]
At the limiting profile the energy form in these coordinates is
bounded above by $B$. Its scalar potential is $-Q^2$, and
$\mathcal B_Q=1-2QR_QQ\le1$ by nonnegativity of the reduced inverse.
For the Schr\"odinger dissipation form, the radial component orthogonal
to $Q$ has gradient coefficient $2/5$ by
\eqref{eq:scalar-virial-gap}. On the modes of nonzero angular
order, \eqref{eq:angular-potentials} and
\eqref{eq:analytic-mass-upper} give the gradient coefficient
$2(1-\|Q\|_2^2/(4\pi))\ge2/15$. The density form has gap $1/25$
on the limiting translation orthogonality conditions by
Corollary~\ref{cor:full-centered-density}. The velocity form
including scale and translation work is bounded below by
$c_0^{-2}\|\epsv\|_2^2/50$ by
Proposition~\ref{prop:center-velocity-gap}. These three bounds
give at least $B/25$ for the limiting derivative form.

The static forms converge in the corrected coordinates at unit
speed, uniformly in angular order, by Proposition~\ref{prop:centered-linear}.
In the velocity-only terms the factor $a^{-1}$
from the scale or translation rate pairs with an $O(a)$ profile test,
so the velocity form converges in operator norm. The mixed
density products tend to zero by \eqref{eq:center-mixed-errors}.
Projecting the density onto the limiting constraints changes it
by $o_a(1)\|\epsn\|_2$ in $L^2$, using the three tests
$Q^2,Q\partial_1Q,Q\partial_2Q$ as in that proof. For the scalar
constraint, \eqref{eq:nonlinear-moment-bound} and local profile
convergence give
\[
 |\pair zQ|
 \le o_a(1)\|z\|_{\mathrm{loc}}+C_I\|\veps\|_Y^2.
\]
Remove the $Q$ component and apply \eqref{eq:constrained-anchoring}.
The uniform bound on $R_aP_a$ and \eqref{eq:response-comparison}
give $\|\eps-z\|_{\mathrm{loc}}\le C\|\epsn\|_2$ after fixing $I$
and shrinking the scale and error thresholds, with $C$ independent
of $I$. After absorbing $o_a(1)$ and $C_I\|\veps\|_Y$, this gives
$\|\veps\|_Y^2\lesssim B$ with a fixed constant. The scalar projection
also changes the limiting scalar form and the gradient square by
at most $[o_a(1)+C_I\|\veps\|_Y]\|\veps\|_Y^2$.

The errors in \eqref{eq:energy-comparison}, including $H_3$ and
$B_\beta$, are bounded by
$C_I(\beta+\|\veps\|_Y+\|\eps\|_2)\|\veps\|_Y^2$.
The energy derivative with translation modulation has principal term
$-\beta\mathcal D_a^c(z,\epsn,\epsv)$, with the same bound multiplied
by $\beta$. Thus for any $0<\theta<1/25$
we can first choose the branch bound, then fix $I$ around $a_0$
and shrink the scale and error thresholds to obtain
\[
 \E\le(1+\theta)B,\qquad -\E_s\ge(1/25-\theta)\beta B .
\]
The static comparison uses a single branch bound. The nonlinear
thresholds may depend on $I$. In particular, the amplitude
pivot involves $M'(a)\asymp a$ by \eqref{eq:mass-slope-expansion}.

Strengthen the initial scale-rate margin to
$K_I\sqrt{\mathcal E(s_0)}\le\theta c_0a_-$ in \eqref{eq:bootstrap-margins}.
Monotonicity of $\E$ then gives
$-\lambda_s/\lambda\le(1+\theta)\beta$ on the bootstrap interval.
Consequently $\E\lambda^{-2\sigma}$ is nonincreasing if
\[
 2\sigma(1+\theta)^2\le 1/25-\theta .
\]
With $\theta=1/50$, the coarser bounds $\E\le2B$,
$-\E_s\ge\beta B/50$ and $-\lambda_s/\lambda\le3\beta/2$
already allow $\sigma=1/300$. For any fixed $0<\sigma<1/50$,
take $\theta$ sufficiently small in the displayed condition.
The integrated inequality and the subsequent parameter rates
follow as in Proposition~\ref{prop:quantitative-rates}.
The full-space comparison also applies to radial solutions, for which
the translation work vanishes.

\appendix
\section{Profile regularity and operator estimates}
\label{app:branch}

We establish the profile-mass monotonicity and the operator
estimates used in the energy argument: parameter regularity,
spatial decay, localized tests, and weighted bounds for the
reduced inverse and the correction tests. Parameter
differentiation is always at positive $a$.

\begin{proposition}[Properties of the small positive branch]\label{prop:profile}
There is a small interval $(0,a_*)$ of positive profiles satisfying
\eqref{eq:profiles}--\eqref{eq:wave-profiles}. On every compact interval
$I\Subset(0,a_*)$, the branch is $C^5$ in $a$ in every finite spatial
Sobolev norm used below. For $0\le j\le5$, every finite $\ell$, and
$0<d<1$,
\begin{align}
 \abs{\partial_a^j\partial_r^\ell P_a(r)}
 &\le C_{I,j,\ell,d}\e^{-dr},\label{eq:P-tail}\\
 \abs{\partial_a^j\partial_r^\ell N_a(r)}
 &\le C_{I,j,\ell}\langle r\rangle^{-3-\ell}.
 \label{eq:N-tail}
\end{align}
The fields are smooth radial functions at the origin. Moreover
\begin{equation}\label{eq:mass-slope}
 M(a)=\norm{P_a}_2^2,\qquad \inf_{a\in I}M'(a)>0.
\end{equation}
\end{proposition}

Proposition~\ref{prop:slope} proves the mass slope. We prove
regularity and decay by first estimating the stationary wave
operator, then differentiating the stationary profile equation.

\subsection{Monotonicity of the profile mass}
\label{app:mass-slope}

The condition selecting $a$ removes the error component along
the ground state of $-\Delta+N$. At zero density error and
zero quadratic phase, its derivative in the branch direction is
$-M'(a)/(2\sqrt{M(a)})$. The required nondegeneracy is therefore
$M'(a)>0$.
The following argument also identifies its leading coefficient
as $a\downarrow0$.

The limiting NLS operator is denoted here by
$L_+=-\Delta+1-3Q^2$ on real radial functions. Radial
nondegeneracy is \cite[Proposition~4.1]{GM1}. Since its potential
decays and its essential spectrum is $[1,\infty)$, this gives a
bounded inverse $L_+^{-1}:L^2_\rad\to H^2_\rad$ by elliptic
regularity. For the classical NLS quadratic forms see
\cite[Section~2]{WeinsteinModulation}. We recall the radial
nondegeneracy argument in Appendix~\ref{app:energy-coercivity}.
For a smooth exponentially decreasing radial source $g$, define
$\mathcal R_ag$ as the solution, regular at the sonic radius $r=1/a$,
of
\begin{equation}\label{eq:sonic-ode}
 (a^2r^2-1)\partial_r(\mathcal R_ag)
                   +3a^2r\mathcal R_ag=\partial_r g.
\end{equation}
This stationary wave solution operator gives
$N_a=\mathcal R_a(P_a^2)$. Its kernel and estimates are proved in
Appendix~\ref{app:branch}. The coefficient of $a^2$ in
$\mathcal R_ag+g$ is the linear potential
\begin{equation}\label{eq:green-coefficient}
 \Xi[g](r)=\int_r^\infty t g(t)\dd t-r^2g(r),\qquad
 \Xi[g]'=-r(3g+rg').
\end{equation}
\begin{proposition}[Derivative of the profile mass]
\label{prop:slope}
On the small positive stationary branch,
\begin{equation}\label{eq:mass-slope-expansion}
 \frac{M'(a)}a\longrightarrow\int_{\R^2}|y|^2Q(y)^4\dd y>0
                     \qquad(a\downarrow0).
\end{equation}
In particular $M'(a)>0$ for every sufficiently small $a>0$.
\end{proposition}

\begin{proof}
\paragraph{The stationary linearized inverse.}
Write $h_a=P_a-Q$ in this proof. The fixed-point map used in
\cite[Section~IV.1, Theorem~4.2]{GM1}, with our sign convention, is
\[
 \mathcal T_a(h)=h-L_+^{-1}\bigl\{
   -\Delta(Q+h)+(Q+h)+(Q+h)\mathcal R_a((Q+h)^2)\bigr\}.
\]
It acts on a small ball in the real radial space $H^2$, and
$\mathcal T_a(h_a)=h_a$. The estimates on its linear, quadratic
and cubic terms in \cite[Propositions~4.7 and~4.9]{GM1} give
\[
 \begin{aligned}
 \|\mathcal T_a(h)-\mathcal T_a(k)\|_{H^2}
 &\le C(a^2+\delta+\delta^2)\|h-k\|_{H^2},\qquad
                         \|h\|_{H^2},\|k\|_{H^2}\le\delta,\\
 \|h_a\|_{H^2}&\le Ca^2.
 \end{aligned}
\]
Choose $\delta$ and then $a_*$ so that the first coefficient is at
most $1/2$. Dependence on $h$ is polynomial in the continuous
multilinear maps in those estimates. Write $D\mathcal T_a(h_a)$
for its Fr\'echet derivative with respect to $h$ at the fixed point;
its $H^2\to H^2$ norm is at most $1/2$.

Differentiating \eqref{eq:profiles} with respect to $P$, with
$N=\mathcal R_a(P^2)$, gives the real radial operator
\[
 \mathcal L_a\xi=(-\Delta+1+N_a)\xi+2P_a\mathcal R_a(P_a\xi).
\]
The definition of $\mathcal T_a$ now gives the exact factorization
\[
 \begin{aligned}
 1-D\mathcal T_a(h_a)&=L_+^{-1}\mathcal L_a,\\
 \mathcal L_a^{-1}
  &=(1-D\mathcal T_a(h_a))^{-1}L_+^{-1}
    =\sum_{j=0}^\infty(D\mathcal T_a(h_a))^jL_+^{-1},\\
 \|\mathcal L_a^{-1}\|_{L^2_\rad\to H^2_\rad}
  &\le\sum_{j=0}^\infty2^{-j}
           \|L_+^{-1}\|_{L^2_\rad\to H^2_\rad}
   =2\|L_+^{-1}\|_{L^2_\rad\to H^2_\rad}.
 \end{aligned}
\]
The resulting inverse bound is uniform as $a\downarrow0$.

\paragraph{The differentiated stationary equation.}
At positive $a$, differentiating $N_a=\mathcal R_a(P_a^2)$ means
\[
 \partial_aN_a=(\partial_a\mathcal R_a)(P_a^2)
                         +2\mathcal R_a(P_a\partial_aP_a).
\]
Here the operator derivative in the first term holds its argument fixed.
Differentiate $(-\Delta+1+N_a)P_a=0$ and insert this formula:
\[
 \mathcal L_a\partial_aP_a=-P_a(\partial_a\mathcal R_a)(P_a^2).
\]
The positive-parameter differentiability and source spaces are
proved in Appendix~\ref{app:branch}. Division by $2a$ gives
\begin{equation}\label{eq:profile-derivative-direct}
 \mathcal L_a\frac{\partial_aP_a}{2a}
       =-P_a\frac{(\partial_a\mathcal R_a)(P_a^2)}{2a}.
\end{equation}
The spatial bounds and convergence, uniform for small positive $a$,
are those of \cite[Appendix~B, Propositions~B.1, B.2 and~B.7]{GM2}.
They meet the hypotheses of Lemma~\ref{lem:green-expansion}, which gives
\[
 \left\|P_a\left[
 \frac{(\partial_a\mathcal R_a)(P_a^2)}{2a}-\Xi[P_a^2]
                         \right]\right\|_2\le Ca^2.
\]
The formula for $\Xi$ bounds $\Xi[P_a^2](r)$ by
$C(1+r)^2\e^{-dr}$ for some fixed $d>0$. Local convergence and
this common integrable bound therefore give
\[
 \begin{aligned}
 &\left\|P_a\frac{(\partial_a\mathcal R_a)(P_a^2)}{2a}
                                      -Q\Xi[Q^2]\right\|_2\\
 &\qquad\le Ca^2+\|P_a\Xi[P_a^2]-Q\Xi[Q^2]\|_2
                                      \longrightarrow0.
 \end{aligned}
\]

\paragraph{The limiting equation.}
Let $v$ be the unique real radial solution
\[
 L_+v=-Q\Xi[Q^2],\qquad v=-L_+^{-1}(Q\Xi[Q^2]).
\]
It is smooth by elliptic regularity. It also has exponential decay:
for $0<d<1$, choose $R$ so that $3Q^2\le(1-d^2)/2$ outside $R$.
Testing the equation with a truncated exponential times $v$ absorbs
the exterior potential; weighted elliptic regularity then gives
\[
 \|\e^{dr}v\|_{H^2}
 \le C_{d,R}\bigl(\|v\|_{H^2(|y|<2R)}
                         +\|\e^{dr}Q\Xi[Q^2]\|_2\bigr)<\infty.
\]
Higher spatial derivatives follow from the same equation.
To compare the operators, use the exact expansion
\[
 (\mathcal L_a-L_+)v
  =(N_a+P_a^2)v+3(Q^2-P_a^2)v
                  +2P_a\{\mathcal R_a(P_av)+P_av\}.
\]
Lemma~\ref{lem:green-expansion}, first with $g_a=P_a^2$ and then
with $g_a=P_av$, gives
\[
 \|N_a+P_a^2\|_2\le Ca^2,\qquad
 \|\mathcal R_a(P_av)+P_av\|_2\le Ca^2.
\]
Since $P_a\to Q$ in $H^2$ and $v,P_a$ are bounded,
\[
 \|(\mathcal L_a-L_+)v\|_2
 \le C_v\bigl(a^2+\|P_a^2-Q^2\|_2\bigr)\longrightarrow0.
\]
Subtract the displayed equation for $v$ from
\eqref{eq:profile-derivative-direct} and use the uniform inverse:
\[
 \begin{aligned}
 \left\|\frac{\partial_aP_a}{2a}-v\right\|_{H^2}
 \le C\biggl(&\left\|P_a\frac{(\partial_a\mathcal R_a)(P_a^2)}{2a}
                                      -Q\Xi[Q^2]\right\|_2\\
             &+\|(\mathcal L_a-L_+)v\|_2\biggr)
                                      \longrightarrow0.
 \end{aligned}
\]
Thus
\begin{equation}\label{eq:profile-slope-limit}
 \frac{\partial_aP_a}{2a}\longrightarrow v\quad\text{in }H^2.
\end{equation}

\paragraph{The sign of the limit.}
Differentiating the NLS scaling gives $L_+\Lambda Q=-2Q$.
Self-adjointness and the equation for $v$ imply
\[
 \begin{aligned}
 \frac{M'(a)}a
 &=4\left(P_a,\frac{\partial_aP_a}{2a}\right)
   \longrightarrow4(Q,v)\\
 &=-2(L_+\Lambda Q,v)
   =2(Q\Lambda Q,\Xi[Q^2])=(\Dtwo Q^2,\Xi[Q^2]).
 \end{aligned}
\]
Put $g=Q^2$ in the following radial integration, including its factor
$2\pi$. Since $\Xi[g]'=-r(3g+rg')$,
\begin{align*}
 (\Dtwo g,\Xi[g])
 &=-2\pi\int_0^\infty r^2g\,\Xi[g]'\dd r\\
 &=2\pi\int_0^\infty(3r^3g^2+r^4gg')\dd r\\
 &=2\pi\int_0^\infty r^3g^2\dd r
    +\pi[r^4g^2]_{0}^{\infty}
  =2\pi\int_0^\infty r^3g^2\dd r>0.
\end{align*}
The exponential decay and smoothness of $Q$ make all boundary
terms zero. The last integral is the constant in
\eqref{eq:mass-slope-expansion}.
\end{proof}

\subsection{Green operator for the stationary wave equation}

Recall that $\mathcal R_ag$ is the solution of
\eqref{eq:sonic-ode} regular at the sonic radius $r_a=1/a$.
Regularity excludes the singular homogeneous solution. With
$s_t=r_a+t(r-r_a)$, the regular solution has the average formula
\begin{equation}\label{eq:sonic-average}
 \mathcal R_ag(r)=a^{-2}(r+r_a)^{-3/2}
       \int_0^1 t^{1/2}(s_t+r_a)^{1/2}g'(s_t)\dd t.
\end{equation}
For $r<r_a$ and $r>r_a$, integration by parts gives respectively
\begin{align}
 \mathcal R_ag(r)&=-\frac{g(r)}{1-a^2r^2}
 +(1-a^2r^2)^{-3/2}
       \int_{ar}^1\frac{xg(x/a)}{\sqrt{1-x^2}}\dd x,
       \label{eq:sonic-interior}\\
 \mathcal R_ag(r)&=\frac{g(r)}{a^2r^2-1}
 -\frac{a^2}{(a^2r^2-1)^{3/2}}
       \int_{r_a}^r\frac{tg(t)}{\sqrt{a^2t^2-1}}\dd t.
       \label{eq:sonic-exterior}
\end{align}
The average formula controls the sonic radius. The integrated
formulas control the origin and the exterior, respectively.

\begin{lemma}[Finite-loss parameter estimates]\label{lem:green}
For $I\Subset(0,a_*)$, finite integers $m,k\ge0$, and smooth radial
$g$,
\begin{equation}\label{eq:green-loss}
 \norm{\partial_a^k\mathcal R_ag}_{H^m}
 \le C_{I,m,k}\norm g_{H^{m+k+1}}.
\end{equation}
The maps are strongly continuous in $a$. Operator-norm continuity
holds after reserving one further derivative of the source.
\end{lemma}

\begin{proof}
Cover the radial axis by a ball separated from $r_a$, a sonic annulus,
and an exterior region, with cutoffs uniform on $I$.
\paragraph{The sonic annulus.}
Differentiate \eqref{eq:sonic-average}. A spatial
derivative adds at most one derivative of $g$; a parameter derivative
does likewise since $\partial_as_t=(1-t)\partial_ar_a$ is bounded
on $I$. The radial $L^2$ change of variable on this annulus gives
\[
 \|g^{(q)}(r_a+t(r-r_a))\|_{L^2(r\dd r;\,\mathrm{annulus})}
       \le C_I t^{-1/2}\|g^{(q)}\|_{L^2(r\dd r)}.
\]
Minkowski's inequality therefore bounds every differentiated term by
\[
 C_I\sum_{q\le m+k+1}\int_0^1t^{1/2}t^{-1/2}
                      \|g^{(q)}\|_{L^2(r\dd r)}\dd t
 \le C_{I,m,k}\|g\|_{H^{m+k+1}}.
\]
All coefficient derivatives are bounded here because $a$ stays
in $I$ and the annulus is separated from zero.

\paragraph{The origin.}
Use \eqref{eq:sonic-interior} and split its integral into
the scalar annular functional and its primitive from zero to $ar$.
The latter has the form
\[
 \int_0^r\frac{tg(t)}{\sqrt{1-a^2t^2}}\dd t
 =r^2\int_0^1\frac{v g(rv)}{\sqrt{1-a^2r^2v^2}}\dd v.
\]
Its Cartesian derivatives are controlled by differentiating $g(vy)$;
the two-dimensional $L^2$ change of variable costs $v^{-1}$,
cancelled by the factor $v$ in the integral. The coefficients are
smooth even functions on this ball. For the scalar annular functional,
first change variable to $x=at$; the integrable endpoint weight
$x/\sqrt{1-x^2}$ is fixed, and parameter derivatives differentiate
$g(x/a)$. The one-dimensional Sobolev estimate away from zero
controls these derivatives with the reserve in \eqref{eq:green-loss}.

\paragraph{The exterior and continuity in $a$.}
Use \eqref{eq:sonic-exterior}, separating the sonic annulus first. For the portion separated from the sonic annulus,
Cauchy--Schwarz gives explicitly
\[
 \begin{aligned}
 \left|\int_{2r_a}^r\frac{tg(t)}{\sqrt{a^2t^2-1}}\dd t\right|
 &\le C_I\int_{2r_a}^r|g(t)|\dd t\\
 &\le C_I\left(\int_{2r_a}^r t|g(t)|^2\dd t\right)^{1/2}
              \left(\int_{2r_a}^r\frac{\dd t}{t}\right)^{1/2}
 \le C_I\sqrt{\log(2+r)}\|g\|_{L^2(r\dd r)}.
 \end{aligned}
\]
Its coefficient is $O_I(r^{-3})$, and
$\int_{2r_a}^\infty r^{-6}\log(2+r)r\dd r<\infty$. Spatial and parameter derivatives have the same
integrable type, with the indicated additional derivatives of the source.
The cutoff derivatives are lower-order terms on the overlap regions.
This proves \eqref{eq:green-loss}. Smooth approximation and this
uniform estimate prove strong continuity. The mean-value formula
with one extra derivative of the source proves operator-norm continuity.
\end{proof}

\subsection{Parameter regularity of the profile branch}

The profile construction and its radial nondegeneracy are those
of \cite[Proposition~4.1 and Theorem~4.2]{GM1}. The proof of
Proposition~\ref{prop:slope} writes its fixed-point map and
Neumann-series inverse explicitly. Here we justify differentiation
on positive parameter intervals and the higher spatial bounds.
The stationary map is
\begin{equation}\label{eq:stationary-map}
 \mathcal S(a,P)=-\Delta P+P+P\mathcal R_a(P^2).
\end{equation}
Lemma~\ref{lem:green} makes it $C^1$ from
$I\times H^{m+2}_\rad$ to $H^m_\rad$. Continuity of its scalar
parameter derivative follows from strong continuity at a fixed
argument and the uniform bound for varying arguments.
Starting from $P_a\in H^2$, if $P_a\in H^{m+1}$ with $m\ge1$,
then $\mathcal R_a(P_a^2)\in H^m$ by \eqref{eq:green-loss}.
Multiplication by $P_a$ preserves $H^m$, so the stationary
equation gives $P_a\in H^{m+2}$. This places $P_a$ in the
domain of the stationary map at every finite level.

Its $P$ derivative is
\begin{equation}\label{eq:stationary-linearization}
 \mathcal L_a\xi=L_a\xi+2P_a\mathcal R_a(P_a\xi).
\end{equation}
The factorization proved in Proposition~\ref{prop:slope} gives
\[
 \|\mathcal L_a^{-1}\|_{L^2_\rad\to H^2_\rad}
       \le2\|L_+^{-1}\|_{L^2_\rad\to H^2_\rad}.
\]
This supplies the base $H^2$ estimate for the higher regularity argument.
The nonlocal second term loses at most one derivative. Thus for
$F\in H^m$, the equation $\mathcal L_a\xi=F$ raises the base
$H^2$ regularity one derivative at a time, up to $H^{m+2}$.
At that level, elliptic regularity and Lemma~\ref{lem:green} give
\[
 \begin{aligned}
 \|\xi\|_{H^{m+2}}
 &\le C_{I,m}(\|F\|_{H^m}
           +\|P_a\mathcal R_a(P_a\xi)\|_{H^m}+\|\xi\|_2)\\
 &\le C_{I,m}(\|F\|_{H^m}+\|\xi\|_{H^{m+1}})\\
 &\le\tfrac12\|\xi\|_{H^{m+2}}
                    +C'_{I,m}(\|F\|_{H^m}+\|\xi\|_2).
 \end{aligned}
\]
Absorb the first term and use the already established $L^2\to H^2$
inverse to bound $\|\xi\|_2\le C_I\|F\|_2$. This proves
$\|\mathcal L_a^{-1}F\|_{H^{m+2}}\le C_{I,m}\|F\|_{H^m}$ at
each finite level.
The $C^1$ implicit-function theorem at each level identifies the
same branch by its $H^2$ uniqueness.

Higher parameter derivatives follow by difference quotients in higher
Sobolev spaces. Differentiation once gives
$\mathcal L_a\partial_aP_a=-P_a(\partial_a\mathcal R_a)(P_a^2)$,
where the derivative of the operator holds its argument fixed.
Subtract this equation at $a+h$ and $a$, divide by $h$, and apply
the fixed inverse at $a$. The exact difference equation is
\[
 \begin{aligned}
 \mathcal L_a\frac{\partial_aP_{a+h}-\partial_aP_a}{h}
 ={}&-\frac{P_{a+h}(\partial_a\mathcal R_{a+h})(P_{a+h}^2)
                  -P_a(\partial_a\mathcal R_a)(P_a^2)}{h}\\
 &-\frac{\mathcal L_{a+h}-\mathcal L_a}{h}\partial_aP_{a+h}.
 \end{aligned}
\]
Here the parameter derivative in $\partial_a\mathcal R_{a+h}$
is evaluated at $a+h$, with the source held fixed.
Use \eqref{eq:green-loss} with the additional derivatives of the source
required by the operator quotient, and its next parameter order
for the forcing quotient. Both then converge in $H^m$. The fixed inverse then gives
\[
 \partial_a^2P_a
 =-\mathcal L_a^{-1}\bigl(
       \partial_a\{P_a(\partial_a\mathcal R_a)(P_a^2)\}
       +(\partial_a\mathcal L_a)\partial_aP_a\bigr)
       \quad\hbox{in }H^{m+2}.
\]
Repeat at each higher Sobolev level, taking enough input
derivatives for the next difference quotient. This proves every
finite parameter order on a positive compact interval.

The differentiated Schr\"odinger equation has an exponentially
localized right side, because every wave output is multiplied by a
profile derivative. Weighted elliptic estimates give
\eqref{eq:P-tail} for the parameter derivatives.
In the exterior Green formula, its limiting coefficient can be
written with the fixed lower endpoint $1$ using $z=at$.
Differentiation then preserves the $r^{-3}$ tail; spatial derivatives
give further inverse powers. This proves \eqref{eq:N-tail} and the
corresponding $r^{-2}$ velocity tails on positive compact intervals.
The uniform spatial bounds and convergence at zero come from
\cite[Appendix B, Propositions B.2 and B.7, Lemmas B.3--B.6]{GM2}.
The exponential rate in those bounds may be increased to any fixed
$d<1$: choose an exterior radius where $\abs{N_a}<1-d_1^2$
uniformly, with $d<d_1<1$, and compare the positive solution with
$\e^{-d_1r}$. Elliptic estimates then give spatial derivatives.
The differentiated equations and their localized right sides give
the same assertion for the positive-parameter derivatives on $I$.

\subsection{Estimates for the stationary wave defect}

The next lemma controls both the stationary wave expansion and
its parameter derivative. In $\partial_a\mathcal R_a$, the input
is held fixed during this operator differentiation.

\begin{lemma}[Small-parameter Green expansions]
\label{lem:green-expansion}
Let $0<a<a_*$, and let real radial functions $g_a,p_a$ satisfy
$\abs{\partial_r^j g_a(r)}\le C\e^{-dr}$ for $0\le j\le4$ and
$\abs{p_a(r)}\le C\e^{-dr}$, with $C,d>0$ independent of $a$.
Then
\begin{equation}\label{eq:green-undifferentiated}
\begin{aligned}
 &\norm{\mathcal R_ag_a+g_a-a^2\Xi[g_a]}_2\\
 &\qquad+\norm{\Dtwo(\mathcal R_ag_a+g_a-a^2\Xi[g_a])}_2
 \le C'a^4,
\end{aligned}
\end{equation}
and the localized parameter derivative satisfies
\begin{equation}\label{eq:green-differentiated}
 \left\|p_a\left[\frac{(\partial_a\mathcal R_a)g_a}{2a}
                                      -\Xi[g_a]\right]\right\|_2
 \le C' a^2.
\end{equation}
The estimate involves only the spatial norms of $g_a$ and $p_a$.
\end{lemma}

\begin{proof}
On $r\le(2a)^{-1}$, formula \eqref{eq:sonic-interior} is
\[
 \mathcal R_ag=-\frac{g}{1-a^2r^2}
       +\frac{a^2 I_a(r)}{(1-a^2r^2)^{3/2}},\qquad
 I_a(r)=\int_r^{1/a}\frac{t g(t)}{\sqrt{1-a^2t^2}}\dd t.
\]
To differentiate the integral, separate the region $at\le3/4$ with
an overlapping smooth cutoff. On this region expand the denominator;
on its complement use $x=at$, with fixed endpoint $x=1$.
The exponential bounds on $g,g'$ give
\begin{align*}
 I_a(r)&=\int_r^\infty tg(t)\dd t
          +O\bigl(a^2(1+r)^3\e^{-dr}\bigr)
          +O(a^{-K}\e^{-d/(2a)}),\\
 (\partial_a I_a)(r)&=O\bigl(a(1+r)^3\e^{-dr}\bigr)
          +O(a^{-K}\e^{-d/(2a)}).
\end{align*}
Here $K$ is a fixed finite integer. In the second estimate the
singular endpoint is never differentiated: $g(x/a)$ is differentiated
under the fixed integrable weight $x/\sqrt{1-x^2}$. Cutoff derivatives
are supported at $t$ comparable to $a^{-1}$ and have exponentially
small bounds. Differentiating the displayed expression for
$\mathcal R_ag$ consequently gives
\[
 \left|\frac{(\partial_a\mathcal R_a)g}{2a}-\Xi[g]\right|
 \le C a^2(1+r)^4\e^{-dr}+C a^{-K-2}\e^{-d/(2a)}
 \quad\text{on }r\le(2a)^{-1}.
\]
The first term is square integrable uniformly after multiplication
by $p_a$; the second is smaller than every power of $a$.
Before differentiation in $a$, the two core remainders are
\[
 \begin{aligned}
 \mathcal R_ag+g-a^2\Xi[g]
 ={}&-\frac{a^4r^4g}{1-a^2r^2}
       +a^2\{(1-a^2r^2)^{-3/2}-1\}\int_r^\infty tg(t)\dd t\\
 &+a^2(1-a^2r^2)^{-3/2}
                \left(I_a(r)-\int_r^\infty tg(t)\dd t\right).
 \end{aligned}
\]
Here $ar\le1/2$ and
$|(1-a^2r^2)^{-3/2}-1|\le Ca^2r^2$.
Also $I_a'(r)=-rg(r)(1-a^2r^2)^{-1/2}$, so differentiating this
identity once in $r$ gives the same estimate with one additional
polynomial factor. Consequently
\[
 \abs{\mathcal R_ag+g-a^2\Xi[g]}
 +\abs{r\partial_r(\mathcal R_ag+g-a^2\Xi[g])}
 \le Ca^4(1+r)^6\e^{-dr}+Ca^{-K-2}\e^{-d/(2a)}
\]
on this core region. Its $L^2$ norm is $O(a^4)$.

On $r\ge(2a)^{-1}$ use the nonsingular average
\eqref{eq:sonic-average}. Every sample point $s_t$ is at least
$(2a)^{-1}$, and $s_t+r_a$ is bounded by a fixed multiple of
$r+r_a$. Also $\partial_as_t=-a^{-2}(1-t)$.
Direct differentiation of that average yields, after increasing a
fixed power of $a^{-1}$ if necessary,
\[
 |(\partial_a\mathcal R_a)g(r)|
       \le C a^{-4}(r+r_a)^{-1}\e^{-d/(2a)}.
\]
The analogous bound for $\Xi[g]$ follows immediately from its
integral. Multiplication by $p_a$ makes the exterior $L^2$ norm
smaller than every power of $a$. This proves
\eqref{eq:green-differentiated} on the whole plane.

For \eqref{eq:green-undifferentiated} the unweighted exterior norm
needs a further split. On $(2a)^{-1}\le r\le2/a$, the average
formula and one spatial derivative give an $L^2$ bound
$Ca^{-3}\e^{-d/(2a)}$ for $\mathcal R_ag$ and its dilation.
For $r\ge2/a$, the integral in \eqref{eq:sonic-exterior} is bounded by
\[
 \int_{1/a}^\infty\frac{t\abs{g(t)}}{\sqrt{a^2t^2-1}}\dd t
 \le Ca^{-3/2}\e^{-d/a}.
\]
Indeed set $t=a^{-1}+x$ and bound the denominator below by
$\sqrt{2ax}$. Thus
\[
 \abs{\mathcal R_ag(r)}
 \le Ca^{-2}r^{-2}\e^{-dr}
                +Ca^{-5/2}r^{-3}\e^{-d/a}.
\]
The radial equation \eqref{eq:sonic-ode} implies
$\abs{r\partial_r\mathcal R_ag}\le C\abs{\mathcal R_ag}
 +Ca^{-2}r^{-1}\e^{-dr}$ there. These bounds are square integrable
and smaller than every power of $a$ on $r\ge2/a$.
The corresponding tails of $g$ and $\Xi[g]$ are exponential.
Together with the core estimate, this proves
\eqref{eq:green-undifferentiated}.
\end{proof}

\begin{corollary}[The stationary wave defect]\label{cor:profile-defect}
For all sufficiently small $a>0$,
\begin{equation}\label{eq:profile-defect-sharp}
 \norm{N_a+P_a^2}_2+\norm{\Dtwo(N_a+P_a^2)}_2\le Ca^2.
\end{equation}
Moreover $a^{-2}(N_a+P_a^2)\to\Xi[Q^2]$ in the norm
$\norm g_2+\norm{\Dtwo g}_2$.
\end{corollary}

\begin{proof}
Apply \eqref{eq:green-undifferentiated} to $g_a=P_a^2$, using
$N_a=\mathcal R_a(P_a^2)$. By the explicit formula for $\Xi$,
\[
 |\Xi[g](r)|\le\int_r^\infty t|g(t)|\dd t+r^2|g(r)|,
 \qquad r\Xi[g]'=-r^2(3g+rg').
\]
The uniform exponential bounds therefore dominate $\Xi[P_a^2]$
and $\Dtwo\Xi[P_a^2]$ by an $L^2$ function independent of $a$.
Dominated convergence gives their respective limits with $Q^2$.
Finally,
\[
 a^{-2}(N_a+P_a^2)
 =\Xi[P_a^2]+a^{-2}
       \{\mathcal R_a(P_a^2)+P_a^2-a^2\Xi[P_a^2]\},
\]
and the last term has norm $O(a^2)$, also after $\Dtwo$. This proves
both claims.
\end{proof}

\subsection{Localized estimates and ground-state perturbation}\label{sec:spectral}

The following estimates let local errors act on exponentially
decaying tests and integrable polynomial tails. They are the reason
all source pairings except the cubic interaction can be
controlled without global Schr\"odinger mass.

\begin{lemma}[Local estimates]\label{lem:anchoring}
The weighted estimates and product tests below hold on $X$ with or
without radial symmetry. For radial ${\eps}\in X$, for almost every $r>0$,
\begin{equation}\label{eq:log-bound}
 \abs{{\eps}(r)}^2\le C(1+\abs{\log r})\norm {\eps}_{\mathrm{loc}}^2.
\end{equation}
For arbitrary ${\eps}\in X$, $d>0$ and $p>2$,
\begin{equation}\label{eq:weighted-local}
 \norm{\e^{-dr}{\eps}}_{H^1}\le C_d\norm {\eps}_{\mathrm{loc}},\qquad
 \int\langle r\rangle^{-p}\abs {\eps}^2\le C_p\norm {\eps}_{\mathrm{loc}}^2.
\end{equation}
If $\pair{\eps}Q=0$, then also
\begin{equation}\label{eq:constrained-anchoring}
 \norm {\eps}_{\mathrm{loc}}\le C\norm{\nabla {\eps}}_2.
\end{equation}
If $p_{\eps}$ is an exponentially localized $H^1$ scalar test, then
\begin{equation}\label{eq:fh-test}
 \abs{\pair{{\epsn}{\eps}}{p_{\eps}}}
 \le C_d\norm{\epsn}_2\norm {\eps}_{\mathrm{loc}}\norm{\e^{dr}p_{\eps}}_{H^1}.
\end{equation}
If $p_v$ is an exponentially localized $H^1$ vector test, or is smooth
with $\abs{p_v}\le C\langle r\rangle^{-2}$ and
$\abs{\diver p_v}\le C\langle r\rangle^{-3}$, then
\begin{equation}\label{eq:h2-test}
 \abs{\pair{\nabla\abs {\eps}^2}{p_v}}\le C_{p_v}\norm {\eps}_{\mathrm{loc}}^2.
\end{equation}
\end{lemma}

\begin{proof}
\paragraph{Radial and angular control.}
The weighted local $L^2$ norm provides an anchor point in $1<r<2$. Radial Cauchy--Schwarz gives
\[
 \abs{{\eps}(r)-{\eps}(t)}^2\le(2\pi)^{-1}\norm{\nabla {\eps}}_2^2
                              \abs{\log(r/t)}.
\]
This proves \eqref{eq:log-bound}. Integration gives the radial
case of \eqref{eq:weighted-local}.
For a general function, let
$\bar\eps(r)=(2\pi)^{-1}\int_0^{2\pi}\eps(r,\theta)\dd\theta$.
Its local norm is bounded by that of $\eps$, and it satisfies the
radial estimate. On each circle,
the Poincar\'e inequality gives
\[
 \int_0^{2\pi}|\eps-\bar\eps|^2\dd\theta
       \le\int_0^{2\pi}|\partial_\theta\eps|^2\dd\theta,
 \qquad
 \int\langle r\rangle^{-p}|\eps-\bar\eps|^2
       \le C_p\norm{\nabla\eps}_2^2.
\]
Indeed, with $h=\eps-\bar\eps$,
\[
 \begin{aligned}
 \int\langle r\rangle^{-p}|h|^2
 &\le\sup_{r>0}\{r^2\langle r\rangle^{-p}\}
       \int_0^\infty\!\int_0^{2\pi}
                         r^{-2}|\partial_\theta h|^2r\dd\theta\dd r
 \le C_p\|\nabla\eps\|_2^2,\\
 \|\nabla(e^{-dr}\eps)\|_2
 &\le\|\nabla\eps\|_2+C_d\|e^{-dr}\eps\|_2.
 \end{aligned}
\]
The smoothed exponential weight from Section~\ref{sec:notation}
is used in the second line. Together with the radial mean estimate,
these inequalities prove \eqref{eq:weighted-local}.
Under $\pair{\eps}Q=0$, the angular mean has the same zero moment.
For that mean, integrate
$\bar\eps(r)-\bar\eps(1)$ against $Q$ in the displayed Cauchy--Schwarz estimate.
Since $\int Q>0$ and $\int Q\sqrt{\abs{\log r}}<\infty$, this gives
$|\bar\eps(1)|\le C\norm{\nabla\eps}_2$. Together with the
mean-zero estimate, this proves \eqref{eq:constrained-anchoring}.
\paragraph{Product tests.}
To prove \eqref{eq:fh-test}, write
${\eps}p_{\eps}=(\e^{-dr}{\eps})(\e^{dr}p_{\eps})$ and use $H^1(\R^2)\hookrightarrow L^4$.
For an exponential vector test,
\[
 |\pair{\nabla|\eps|^2}{p_v}|
 =\left|\int|\eps|^2\diver\bar p_v\right|
 \le\|e^{-dr/2}\eps\|_4^2\|e^{dr}\diver p_v\|_2
 \le C_d\|\eps\|_{\mathrm{loc}}^2\|p_v\|_{H^1_d}.
\]
For a polynomial-tail test the estimate is instead
\[
 |\pair{\nabla|\eps|^2}{p_v}|
 \le C\int\langle r\rangle^{-3}|\eps|^2
 \le C\|\eps\|_{\mathrm{loc}}^2 .
\]
With a smooth cutoff at radius $R$,
the error is bounded by
$C\int_{R<r<2R}r^{-3}|\eps|^2$, which tends to zero. This also
defines the distributional pairing for general $\eps$.

\paragraph{Approximation in $X$.}
For density of compactly supported functions, first consider the
radial part. Its finite tail gradient implies
$\eps(r)=o(\sqrt{\log r})$ as $r\to\infty$: anchor at
an arbitrarily large fixed radius and apply Cauchy--Schwarz to the
remaining, arbitrarily small, gradient tail. Choose a radial cutoff
equal to one on $r\le R$, zero on $r\ge R^2$, with derivative bounded
by $C/(r\log R)$. Its extra gradient has squared norm bounded by
\[
 \frac C{(\log R)^2}\int_R^{R^2}\frac{\abs{{\eps}(r)}^2}{r}\dd r=o(1).
\]
The other cutoff errors tend to zero by integrability. For the
mean-zero part, the angular Poincar\'e inequality gives
$\int r^{-2}|\eps-\bar\eps|^2\le\norm{\nabla\eps}_2^2$; it makes the
extra cutoff gradient tend to zero as well. Mollification proves
density of smooth compactly supported functions in $X$, preserving
radial symmetry when present. Fixed compactly supported tests correct
any finitely many independent profile moments, with coefficients
tending to zero.
\end{proof}

We use the weighted spaces and the reserve exponents fixed in
Section~\ref{sec:notation}.

\begin{lemma}[The reduced inverse on $X$]\label{lem:RX}
With or without radial symmetry, $R_a$ extends to a bounded map
$X\to X$, uniformly for $a\in I$. For a smooth exponentially localized multiplier $p$,
\begin{equation}\label{eq:local-R-regularity}
 \norm{pR_a{\eps}}_{H^3_d}+\norm{pR_a^2{\eps}}_{H^5_d}
 \le C_{I,p,d}\norm {\eps}_{\mathrm{loc}},
\end{equation}
provided the decay of $p$ is faster than $d$. Parameter derivatives of
the finite orders used below have the same mapping properties.
\end{lemma}

\begin{proof}
We first bound the free inverse on $X$, then treat $N_a$ as a
compact perturbation. Let $J$ be the Bessel kernel of
$(1-\Delta)^{-1}$. Its origin singularity is logarithmic and its
tail is exponential, so weighted Young's inequality gives
\[
 \begin{aligned}
 \norm{\e^{-r/2}(1-\Delta)^{-1}h}_2
 &\le\norm{\e^{r/2}J}_{L^1}\norm{\e^{-r/2}h}_2,\\
 \norm{\nabla(1-\Delta)^{-1}h}_2&\le\norm{\nabla h}_2.
 \end{aligned}
\]
The kernel norm is finite, and the last estimate follows by commuting
the gradient with $(1-\Delta)^{-1}$.
Lemma~\ref{lem:anchoring} gives $(1-\Delta)^{-1}:X\to X$.

\paragraph{The profile perturbation.}
Multiplication $N_a:X\to L^2$ is compact. On a fixed ball this follows
from Rellich compactness. Outside that ball the operator norm tends to
zero by $\abs{N_a}\le C\langle r\rangle^{-3}$ and
\eqref{eq:weighted-local}, with any exponent between $2$ and $6$.
The equation for the extension of $R_a{\eps}$ is
\begin{equation}\label{eq:fredholm}
 \begin{aligned}
 (1+(1-\Delta)^{-1}N_a)z
 &=(1-\Delta)^{-1}\left({\eps}-\frac{\pair {\eps}{P_a}}M P_a\right),\\
 \pair z{P_a}&=0.
 \end{aligned}
\end{equation}
The homogeneous kernel is $\C P_a$. Indeed, if $L_az=0$ in $X$,
then $N_az\in L^2$ by the preceding multiplication bound. The
anchoring estimate makes $z$ a tempered distribution, so
$(1-\Delta)z=-N_az$ implies $z=-(1-\Delta)^{-1}(N_az)\in H^2$.
The ground-state transform then gives $z\in\C P_a$.

The range of $1+(1-\Delta)^{-1}N_a$ is the kernel of the continuous functional
$v\mapsto\pair v{(1-\Delta)P_a}$. Its value on $(1-\Delta)^{-1}{\eps}$ is
$\pair {\eps}{P_a}$, and $(1-\Delta)P_a=-N_aP_a$ is localized.
To verify the range assertion, for $z\in X$,
\[
 \begin{aligned}
 \pair{(1+(1-\Delta)^{-1}N_a)z}{(1-\Delta)P_a}
 &=\pair z{(1-\Delta)P_a}+\pair{N_az}{P_a}
   =\pair z{L_aP_a}=0,\\
 \pair{(1-\Delta)^{-1}(\eps-M^{-1}\pair\eps{P_a}P_a)}
       {(1-\Delta)P_a}&=0.
 \end{aligned}
\]
The compact perturbation of identity has Fredholm index zero and
one-dimensional kernel, so this annihilating functional describes
its entire range. The normalization $\pair z{P_a}=0$ removes the kernel.
The bounded inverse theorem then gives $\|R_a\eps\|_X\le C_I\|\eps\|_X$. Elliptic estimates on unit balls, summed
with a decreasing exponential weight, give local $H^3$ control from
the weighted $H^1$ bounds for the source and the solution. Multiplication
by $p$ gives the first estimate in \eqref{eq:local-R-regularity}.
For $R_a^2\eps$, apply the same argument with source $R_a\eps$,
whose weighted local $H^3$ norm was just bounded; elliptic regularity
then gives $H^5$. The spare decay of $p$ absorbs the decreasing
weight in both estimates.

Differentiate the equation \eqref{eq:fredholm} in $a$.
The derivatives of $N_a$ have the same compact-tail bounds, while
those of $P_a$ are localized. The inverse remains bounded on the
compact interval $I$. This proves the asserted parameter bounds.
\end{proof}

Only in the following eigenpair calculation, let $S_N$ denote
the reduced inverse of $-\Delta+N-\mu(N)$ on $\phi(N)^\perp$.
It differs from the inverse $R_N$ of $-\Delta+1+N$ used in
Section~\ref{sec:nonlinear-cancellation}. Analytic dependence of a simple isolated
eigenpair is classical \cite[Chapter~VII]{Kato}. The next lemma
also gives the weighted estimates and the distributional direction
needed for the nonlinear proof, which we verify below.

\begin{lemma}[Derivatives of the ground state]\label{lem:spectral}
On a uniform $L^2$ neighborhood of $N_a$, the eigenpair is real analytic
with $\phi(N)\in H^2_{d'}$, after shrinking the neighborhood so that
$(d')^2<-\mu(N)$. For real ${\epsn}$,
\begin{equation}\label{eq:phi-derivative}
 \phi'_N[{\epsn}]=-S_N({\epsn}\phi),\qquad \mu'_N[{\epsn}]=({\epsn},\phi^2).
\end{equation}
There is a continuous distributional extension for one velocity
divergence. In the next two estimates, $f_i$ denote arbitrary real $L^2$
density directions. For every fixed finite $j$,
\begin{align}
 \norm{\phi_N^{(j)}[f_1,\ldots,f_j]}_{H^2_d}
 &\le C_{I,j}\prod_{i=1}^j\norm{f_i}_2,\label{eq:pure-spectral}\\
 \norm{\phi_N^{(j+1)}[f_1,\ldots,f_j,\diver {\epsv}]}_{H^1_d}
 &\le C_{I,j}\norm {\epsv}_2\prod_{i=1}^j\norm{f_i}_2.
 \label{eq:mixed-spectral}
\end{align}
The mixed differential in \eqref{eq:mixed-spectral} is defined by
the differentiated eigenvalue equation with its normalization and
an $H^{-1}$ source.
\end{lemma}

\begin{proof}
The normalization fixes the ground-state component, while the
eigenvalue variable fills the missing range direction. Apply the
implicit-function theorem to
\[
 (\phi,\mu)\longmapsto
 \bigl((-\Delta+N-\mu)\phi,\tfrac12(\norm\phi_2^2-1)\bigr).
\]
At a normalized eigenpair its derivative is
\[
 (\eta,\nu)\longmapsto
       ((-\Delta+N-\mu)\eta-\nu\phi,(\phi,\eta)).
\]
For a right-hand side $(F,b)\in L^2\times\R$, the inverse is
explicitly
\[
 \nu=-(F,\phi),\qquad \eta=S_NF+b\phi.
\]
The inverse at \eqref{eq:phi0} is bounded $L^2\times\R\to H^2\times\R$
by the simple ground-state gap. Normalization fixes the kernel
component, and $\nu$ supplies the missing range direction. The product
$L^2\times H^2\to L^2$ is bounded in dimension two. The analytic
implicit-function theorem gives the eigenpair and
\eqref{eq:phi-derivative}. For the weighted implicit-function argument at $N_a$, use
$S_{N_a}=R_a$. The exterior estimate in the proof of
Lemma~\ref{lem:current-inverse} for this reference operator
uses only the gap and the tail of $N_a$. It gives the normalized
inverse from $L^2_{d'}\times\R$ to $H^2_{d'}\times\R$.
Multiplication by an $L^2$ potential is bounded from
$H^2_{d'}$ to $L^2_{d'}$, so the weighted implicit-function
theorem gives the same eigenpair by uniqueness. For the
weighted form inverse at nearby $N$, add the projection onto
$\phi(N)$ to $-\Delta+N-\mu(N)$. This operator is a small
perturbation of $L_a+1-\Pi_a$ from $H^1_{d'}$ to $H^{-1}_{d'}$,
using the weighted eigenfunction bound and
$L^2\times H^1\to H^{-1}$. Its inverse on the complementary
projection gives $S_N:H^{-1}_{d'}\to H^1_{d'}$.

\paragraph{Higher density derivatives.}
For real $f,g\in L^2$, differentiating the eigenvalue equation twice gives
\[
 \begin{aligned}
 (-\Delta+N-\mu(N))\phi''_N[f,g]
 ={}&-(f-\mu'_N[f])\phi'_N[g]
     -(g-\mu'_N[g])\phi'_N[f]\\
    &+\mu''_N[f,g]\phi,\\
 (\phi,\phi''_N[f,g])={}&-(\phi'_N[f],\phi'_N[g]),\qquad
 \mu''_N[f,g]=2(f\phi,\phi'_N[g]).
 \end{aligned}
\]
The normalization specifies the ground-state component; the reduced
inverse solves the orthogonal component. Pure density derivatives
are controlled in $H^2_d$: each new $L^2$ direction multiplies a
weighted $H^2$ output, which is bounded in dimension two and hence
gives an $L^2_d$ source.

\paragraph{One distributional direction.}
For a divergence direction, write
\[
 \phi\diver w=\diver(\phi w)-\nabla\phi\cdot w.
\]
For an $H^1$ test the first term is controlled by
$\norm\phi_\infty\norm w_2$, and the second by
$\norm{\nabla\phi}_4\norm w_2$, using $H^1\hookrightarrow L^4$.
The same estimates with exponential weights put the source in
$H^{-1}_d$, meaning that the source multiplied by $e^{dr}$
belongs to $H^{-1}$. For a pure-density output $g\in H^2_{d'}$ and a test $p$ with
$\|p\|_{H^1_{-d}}:=\|e^{-dr}p\|_{H^1}<\infty$, the source
estimate is
\[
 \begin{aligned}
 |\langle g\diver w,p\rangle|
 &=|(w,g\nabla p+p\nabla g)|\\
 &\le\|w\|_2
   \bigl(\|e^{dr}g\|_\infty\|e^{-dr}\nabla p\|_2
         +\|e^{dr}\nabla g\|_4\|e^{-dr}p\|_4\bigr)\\
 &\le C_d\|w\|_2\|g\|_{H^2_{d'}}\|p\|_{H^1_{-d}}.
 \end{aligned}
\]
The reserve $d'>d$ controls the weight derivatives. Thus the
source lies in $H^{-1}_d$, and the normalized inverse maps it to
$H^1_d$ by the weighted elliptic form estimate.
Its eigenvalue component is also defined distributionally:
\[
 \mu'_N[\diver w]=-(w,\nabla\phi^2),\qquad
 |\mu'_N[\diver w]|\le\|w\|_2\|\nabla\phi^2\|_2.
\]

The displayed second-order recursion remains valid with exactly one
divergence argument, interpreted this way. Every further density
differentiation produces either a pure-density term already in
$H^2_d$, or an $L^2$ direction multiplying a mixed $H^1_d$ output.
The latter product lies in $H^{-1}_d$, since
\[
 |\pair{fh}p|\le C\norm f_2\norm h_{H^1}\norm p_{H^1}.
\]
If the divergence acts on a higher pure-density output, that output
is in $H^2_d$ and the same divergence decomposition applies.
Normalization and eigenvalue terms contain at most one mixed
factor and satisfy the same bounds. Induction in the number of
density directions proves \eqref{eq:pure-spectral} and
\eqref{eq:mixed-spectral}, with constants depending on the finite
order and $I$, uniformly in the scale.
\end{proof}

\subsection{Reduced inverse of the density-dependent Schr\"odinger operator}
\label{app:current-inverse}

\begin{lemma}[The inverse and its differentiated tests]
\label{lem:current-inverse}
For $N$ in a sufficiently small $L^2$ neighborhood of $N_a$,
$a\in I$, the operators in \eqref{eq:current-inverse} are
well defined and depend smoothly on $N$. They map localized
$L^2$ sources to $H^2_d$ and localized $H^{-1}$ sources to
$H^1_d$, uniformly on this neighborhood. For fixed $g$ and a real $L^2$ direction $h$,
\begin{equation}\label{eq:inverse-derivative}
 \begin{aligned}
 R'_N[h]g={}&-R_N(hR_Ng)
       -\pair g\phi R_N\phi'_N[h]
       -\pair{R_Ng}{\phi'_N[h]}\phi .
 \end{aligned}
\end{equation}
For $h=\diver w$, $w\in L^2$, this is a continuous tested
identity with values in $H^1_d$. For $0\le j\le3$, real
$f_i\in L^2$, and localized $g\in L^2_{d'}$,
\begin{equation}\label{eq:inverse-mixed-tests}
 \begin{aligned}
 \|R_N^{(j)}[f_1,\ldots,f_j]g\|_{H^2_d}
       &\le C_{I,j}\|g\|_{L^2_{d'}}\prod_i\|f_i\|_2,\\
 \|R_N^{(j)}[f_1,\ldots,f_j]g\|_{H^1_d}
       &\le C_{I,j}\|g\|_{H^{-1}_{d'}}\prod_i\|f_i\|_2,\\
 \|R_N^{(j+1)}[f_1,\ldots,f_j,\diver w]g\|_{H^1_d}
       &\le C_{I,j}\|w\|_2\|g\|_{L^2_{d'}}\prod_i\|f_i\|_2 .
 \end{aligned}
\end{equation}
The reserve exponential weight also permits the polynomial
factors in the parameter and dilation identities below.
\end{lemma}

\begin{proof}
At $N=N_a$, $L_a$ has a simple zero eigenvalue and a positive
gap on $P_a^\perp$. The augmented operator
\[
 L_N+\pair{\,\cdot\,}\phi\phi
\]
has an inverse near $N_a$. Indeed multiplication by $N-N_a$
maps $H^2$ to $L^2$ with norm at most $C\|N-N_a\|_2$,
and $\phi(N)$ is analytic by Lemma~\ref{lem:spectral}.
The Neumann series about the augmented operator at $N_a$
therefore converges. Its inverse followed by $\Pi_N$ is $R_N$.
For the weighted estimates the reference inverse is explicitly
\[
 (L_a+1-\Pi_a)^{-1}=R_a+1-\Pi_a.
\]
For $u=R_a g$, the gap and the bounded potential give
$\|u\|_{H^1}\le C_I\|g\|_{H^{-1}}$. Test $L_au=\Pi_ag$ with
$w^2u$, where $w$ is a smooth bounded truncation of $\e^{dr}$
with $|\nabla\log w|\le d$ outside a fixed ball. Taking the real
part gives
\[
 \|\nabla(wu)\|_2^2+
 \int (1+N_a-|\nabla\log w|^2)|wu|^2
      =\Rea\pair{\Pi_ag}{w^2u}.
\]
Outside a fixed larger ball, the coefficient is at least
$(1-d^2)/2$, uniformly on $I$. The unweighted bound controls
the interior term. Since $w\e^{-dr}$ and its first derivatives
are uniformly bounded, the right-hand side is at most
$C_I\|g\|_{H^{-1}_d}\|wu\|_{H^1}$, using the weighted bound on
the projection onto $P_a$. Absorption and removal of the truncation
give $R_a:H^{-1}_d\to H^1_d$. Conjugating $L_a$ by $\e^{dr}$
and applying the $H^2$ elliptic estimate then give
$R_a:L^2_d\to H^2_d$. The argument also applies with $d'$
in place of $d$.

The perturbation is small in both source norms, since
\[
 \begin{aligned}
 \|\epsn u\|_{L^2_d}&\le C\|\epsn\|_2\|u\|_{H^2_d},\\
 \|\epsn u\|_{H^{-1}_d}&\le C\|\epsn\|_2\|u\|_{H^1_d}.
 \end{aligned}
\]
The difference of the two rank-one projections has the same
small operator bounds, by $\|\phi-\phi(N_a)\|_{H^2_d}
\le C_I\|\epsn\|_2$. Factoring the augmented operator by its
reference inverse therefore gives a convergent Neumann series
in each weighted space. Intermediate rates between $d$ and
$d'$ absorb the finitely many polynomial factors below.

Differentiate $L_NR_Ng=\Pi_Ng$ and
$\pair{R_Ng}\phi=0$. Projection onto $\phi^\perp$ and onto
$\C\phi$, respectively, gives
\[
 \begin{aligned}
 \Pi_NR'_N[h]g&=-R_N(hR_Ng)-\pair g\phi R_N\phi'_N[h],\\
 \pair{R'_N[h]g}\phi&=-\pair{R_Ng}{\phi'_N[h]}.
 \end{aligned}
\]
Their sum is \eqref{eq:inverse-derivative}.
The same inverse differentiation in the form spaces
$H^{-1}_{d'}\to H^1_d$, with multiplication
$L^2\times H^1\to H^{-1}$, proves the pure density
derivative bounds for negative-Sobolev sources as well.

For the distributional direction, if $u\in H^2_{d'}$ and
$\eta\in H^1$ is a test, integration by parts gives
\[
 \begin{aligned}
 |\pair{(\diver w)u}\eta|
 &=\left|\int w\cdot\nabla(u\bar\eta)\right|\\
 &\le\|w\|_2
       (\|\nabla u\|_4\|\eta\|_4
                       +\|u\|_\infty\|\nabla\eta\|_2)
 \le C\|w\|_2\|u\|_{H^2}\|\eta\|_{H^1}.
 \end{aligned}
\]
Inserting exponential weights gives the localized $H^{-1}$
estimate. Lemma~\ref{lem:spectral} bounds
$\phi'_N[\diver w]$ in $H^1_d$, so each term of
\eqref{eq:inverse-derivative} has the asserted regularity.
Repeated differentiation produces products of these inverses,
ground-state derivatives and multiplication by the $f_i$.
Before the divergence direction, the factors have $H^2$
regularity; after that direction, multiplication is estimated by
\[
 \|fu\|_{H^{-1}}\le C\|f\|_2\|u\|_{H^1}.
\]
Every subsequent inverse restores $H^1$ regularity. There is
exactly one divergence argument in each product. These bounds,
the pure ground-state bounds in \eqref{eq:pure-spectral} and
the mixed bounds in \eqref{eq:mixed-spectral}, prove
\eqref{eq:inverse-mixed-tests} term by term.
\end{proof}

\subsection{Distributional bounds for the linear correction}

\begin{lemma}[Localized inverse compositions]\label{lem:linear-tested}
For real radial $\epsn\in L^2$ and radial no-swirl
$\epsv\in L^2(\R^2;\R^2)$,
\begin{equation}\label{eq:linear-slow-tests}
 \begin{aligned}
 \|R_a\mathcal M_a\Dtwo\epsn\|_{H^1_d}
  +\|R_a\mathcal M_a\diver\epsv\|_{H^1_d}
      &\le C_I(\|\epsn\|_2+\|\epsv\|_2),\\
 \|R_a^2(P_a\diver\Dtwo\epsv)\|_{H^1_d}
  +\|R_a^2(P_a\Delta\epsn)\|_{H^1_d}
      &\le C_I(\|\epsn\|_2+\|\epsv\|_2).
 \end{aligned}
\end{equation}
The constants are uniform for $a\in I$. All sources are interpreted
as distributions.
\end{lemma}

\begin{proof}
Let $p$ be one of the smooth exponentially decreasing profile
multipliers occurring in $\mathcal M_a$. In distributions,
\[
 \begin{aligned}
 p\Dtwo\epsn
    &=\diver(yp\epsn)-(y\cdot\nabla p)\epsn,\\
 p\diver\epsv
    &=\diver(p\epsv)-\nabla p\cdot\epsv,\\
 p\Delta\epsn
    &=\Delta(p\epsn)-2\diver(\epsn\nabla p)
                         +(\Delta p)\epsn .
 \end{aligned}
\]
Thus the first two expressions are bounded in $H^{-1}_d$ and
the third in $H^{-2}_d$ by the corresponding $L^2$ error norms.
Here $H^{-j}_d$ consists of distributions whose product with the
fixed exponential weight belongs to $H^{-j}$. For a real test $\eta$,
two integrations by parts also give
\[
 |\pair{p\diver\Dtwo\epsv}{\eta}|
   =|(\epsv,y\cdot\nabla\nabla(p\overline\eta))|
   \le C_I\|\epsv\|_2\|\eta\|_{H^2}.
\]
The same estimate holds with the exponential weight inserted.
All powers of $y$ multiply localized profiles; the reserve
exponent $d'>d$ absorbs them.

Expansion of \eqref{eq:Fa} with
\eqref{eq:R-commutator}--\eqref{eq:K-commutator} gives
\[
 \begin{aligned}
 R_a\mathcal M_a={}&R_a[(2-\Lambda)P_a]-R_aA_aR_aP_a
                 +R_a^2A_aR_aP_a\\
 &+R_a^2[(y\cdot\nabla-1)P_a]
                 -\frac{(R_a\chi_a)\otimes P_a^2}{M}.
 \end{aligned}
\]
Profiles adjacent to an inverse denote multiplication; the last
term is a rank-one operator. For example,
$[(R_a\chi_a)\otimes P_a^2]\epsn
=(\epsn,P_a^2)R_a\chi_a$.

The cases $j=0,1$ of $R_a:H^{-j}_d\to H^{2-j}_d$ were proved
in Lemma~\ref{lem:current-inverse}. For $j=2$, self-adjointness
and the positive-weight $L^2_d\to L^2_d$ bound first give
$R_a:L^2_{-d}\to L^2_{-d}$. Conjugating $L_a R_a=\Pi_a$ by
$\e^{-dr}$ gives the $L^2_{-d}\to H^2_{-d}$ bound by elliptic
regularity, since $\Pi_a$ is bounded on $L^2_{-d}$ by the
exponential decay of $P_a$. Dualizing this last estimate gives
$R_a:H^{-2}_d\to L^2_d$. For a term with a second-order middle
operator,
\[
 P_a\Dtwo\epsn\in H^{-1}_d
 \xrightarrow{R_a}H^1_d
 \xrightarrow{A_a}H^{-1}_d
 \xrightarrow{R_a}H^1_d .
\]
All other summands in $R_a\mathcal M_a$ are no less regular,
and the rank-one term tests a distribution against a smooth
profile. This proves the first line of
\eqref{eq:linear-slow-tests}. For its second line, use
\[
 \begin{gathered}
 \|P_a\diver\Dtwo\epsv\|_{H^{-2}_d}
       +\|P_a\Delta\epsn\|_{H^{-2}_d}
            \le C_I(\|\epsn\|_2+\|\epsv\|_2),\\
 H^{-2}_d\xrightarrow{R_a}L^2_d
                \xrightarrow{R_a}H^2_d\hookrightarrow H^1_d .
 \end{gathered}
\]
These arguments first apply to smooth errors and extend by
density in $L^2$. They also apply without radial symmetry, using
the full-space reduced inverse.
\end{proof}

\subsection{Weighted estimates for the correction tests}\label{app:tests}\label{sec:form-domain}

\begin{lemma}[Derivatives of the quadratic correction]\label{lem:form-gradients}
The real gradients of $B_\beta$ have representations
\begin{equation}\label{eq:gradient-class}
 \begin{aligned}
 d_{\eps}B_\beta({\veps})&=\beta p_{\eps},\\
 d_{\epsn}B_\beta({\veps})
   &=\beta p_n-\frac{2\beta(\Ima\eps,\chi_a)}{M'}\partial_a(N_a+P_a^2),\\
 d_{\epsv}B_\beta({\veps})
   &=\beta p_v-\frac{2\beta(\Ima\eps,\chi_a)}{c_0^2M'}\partial_aV_a.
 \end{aligned}
\end{equation}
Here $p_{\eps},p_n,p_v$ have $H^1_d$ norms bounded by $C_I\norm{\veps}_Y$.
The remaining profile tests satisfy
\begin{equation}\label{eq:polynomial-tests}
\begin{aligned}
 |\partial_a(N_a+P_a^2)|+|y\cdot\nabla\partial_a(N_a+P_a^2)|
       &\le C_I\langle r\rangle^{-3},\\
 |\partial_aV_a|+|y\cdot\nabla\partial_aV_a|
       &\le C_I\langle r\rangle^{-2},\\
 |\diver\partial_aV_a|&\le C_I\langle r\rangle^{-3}.
\end{aligned}
\end{equation}
Their common coefficient satisfies
$|(\Ima\eps,\chi_a)/M'|\le C_I\|\veps\|_Y$.
The same estimates hold with a slightly larger exponential weight,
and after the finite parameter differentiations needed below.
All these gradients are linear in ${\veps}$.
\end{lemma}

\begin{corollary}\label{cor:nonlinear-pairings}
For $u\in Y$, real $f\in L^2$ and $h,h_1,h_2\in X$,
\begin{equation}\label{eq:polarized-cubic}
\begin{aligned}
 &\left|2\beta^{-1}B_\beta(u,(-\ii f h,0,-c_0^2\lambda\nabla\Rea(h_1\overline{h_2})))\right|\\
 &\qquad\le C_I\norm u_Y
       (\norm f_2\norm h_{\mathrm{loc}}
                     +\lambda\norm{h_1}_{\mathrm{loc}}\norm{h_2}_{\mathrm{loc}}).
\end{aligned}
\end{equation}
In particular, on energy data,
\begin{equation}\label{eq:B-cubic}
 \abs{2B_\beta({\veps},(-\ii {\epsn}{\eps},0,-c_0^2\lambda\nabla\abs {\eps}^2))}
 \le C_I \beta(1+\lambda)\norm{\veps}_Y^3.
\end{equation}
The linear wave and dilation sources admit the same tested
interpretation. The operator estimates for these gradients are
proved in Appendix~\ref{app:tests}.
\end{corollary}

\begin{proof}
For an exponentially localized $H^1$ scalar test $p$, weighted
$H^1\hookrightarrow L^4$ gives
\[
 |\pair{fh}p|\le
 C\norm f_2\norm{\e^{-dr}h}_{H^1}\norm p_{H^1_d}
 \le C\norm f_2\norm h_{\mathrm{loc}}\norm p_{H^1_d}.
\]
For a vector test $p$, the complete product estimate is
\[
 \begin{aligned}
 |\pair{\nabla\Rea(h_1\overline{h_2})}p|
 &=\left|\int\Rea(h_1\overline{h_2})\diver\bar p\right|\\
 &\le\|e^{-dr/2}h_1\|_4\|e^{-dr/2}h_2\|_4
                                      \|e^{dr}\diver p\|_2\\
 &\le C_d\|h_1\|_{\mathrm{loc}}\|h_2\|_{\mathrm{loc}}
                                      \|p\|_{H^1_d}.
 \end{aligned}
\]
The weighted $H^1$ bound of Lemma~\ref{lem:anchoring} and
$H^1(\R^2)\hookrightarrow L^4$ give the last line.
For the polynomial-tail test $\partial_aV_a$, the estimate is
\[
 \begin{aligned}
 |\pair{\nabla\Rea(h_1\overline{h_2})}{\partial_aV_a}|
 &\le C_I\int\langle r\rangle^{-3}|h_1h_2|\\
 &\le C_I\left(\int\langle r\rangle^{-3}|h_1|^2\right)^{1/2}
           \left(\int\langle r\rangle^{-3}|h_2|^2\right)^{1/2}\\
 &\le C_I\|h_1\|_{\mathrm{loc}}\|h_2\|_{\mathrm{loc}}.
 \end{aligned}
\]
A cutoff at radius $R$ produces an error bounded by the same
weighted integral on $R<r<2R$, which tends to zero. Thus this
integration by parts also defines the pairing on $X$.
Insert the gradients from Lemma~\ref{lem:form-gradients} to obtain
\eqref{eq:polarized-cubic}; take $u=\veps$, $f=\epsn$ and
$h=h_1=h_2=\eps$ for \eqref{eq:B-cubic}.

The linear source pairings are justified by
\[
 \begin{aligned}
 (\diver w,p_n)&=-(w,\nabla p_n),&
 (\nabla f,p_v)&=-(f,\diver p_v),\\
 (\Dtwo f,p_n)&=-(f,y\cdot\nabla p_n),&
 (\Dtwo w,p_v)&=-(w,y\cdot\nabla p_v),\\
 \pair{\Lambda h}{p_\eps}
   &=\pair h{p_\eps}+\pair{\nabla h}{y p_\eps},&
 \pair{-\Delta h}{p_\eps}&=\pair{\nabla h}{\nabla p_\eps}.
 \end{aligned}
\]
For example the last two right sides are bounded by
\[
 C_d\|h\|_{\mathrm{loc}}\|p_\eps\|_{H^1_d}
       +\|\nabla h\|_2\|y p_\eps\|_2,
 \qquad \|\nabla h\|_2\|\nabla p_\eps\|_2,
\]
respectively. For wave dilations, the reserve exponential weight
controls $y\nabla p_n,y\nabla p_v$ in $L^2$; the stated
polynomial-tail bounds control the finite profile tests. These are
continuous tested dualities in the stated local norms.
\end{proof}

\begin{proof}[Proof of Lemma~\ref{lem:form-gradients}]
Write the wave-linear test in \eqref{eq:direct-test} as $p$.
The Schr\"odinger gradient of $B_\beta/\beta$ is an
imaginary multiple of $p$. Its density and velocity gradients
are the adjoint tests of this wave-linear map applied to
$\Ima\eps$. We give their operator bounds explicitly.

\paragraph{The density operator.}
The commutator expansion from \eqref{eq:Fa} is
\[
\begin{aligned}
 R_a\mathcal M_a={}&R_a[(2-\Lambda)P_a]-R_aA_aR_aP_a
       +R_a^2A_aR_aP_a\\
 &+R_a^2[(y\cdot\nabla-1)P_a]
       -\frac{(R_a\chi_a)\otimes P_a^2}{M}.
\end{aligned}
\]
Profiles adjacent to an inverse mean multiplication; the last
term is rank one. This maps $L^2$ to $H^2_d$. The adjoint
terms are localized multipliers times elliptic images of
$\eps$ or $R_a\eps$. For example,
\[
 (R_a[(2-\Lambda)P_a])^*\eps=[(2-\Lambda)P_a]R_a\eps,
 \qquad (R_aA_aR_aP_a)^*\eps=P_aR_aA_aR_a\eps.
\]
Use
\begin{equation}\label{eq:AR-local}
 A_aR_a\eps=2\Pi_a\eps-[2+2N_a+y\cdot\nabla N_a]R_a\eps
\end{equation}
and Lemma~\ref{lem:RX} to bound these adjoints in $H^1_d$
by $C_I\|\eps\|_{\mathrm{loc}}$. The rank-one tests are
localized and satisfy the same bound.

\paragraph{The velocity operator.}
The inverse composition and its adjoint satisfy
\begin{equation}\label{eq:S-adjoint}
 [-R_a^2(P_a\diver\,\cdot)]^*\eps=\nabla(P_aR_a^2\eps).
\end{equation}
The composition maps $L^2$ vectors to $H^3_d$, since a divergence
costs one derivative and the two inverses restore four.
Lemma~\ref{lem:RX} gives local $H^5$ control of $R_a^2\eps$.
Multiplication by $P_a$ therefore gives the adjoint and its
divergence the required weighted norms.

\paragraph{The branch term.}
The remaining part is
$\chi_a\mathcal J_a(\epsn,\epsv)/M'$. Its adjoint tests are
the common coefficient $-2(\Ima\eps,\chi_a)/M'$ times
$\partial_a(N_a+P_a^2)$ and $c_0^{-2}\partial_aV_a$, respectively,
for the gradients of $B_\beta/\beta$.
These smooth profiles have, respectively, $r^{-3}$ and
$r^{-2}$ tails; their differentiated tails have the orders in
\eqref{eq:polynomial-tests}. The coefficient is bounded by
$C_I\|\eps\|_{\mathrm{loc}}$.

All constants are uniform for $a\in I$, since $M'$ is bounded
away from zero there. The finite parameter derivatives use the
same profile estimates and differentiated reduced inverse.
The prefactor $\beta$ in \eqref{eq:B} is explicit, so these
bounds prove \eqref{eq:gradient-class} directly.
\end{proof}

\section{Spectral estimates and their verification}\label{app:certificate}

This appendix supplies the common energy and dissipation
inequalities used in Sections~\ref{sec:linear-energy},
\ref{sec:energy} and~\ref{sec:nonradial}. A radial integral bound
gives Schr\"odinger coercivity. The radial and first angular
density estimates use interval enclosures of their kernels.
For the first angular order, we compress the operator to the
orthogonal complement of the translation directions before
deducing a gap. The higher angular density and velocity estimates
are analytic.
The auxiliary energy inequalities are also proved here. They control
the size of the error, separately from the dissipation inequalities
that control its time derivative.

\subsection{The limiting radial density kernel}

With $A_Q=[L_-,\Lambda]=-2\Delta+2rQQ'$, the density dissipation
operator at the limiting profile is
\[
 \mathcal D_Q=1+QR_QA_QR_QQ-[(2-\Lambda)Q]R_QQ
                         -QR_Q[(2-\Lambda)Q].
\]
Put $M_2=\int_0^\infty rQ^2\dd r$.
For $P(r)=M_2^{-1}\int_0^r tQ^2\dd t$ and $T(r)=1-P(r)$,
the reduced inverse uses the kernel
\[
 H(r,t)=\int_0^\infty
 \frac{[P(v)-\mathbf1_{r<v}][P(v)-\mathbf1_{t<v}]}{vQ(v)^2}\dd v.
\]
Relative to radial measure $t\dd t$, the kernel of $\mathcal D_Q-1$, for $r\le t$,
is
\begin{equation}\label{eq:density-kernel}
\begin{aligned}
 S(r,t)={}&-2[\ell(r)+\ell(t)]Q(r)^2Q(t)^2H(r,t)\\
 &+Q(t)^2P(r)-Q(r)^2T(t)\\
 &-\frac{Q(r)^2Q(t)^2}{4M_2}
          \left[r^2+t^2-\frac{2}{M_2}\int_0^\infty u^3Q(u)^2\dd u\right],
 \qquad \ell=-rQ'/Q.
\end{aligned}
\end{equation}
The last line is the contribution of the ground-state projection.
The interval estimate below gives
\begin{equation}\label{eq:HS}
 \int_0^\infty\!\int_0^\infty rt\abs{S(r,t)}^2\dd r\dd t
 <0.88761<0.9.
\end{equation}
Since $0.88761<9/10<(19/20)^2$, the Hilbert--Schmidt
bound gives $\norm{\mathcal D_Q-1}_{\mathrm{op}}<19/20$ and hence
$\mathcal D_Q\ge1/20$.
This estimate holds on the full real radial $L^2$ space.
For the nonradial theorem, the separate first-mode enclosure in
Appendix~\ref{sec:angular-certificate} is used only after the
projection in Proposition~\ref{prop:low-angular-density}.
It is combined with the analytic higher-mode estimate.

The interval construction below bounds the Schr\"odinger integral
and the density kernel, including the origin and exterior regions.
The decimal ODE enclosure
uses $65$ significant digits and directed rounding. Double integrals
are accumulated in binary64 intervals enlarged outwards after each
elementary operation.

\subsection{Interval enclosure of the NLS ground state}

Start with the rational peak interval
\[
 [2.20620086464,2.20620086466].
\]
For a peak in this interval the
radial integral equation on the complex disk $\abs r\le1/10$
is a contraction in the ball $\abs{Q-Q(0)}\le1/10$.
Indeed $\abs Q\le12/5$ there, and its Lipschitz constant is at most
$(1+3(12/5)^2)/400<1$; the integral map also preserves that ball.
If $Q=\sum_{j\ge0}a_jr^{2j}$, its coefficients obey
\[
 a_{j+1}=\frac{a_j-\sum_{i+k+\ell=j}a_i a_k a_\ell}
                   {4(j+1)^2}.
\]
At $r_0=1/100$, evaluate through $j=20$ and bound the tails by
the Cauchy estimate $\abs{a_j}\le(12/5)10^{2j}$, including the
differentiated geometric sum for $Q'$. This gives rigorous initial
intervals despite the radial equation's apparent origin singularity.

\paragraph{Validated integration away from zero.}
Use steps of length at most $\min(1/100,r/10)$ for
$q'=p$, $p'=q-q^3-p/r$, where $q=Q,p=Q'$.
A Picard tube is an interval rectangle containing its integral
image throughout the step. Check this inclusion, then intersect
with eight integral images. The order-24 Taylor coefficients at
the left endpoint and the same recurrence on the whole tube bound
the endpoint value and the next derivative. In this ODE calculation
write $Y=(q,p)$ and $F(r,Y)=(p,q-q^3-p/r)$; a tube $\mathcal T$
is accepted only if
\[
 Y(r)+[0,h]F([r,r+h],\mathcal T)\subset\mathcal T.
\]
The polynomial coefficients are computed recursively from
\[
 (j+1)Y_{j+1}
   =[z^j]F\left(r+z,\sum_{\ell=0}^jY_\ell z^\ell\right).
\]
Here $[z^j]$ selects a coefficient, and every operation is interval
arithmetic. The next derivative, evaluated on the whole tube,
gives the Taylor remainder
\[
 Y(r+h)\in\sum_{j=0}^{24}Y_jh^j
                  +\frac{h^{25}}{25!}\,Y^{(25)}([r,r+h]).
\]
Thus endpoint enclosures and whole-cell tubes have distinct roles:
the first starts the next step; the second bounds the quadratures.

For the lower peak the integration certifies $Q,Q'>0$ at a
radius below $14.30$, with positivity on the preceding tubes. Thus
the initially decreasing solution has reached a positive minimum.
For the upper peak the integration certifies $Q<0$ before $14.34$,
subdividing cells near the first zero. A positive first minimum
has $0<Q<1$, where the quantity
$\tfrac12(Q')^2-\tfrac12Q^2+\tfrac14Q^4$ is negative. This quantity
is nonincreasing, so a solution with such a minimum cannot later
cross zero. The certified events are therefore the disjoint open
undershoot and overshoot alternatives. A peak between them that
belongs to neither alternative gives a positive decreasing solution
for all radii: $Q$ is bounded by its initial value and the same
monotone quantity bounds $Q'$. Its limit is zero. A positive limit
other than one contradicts $(rQ')'=r(Q-Q^3)$, and convergence to one
from above would make $rQ'$ negative and decreasing, forcing an
unbounded decrease of $Q$. Uniqueness of the positive decaying
ground state \cite{Kwong} identifies this solution with $Q$.

Integrating the full peak interval to $R=8$ gives
\[
 0.0003989159<Q(8)<0.0004233731.
\]
Let $q_-,q_+$ denote the corresponding endpoint enclosures,
$\alpha_-=9/10$ and $\alpha_+=11/10$ in this appendix only. Exterior
comparison for $-\Delta+1-Q^2$ gives, for $x\ge0$,
\begin{equation}\label{eq:cert-tail}
 q_-\e^{-\alpha_+ x}\le Q(R+x)\le q_+\e^{-\alpha_- x},\qquad
 -rQ(r)Q'(r)\le q_+^2\e^{-2\alpha_- x}(r/\alpha_-+1/\alpha_-^2).
\end{equation}
For the upper barrier the residual coefficient is
$1-\alpha_-^2+\alpha_-/r-Q^2>0$, since $Q(R)^2<1-\alpha_-^2$.
For the lower barrier it is
$1-\alpha_+^2+\alpha_+/r-Q^2<0$ at $r\ge8$.
The derivative estimate follows by integrating the radial equation
from infinity.

\paragraph{Radial Schr\"odinger coercivity.}
We use the same ground-state enclosure to bound the potential in
$A_Q$. For comparison, \cite[Section~2.4, Step~1]{FMR} proves
nonnegativity of the radial form of $-\Delta+3rQQ'$ on $Q^\perp$.
In this paragraph set
\[
 M_1=\int_0^\infty rQ(r)\dd r,\qquad
 F(t)=M_1^{-1}\int_0^t rQ(r)\dd r,\qquad
 C(t)=\int_0^t rQ(r)^2\dd r-\tfrac12t^2Q(t)^2.
\]
Here $C(t)=\int_0^t r(-rQQ')\dd r$ and $C(\infty)=M_2$.
For smooth radial tests with $\pair\eps Q=0$, subtract the mean of
$\eps$ with respect to $rQ(r)\dd r/M_1$ from the fundamental theorem
of calculus. This gives
\[
 \eps(r)=\int_0^\infty [F(t)-\mathbf1_{r<t}]\,\eps'(t)\dd t .
\]
Cauchy--Schwarz in $t\dd t$ followed by integration against
$r(-rQQ')\dd r$ bounds the potential term by
\begin{equation}\label{eq:nls-radial-positive}
 \int_0^\infty r(-rQQ')|\eps(r)|^2\dd r
 \le\left(\int_0^\infty
 \frac{F(t)^2(M_2-C(t))+(1-F(t))^2C(t)}{t}\dd t\right)
 \int_0^\infty t|\eps'(t)|^2\dd t.
\end{equation}

The integral in parentheses is bounded using the whole-cell tubes
above. Accumulate interval sums of $rQ$ and $rQ^2$ from zero,
including all partial cells, and use these in $F$ and $C$.
For the missing moments at $R=8$, the weaker consequence
$Q(R+x)\le q_+\e^{-x/2}$ of \eqref{eq:cert-tail} suffices.
The interval sum on $[r_0,R]$ is less than $0.780706904$.
For $t\le r_0$, the origin bounds $F(t)\le Q(0)t^2/(2M_1)$
and $C(t)\le Q(0)^4t^4/8$ give the contribution
\[
 \frac{r_0^4}{4}\left(\frac{Q(0)^2M_2}{4M_1^2}
                                +\frac{Q(0)^4}{8}\right).
\]
For $t=R+x$, the same exponential estimate gives
\[
 \begin{aligned}
 M_2-C(t)&\le q_+^2\e^{-x}(t+1+t^2/2),\\
 1-F(t)&\le q_+\e^{-x/2}(2t+4)/M_1 .
 \end{aligned}
\]
Use $1/t\le1/R$, $F\le1$ and $C\le M_2$ in
\eqref{eq:nls-radial-positive}. The resulting bound on the tail is
\[
 \frac{q_+^2}{R}\left(\frac{R^2}{2}+2R+3+
                 \frac{M_2}{M_1^2}(4R^2+24R+40)\right).
\]
The origin and tail contributions together are below $0.000004827$.
Thus the coefficient in \eqref{eq:nls-radial-positive} is below
$0.780712<4/5$. Since the angular factor $2\pi$ cancels,
\begin{equation}\label{eq:scalar-virial-gap}
 (A_Q\eps,\eps)_\R\ge\tfrac25\norm{\nabla\eps}_2^2
                 \ge c\norm\eps_{\mathrm{loc}}^2,\qquad \pair\eps Q=0 .
\end{equation}
The last step uses \eqref{eq:constrained-anchoring}. Approximation
extends the estimate to radial $X$ on the same constraint.

\subsection{Interval bounds for the radial density kernel}\label{sec:radial-density-certificate}

On every compact radial cell use its whole tube for $Q,Q'$.
Accumulate upper and lower sums for $\int_0^r tQ^2$ and bound the
missing tails by \eqref{eq:cert-tail}. The origin estimates are
\[
 Q(0)-Q(r)\le Q(0)^3r^2/4,\qquad -Q'(r)\le Q(0)^3r/2.
\]
For the density kernel write, for $r\le t$,
\[
 H(r,t)=A(r)-\int_r^t\frac{P(v)T(v)}{vQ(v)^2}\dd v+C(t),
 \quad A(r)=\int_0^r\frac{P(v)^2}{vQ(v)^2}\dd v,\quad
 C(t)=\int_t^\infty\frac{T(v)^2}{vQ(v)^2}\dd v.
\]
To derive and normalize the kernel, divide the inverse output by
the positive ground state: for a compactly supported radial source
$g$, write $R_Qg=Qv$. The equation and the orthogonality
condition are
\[
 -(rQ^2v')'=rQ\left(g-\frac{\int_0^\infty tQ(t)g(t)\dd t}{M_2}Q\right),
 \qquad \int_0^\infty rQ^2v\dd r=0.
\]
Integration first determines $v'$, and the last condition determines
its additive constant. The resulting kernel of $R_Q$, for the
measure $t\dd t$, is $Q(r)Q(t)H(r,t)$. Indeed the factors in the
definition of $H$ have zero mean for the probability measure
$rQ(r)^2\dd r/M_2$. For $r<t$ they give
\[
 r\partial_rH(r,t)=\frac{P(r)}{Q(r)^2},\qquad
 t\partial_tH(r,t)=-\frac{T(t)}{Q(t)^2}.
\]
These formulas also give the derivative jump and hence the reduced
inverse equation.

The finite-rank term in the commutator is determined by the identities
\[
 L_-(r^2Q)=-4\Lambda Q,\qquad
 R_Q\Lambda Q=-\tfrac14(r^2-\frac{1}{M_2}\int_0^\infty u^3Q(u)^2\dd u)Q.
\]
Since $[L_-,\Lambda]=A_Q$ and
$L_-R_Q=R_QL_-=\Pi_{Q^\perp}$,
\[
\begin{aligned}
 R_QA_QR_Q
 &=\Pi_{Q^\perp}\Lambda R_Q-R_Q\Lambda\Pi_{Q^\perp}\\
 &=[\Lambda,R_Q]
 +\frac{(R_Q\Lambda Q)\otimes Q+Q\otimes(R_Q\Lambda Q)}{\|Q\|_2^2}.
\end{aligned}
\]
For the second equality use $\Lambda^*=-\Lambda$ and
$R_QQ=0$. In radial kernels the denominator is $M_2$,
because the common angular factor cancels. Substitution of
the displayed formula for $R_Q\Lambda Q$ gives the last line
of \eqref{eq:density-kernel}, including its coefficient $1/(4M_2)$. The commutator kernel is
\[
 (r\partial_r+t\partial_t+2)\{Q(r)Q(t)H(r,t)\}.
\]
Combining this kernel with the two derivatives of $H$ above
gives \eqref{eq:density-kernel}.
This derivation is first made on smooth compactly supported sources;
the square-integrable kernel bounds below justify the bounded
operator extension used in the spectral inequality.
Cellwise interval sums enclose $A$, the middle integral, and $C$;
the tail probability $T$ is accumulated backwards from its positive
tail rather than obtained by subtracting nearly equal numbers.
For a cell $I_i=[r_i,r_{i+1}]$, its radial area is exactly
$\omega_i=(r_{i+1}^2-r_i^2)/2$. If $[S]_{ij}$ encloses the kernel
on $I_i\times I_j$, then
\[
 \iint_{I_i\times I_j}rt|S(r,t)|^2\dd r\dd t
       \le\omega_i\omega_j\sup\{|z|^2:z\in[S]_{ij}\}.
\]
Symmetry counts $i<j$ twice and $i=j$ once; the latter uses the
whole square, including the diagonal, rather than a point value.

For the exterior estimates let $M_-$ be the enclosed lower bound
for $M_2$. Since $2\alpha_+-2\alpha_-=2/5$ and
$4\alpha_--2\alpha_+=7/5$, we obtain
\begin{align*}
 T(R+x)&\le\frac{q_+^2}{M_-}\e^{-2\alpha_- x}
              \left(\frac{R+x}{2\alpha_-}+\frac1{4\alpha_-^2}\right),\\
 C(R)&\le\frac{q_+^4}{M_-^2q_-^2R}
       \int_0^\infty\left(\frac R{2\alpha_-}+\frac1{4\alpha_-^2}
                         +\frac x{2\alpha_-}\right)^2
                      \e^{-7x/5}\dd x,\\
 \int_R^{R+x}\frac{PT}{rQ^2}\dd r
 &\le\frac{q_+^2}{M_-q_-^2}\left(\frac1{2\alpha_-}+\frac1{4\alpha_-^2R}\right)
                  x\e^{2x/5},\\
 \int_R^{R+x}\frac{P^2}{rQ^2}\dd r
 &\le\frac{x}{Rq_-^2}\e^{2\alpha_+ x}.
\end{align*}
Insert these bounds and \eqref{eq:cert-tail} in the four terms of
\eqref{eq:density-kernel}. When one variable exceeds $R$, its
absolute value is bounded by an explicit polynomial of degree two
times $\e^{-7x/5}$. When both variables exceed $R$, order
them as $R+x\le R+y$, replace $x$ by $y$ in the nonnegative
polynomials, and retain the same exponential in $y$.
The inner radial measure is then $Ry+y^2/2$.
All remaining integrals are finite sums of
$\int_0^\infty x^j\e^{-\alpha x}\dd x=j!\alpha^{-j-1}$
with rational positive $\alpha$. At the origin the only logarithmic
singularity is multiplied by $\ell(r)+\ell(t)=O(r^2+t^2)$;
$t^2\log(r_0/t)\le r_0^2/2$ bounds it.
The separate enclosed contributions to $\norm{\mathcal D_Q-1}_{\mathrm{HS}}^2$
are
\[
\begin{array}{c|c}
 [r_0,8]^2&<0.887232473711\\
 \text{at least one variable below }r_0\text{, both below }8
      &<0.000369279558\\
 \text{exactly one variable above }8&<3.212\,10^{-9}\\
 \text{both variables above }8&<5.958\,10^{-11}.
\end{array}
\]
Summing the outward-rounded interval bounds gives
$\norm{\mathcal D_Q-1}_{\mathrm{HS}}^2<0.887601756539$.

\subsection{Green kernel in the first angular mode}\label{app:angular-green}
We construct the first angular inverse needed for the density
criterion of Proposition~\ref{prop:low-angular-density}.
Write $M_2=\int_0^\infty rQ^2\dd r$. The amplitudes use measure
$r\dd r$; we consider angular order one throughout this subsection.

Write
\[
 L_{-,1}=-\partial_r^2-r^{-1}\partial_r+r^{-2}+1-Q^2.
\]
Its quadratic form is at least
$\norm u_{L^2(r\dd r)}^2+(1-M_2/2)\int_0^\infty(|u'|^2+r^{-2}|u|^2)r\dd r$ by
Lemma~\ref{lem:angular-bound}, so its inverse is well defined.
Let $u_1$ be its regular zero solution, normalized by
$u_1(r)/r\to1$ at the origin. This solution is positive:
a first zero would, after multiplying the equation by $u_1$ and
integrating up to that zero, contradict the preceding coercivity.
The exterior barriers in
Appendix~\ref{sec:angular-certificate} give
$9/10\le u_1'/u_1\le3/2$ for $r\ge8$. In particular the following
integral converges. Define
\begin{equation}\label{eq:angular-green}
 h_1(t)=\int_t^\infty\frac{\dd s}{s u_1(s)^2},\qquad
 \Gamma_1(r,t)=u_1(r)u_1(t)h_1(\max\{r,t\}).
\end{equation}
Its derivative has the jump
$\partial_r\Gamma_1(t+,t)-\partial_r\Gamma_1(t-,t)=-1/t$.
Together with regularity at zero and decay at infinity this proves
that $\Gamma_1$ is the kernel of $L_{-,1}^{-1}$ for the measure $t\dd t$.

\subsection{Projection in the first angular density mode}
\label{app:projected-density-criterion}

The first angular order
requires the translation constraint before taking the operator norm.
Let $\Gamma_1(r,t)$ be the kernel, for measure $t\dd t$, of
$(-\partial_r^2-r^{-1}\partial_r+r^{-2}+1-Q^2)^{-1}$.
Its construction is given in Appendix~\ref{app:angular-green}.

\begin{proposition}[A projected integral criterion]
\label{prop:low-angular-density}
Define
\begin{equation}\label{eq:angular-density-kernel}
 K_1(r,t)=(r\partial_r+t\partial_t)
                     \{Q(r)\Gamma_1(r,t)Q(t)\},\qquad
 I_1=\int_0^\infty\!\int_0^\infty
                   |K_1(r,t)|^2rt\dd r\dd t.
\end{equation}
The integral is finite.
In angular order one, after imposing
$({f},Q\partial_1Q)=({f},Q\partial_2Q)=0$, one has
\begin{equation}\label{eq:centered-angular-hs-gap}
 (\mathcal D_Q {f},{f})\ge(1-\sqrt{I_1-1})\norm {f}_2^2.
\end{equation}
In particular, $I_1<2$ suffices for the first-mode gap under the translation constraints.
\end{proposition}

\begin{proof}
In angular order one the reduced inverse is $R_Q=L_{-,1}^{-1}$. The commutator
$[L_-,\Lambda]=A_Q$ therefore gives
$[\Lambda,R_Q]=R_QA_QR_Q$ on smooth compactly supported tests.
Integration by parts in the measure $t\dd t$ shows that the kernel
of this commutator is
$(r\partial_r+t\partial_t+2)\Gamma_1$. The two mixed terms in \eqref{eq:full-density-form} have kernel
\[
 -\bigl[(Q-rQ')(r)Q(t)+Q(r)(Q-tQ')(t)\bigr]\Gamma_1(r,t).
\]
The inverse--commutator term has kernel
$Q(r)Q(t)(r\partial_r+t\partial_t+2)\Gamma_1(r,t)$.
Their sum is
\[
 \begin{aligned}
 &Q(r)Q(t)(r\partial_r+t\partial_t+2)\Gamma_1\\
 &\quad-\{2Q(r)Q(t)-rQ'(r)Q(t)-tQ(r)Q'(t)\}\Gamma_1\\
 &\qquad=(r\partial_r+t\partial_t)
                           \{Q(r)\Gamma_1(r,t)Q(t)\}=K_1(r,t).
 \end{aligned}
\]
Thus $\mathcal D_Q=1+K_1$ on that mode; here $K_1$ also denotes
the integral operator with kernel \eqref{eq:angular-density-kernel}.
The estimates in Appendix~\ref{sec:angular-certificate} bound its
square integral at the origin and infinity, as well as on compact
sets. Thus $K_1$ is a self-adjoint Hilbert--Schmidt operator, and
$\norm{K_1}_{\mathrm{op}}\le\sqrt{I_1}$.
The form and projection identities therefore extend by density.

For angular order one it is essential to use the translation constraint
before taking the operator norm. Normalize $Q Q'$ to a unit vector
$e$ in $L^2(r\dd r)$, and let $\Pi_{e^\perp}=1-e\otimes e$ in this paragraph.
Lemma~\ref{prop:density-translation} gives
$(K_1e,e)=-1$. The orthogonal decomposition
$\operatorname{span}\{e\}\oplus e^\perp$ gives four mutually
orthogonal Hilbert--Schmidt blocks. Since $K_1$ is self-adjoint,
\[
\begin{aligned}
 I_1&=|(K_1e,e)|^2+2\|\Pi_{e^\perp} K_1e\|_2^2
                     +\|\Pi_{e^\perp} K_1\Pi_{e^\perp}\|_{\mathrm{HS}}^2,\\
 \|K_1e\|_2^2&=|(K_1e,e)|^2+\|\Pi_{e^\perp}K_1e\|_2^2.
\end{aligned}
\]
Subtracting these equalities yields
\[
 \norm{\Pi_{e^\perp} K_1\Pi_{e^\perp}}_{\mathrm{HS}}^2
 =I_1-2\norm{K_1e}_2^2+|(K_1e,e)|^2
 \le I_1-1,
\]
because $\norm{K_1e}_2\ge1$. On $e^\perp$ the density form is
$1+\Pi_{e^\perp}K_1\Pi_{e^\perp}$, which proves \eqref{eq:centered-angular-hs-gap}.
The sine and cosine amplitudes give identical copies of this argument;
their orthogonality gives the assertion for the full first mode.
\end{proof}

Section~\ref{sec:density-coercivity} uses the interval bound on
$I_1$ to obtain the full-space estimate.

\subsection{Interval bounds in the first angular density mode}
\label{sec:angular-certificate}

We bound the kernel norm $I_1$ in
Proposition~\ref{prop:low-angular-density}, in angular order one.
Use the same Townes enclosure and tail estimate
\eqref{eq:cert-tail}, and enclose the regular solution $u_1$
needed in \eqref{eq:angular-green}.

At the origin write $u_1=rz_1$. Then
\[
 z_1''+\frac3rz_1'=(1-Q^2)z_1,\qquad z_1(0)=1.
\]
The Volterra equation on $|r|\le1/10$ has contraction constant at
most $(1+(12/5)^2)/(800)<1/100$ and gives $|z_1|<2$.
If $z_1=\sum b_jr^{2j}$, its recurrence is
\[
 b_0=1,\qquad
 4(j+1)(j+2)b_{j+1}=[(1-Q^2)z_1]_{r^{2j}}.
\]
Evaluate through $j=20$ at $r_0=1/100$ and use the Cauchy
remainder $|b_j|\le2\,10^{2j}$, including its differentiated
geometric sum. Away from zero integrate the regular first-order
system
\begin{equation}\label{eq:angular-ode-verifier}
 u_1'=v_1/r,\qquad
 v_1'=\{1/r+r(1-Q^2)\}u_1,
 \qquad v_1=ru_1',
\end{equation}
together with $Q,Q'$. Every step has a verified Picard tube and a Taylor
expansion of order 24 with its interval derivative remainder.
The scaled variable $v_1$ avoids an $r^{-2}$ coefficient in this
step. All $u_1$ tubes have strictly positive lower endpoints.

At $R=8$ the enclosures give
\[
 u_1'(R)/u_1(R)\in[0.9442,0.9443],\qquad I_1<1.652520.
\]
The bound for $I_1$ includes all origin and exterior contributions,
as detailed next. Set $\alpha_-=9/10$ in this verification. The
logarithmic derivative $p=u_1'/u_1$ obeys
\[
 p'=1-Q^2+1/r^2-p^2-p/r.
\]
For $r\ge8$, its right-hand side is positive at $p=\alpha_-$ and
negative at $p=3/2$: it is bounded below there by
$1-1/100-\alpha_-^2-\alpha_-/8>0$ at the lower barrier and above by
$1+1/64-9/4<0$ at the upper one. Thus
\begin{equation}\label{eq:angular-tail-barriers}
 \alpha_-\le u_1'/u_1\le3/2,\qquad
 0\le h_1(R)\le\frac1{2\alpha_- R u_1(R)^2}.
\end{equation}
The second inequality follows by integrating
$u_1(R+x)\ge u_1(R)\e^{\alpha_- x}$.

Differentiating \eqref{eq:angular-green} gives the quadrature formula,
for $r\le t$,
\begin{equation}\label{eq:angular-kernel-evaluation}
\begin{aligned}
 K_1(r,t)={}&h_1(t)\bigl[r(Qu_1)'(r)Q(t)u_1(t)
                     +Q(r)u_1(r)t(Qu_1)'(t)\bigr]\\
 &-Q(r)u_1(r)\frac{Q(t)}{u_1(t)}.
\end{aligned}
\end{equation}
Use the whole-cell tubes in this expression. Starting from
\eqref{eq:angular-tail-barriers}, accumulate interval sums for
$h_1$ backwards with integrand $1/(r u_1^2)$; on each cell its
enclosure includes every possible partial-cell integral.
Square the resulting interval kernel and multiply by the exact
radial cell areas. Sum the upper triangle twice and each diagonal
square once. Symmetry justifies this rule also in diagonal cells.
The compact-square contribution is less than $1.652517$.

At the origin, for $r\le r_0$, the Volterra equation gives
\[
 99/100\le z_1\le101/100,\qquad |rz_1'|\le1/100,
 \qquad |rQ'|\le Q(0)^3r_0^2/2.
\]
Consequently $|Q(r)u_1(r)|\le a_0r$, $|r(Qu_1)'(r)|\le b_0r$ with the
explicit constants
\[
 a_0=\tfrac{101}{100}Q(0),\qquad
 b_0=Q(0)(\tfrac{101}{100}+\tfrac1{100})
             +\tfrac{101}{200}Q(0)^3r_0^2.
\]
For one variable below $r_0$, insert these bounds in
\eqref{eq:angular-kernel-evaluation} and integrate $r^3$
exactly. For both variables below $r_0$, use
\[
 h_1(t)\le\frac{t^{-2}}{2(99/100)^2}+h_1(r_0).
\]
Thus $|K_1(r,t)|\le C_0r/t+C_1rt$ for $r\le t\le r_0$,
where
\[
 C_0=\frac{a_0b_0}{(99/100)^2}
                 +\frac{a_0Q(0)}{99/100},\qquad
 C_1=2a_0b_0h_1(r_0).
\]
The inequality $(x+y)^2\le2x^2+2y^2$ bounds its full square integral by
\[
 \frac{C_0^2r_0^4}{4}
       +2C_1^2\left(\frac{r_0^{4}}{4}\right)^2.
\]
Including the one-variable regions, the origin error is below
$2.938\,10^{-6}$.

Finally let $q_+$ be the upper bound for $Q(R)$ and put
$C_p=\alpha_-^{-1}+\alpha_-^{-2}R^{-1}$ in this paragraph.
Equation~\eqref{eq:cert-tail} gives
$Q(R+x)\le q_+\e^{-\alpha_- x}$ and
$|Q'(R+x)|\le q_+C_p\e^{-\alpha_- x}$.
For $r\le R\le t$, \eqref{eq:angular-kernel-evaluation} and
\eqref{eq:angular-tail-barriers} give
\[
 |K_1(r,t)|\le\frac{q_+}{u_1(R)}\e^{-2\alpha_-(t-R)}
 \left[|Q(r)u_1(r)|\left(1+\frac{3/2+C_p}{2\alpha_-}\right)
                         +\frac{|r(Qu_1)'(r)|}{2\alpha_- R}\right].
\]
For $R\le r\le t$ they give
\[
 |K_1(r,t)|\le q_+^2
       \left(1+\frac{3/2+C_p}{\alpha_-}\right)\e^{-2\alpha_-(t-R)}.
\]
Integrating their squares uses only polynomial exponential moments.
The combined exterior error is less than $1.174\,10^{-10}$.
Its inclusion gives the displayed bound on $I_1$.

The angular implementation uses the same interval arithmetic and
Townes enclosure. Together, the compact cells, origin estimates and
exterior bounds cover the full square $(0,\infty)^2$.

\begin{proposition}[Quantitative radial coercivity]\label{prop:verify-spectral}
For radial $\eps\in H^1$ with $\pair\eps Q=0$, the form of $A_Q$
is bounded below by $\tfrac25\|\nabla\eps\|_2^2$.
On real radial $L^2$, the density operator satisfies
$\mathcal D_Q\ge 1/20$.
\end{proposition}

\begin{proof}
The Schr\"odinger coercivity estimate follows from \eqref{eq:scalar-virial-gap}.
For the density, let $f$ be real and radial in $L^2$ and set $z=R_Q(Qf)$. At the
limiting profile the three terms in the density inequality are
\[
 \|f\|_2^2+(A_Qz,z)-2(f,(Q-rQ')z)
        =(\mathcal D_Qf,f).
\]
The displayed interval bound \eqref{eq:HS} implies
\[
 (\mathcal D_Qf,f)\ge
   (1-\|\mathcal D_Q-1\|_{\mathrm{op}})\|f\|_2^2
 \ge(1-\|\mathcal D_Q-1\|_{\mathrm{HS}})\|f\|_2^2
 >\tfrac1{20}\|f\|_2^2.
\]
This holds on all real radial $L^2$. On the density constraint
in Proposition~\ref{ass:spectral}, the reduced-inverse response
solves $L_-z=Qf$, as required in \eqref{eq:spectral-density}.
The minimum of $1/20$ and the constant in
\eqref{eq:scalar-virial-gap} gives $\delta_*$.
\end{proof}

\subsection{Coercivity of the limiting energy}
\label{app:energy-coercivity}

Let $R_Q=(L_-|_{Q^\perp})^{-1}\Pi_{Q^\perp}$, where
$\Pi_{Q^\perp}$ is the orthogonal projection in $L^2$.
The two limiting energy estimates are as follows.

\begin{lemma}[Coercivity at the limiting profile]\label{lem:limiting-energy}
There exists $c>0$ such that, for radial $h\in X$,
\begin{equation}\label{eq:scalar-energy-gap}
 \pair hQ=0\quad\Longrightarrow\quad
 \int(\abs{\nabla h}^2-Q^2\abs h^2)\ge c\norm h_{\mathrm{loc}}^2.
\end{equation}
On real radial $L^2$, the density energy operator satisfies
\begin{equation}\label{eq:B0-density}
 \mathcal B_Q=1-2QR_QQ\ge0,
 \qquad \ker \mathcal B_Q=\operatorname{span}\{Q\Lambda Q\}.
\end{equation}
Consequently $(\mathcal B_Q{f},{f})\ge c\norm {f}_2^2$ for real
radial ${f}$ satisfying $({f},Q^2)=0$.
\end{lemma}

\begin{proof}
\paragraph{Schr\"odinger energy.}
Let $Q>0$ be the radial Townes state,
\begin{equation}\label{eq:townes}
 Q''+r^{-1}Q'-Q+Q^3=0,\qquad Q'(0)=0,\quad Q(\infty)=0.
\end{equation}
We prove \eqref{eq:scalar-energy-gap}.
To exclude a second negative radial state and a finite-gradient
zero resonance, recall $\ell=-rQ'/Q$ and put
$U_0=\ell^2-\ell-r^2$ and $W_0=1+2U_0+r^2Q^2$ in this proof. The radial ODE gives
\begin{align*}
 U_0=0&\ \Longrightarrow\ rU_0'={\ell}+(2{\ell}-1)r^2Q^2>0,\\
 W_0=0&\ \Longrightarrow\ rW_0'=2U_0.
\end{align*}
For sufficiently small $r>0$, one has $U_0<0$ and $W_0>0$. The decaying Bessel expansion at infinity gives
$U_0=-1/2+1/(4r)+O(r^{-2})$ and
$W_0=1/(2r)+O(r^{-2})$.
At a zero of $U_0$, the relations ${\ell}^2-{\ell}=r^2$ and ${\ell}>0$
give ${\ell}>1$, so the first crossing derivative is strictly positive.
Such a crossing is incompatible with negativity near both endpoints.
Thus $U_0<0$. Every zero of $W_0$ would then have strictly negative
crossing derivative, incompatible with its positive endpoint signs.
This proves $W_0>0$.

For $d_0=\partial_r-Q'/Q$ and its adjoint on $L^2(r\dd r)$,
\[
 d_0^*d_0=-\Delta_\rad+1-Q^2,\qquad
 d_0d_0^*-1=-\Delta_\rad+W_0/r^2.
\]
The positive partner excludes spectrum in $(0,1)$ for $d_0^*d_0$:
$d_0$ would take such an eigenfunction to a nonzero eigenfunction of
$d_0d_0^*-1$ with negative eigenvalue, contradicting $W_0>0$.
To exclude a zero resonance, suppose
$(-\Delta-Q^2)h=0$ with $h$ radial, regular at zero and
$\int_0^\infty r|h'|^2\dd r<\infty$.
It suffices to treat real $h$. The logarithmic bound
\eqref{eq:log-bound} and exponential decay of $Q$ imply
$\int_0^\infty rQ^2|h|\dd r<\infty$. Integrating the radial
equation gives
\[
 rh'(r)=-\int_0^r tQ(t)^2h(t)\dd t.
\]
The right side has a limit. A nonzero limit would make
$\int r|h'|^2$ diverge logarithmically. Thus
\[
 rh'(r)=\int_r^\infty tQ(t)^2h(t)\dd t,
 \qquad \int_1^\infty|h'(r)|\dd r<\infty,
\]
and $h$ tends to a constant. The function $g=d_0h$ is bounded
at infinity, is $O(r)$ at zero, and solves
\[
 (rg')'=\frac{W_0}{r}g.
\]
The regular solution space at zero is one-dimensional, with
$g(r)=cr+O(r^3)$. If $c\ne0$, change its sign so that $c>0$.
Then $g$ and $g'$ stay positive, since $W_0>0$. For any fixed
$R>0$,
\[
 rg'(r)\ge Rg'(R)>0,\qquad
 g(r)\ge g(R)+Rg'(R)\log(r/R),\quad r\ge R,
\]
contrary to boundedness. Hence $g=0$. Now $d_0h=0$ gives
$h=cQ$, and $(-\Delta-Q^2)Q=-Q$ forces $c=0$.
The real and imaginary parts give the same conclusion for
complex $h$.

Compactness of multiplication by $Q$ in the anchored gradient
space now yields \eqref{eq:scalar-energy-gap}. More explicitly,
the partner argument leaves only the eigenfunction $Q$ below zero
for $-\Delta-Q^2$, so its form is nonnegative on $Q^\perp$.
If the strict local gap failed, take $h_n\in X\cap Q^\perp$ with
$\norm{h_n}_{\mathrm{loc}}=1$ and form tending to zero.
After passage to a subsequence it converges weakly in $X$ to $h$;
the weighted local mass and $Qh_n$ converge strongly by the
anchoring estimates and local compactness. The limit cannot be zero:
otherwise the potential term tends to zero, the vanishing form
forces $\norm{\nabla h_n}_2\to0$, and the weighted local mass tends
to zero as well, contradicting the normalization.
Lower semicontinuity and nonnegativity on $Q^\perp$ give zero form
at $h$. Varying within that constraint yields
$(-\Delta-Q^2)h=\alpha Q$ in distributions. Testing against $Q$,
using $(-\Delta-Q^2)Q=-Q$ and $\pair hQ=0$, gives $\alpha=0$.
The excluded nonzero finite-gradient zero solution is the required
contradiction.

\paragraph{Density energy.}
For real $u\perp Q$, the sharp critical interpolation inequality
\cite{Weinstein} gives $(L_+u,u)\ge0$, where $L_+=L_--2Q^2$;
see also \cite[Proposition~2.7]{WeinsteinModulation}.
Indeed, put
\[
 g_\tau=\frac{\|Q\|_2(Q+\tau u)}{\|Q+\tau u\|_2},\qquad
 \|g_\tau\|_2=\|Q\|_2.
\]
The sharp inequality and expansion at $\tau=0$ give
\[
\begin{aligned}
 0\le\|\nabla g_\tau\|_2^2-\tfrac12\|g_\tau\|_4^4
   &=\tau^2\left(\|\nabla u\|_2^2+\|u\|_2^2-3\int Q^2u^2\right)
       +O_u(|\tau|^3)\\
   &=\tau^2(L_+u,u)+O_u(|\tau|^3).
\end{aligned}
\]
The constant term vanishes by the Pohozaev identities, and the
linear term is $-2(Q,u)=0$. Division by $\tau^2$ and passage to
zero prove the claim.
For $u=R_Q(Qf)\perp Q$, this gives
\[
 2\norm{Qu}_2^2\le(L_-u,u)=(f,Qu)
                         \le\norm{f}_2\norm{Qu}_2.
\]
Hence $\norm{Qu}_2\le\norm{f}_2/2$ and
$(L_-u,u)\le\norm{f}_2^2/2$. This proves nonnegativity in \eqref{eq:B0-density}.
The radial invertibility needed to identify the constrained
nullspace is \cite[Proposition~4.1]{GM1}. We recall its
nondegeneracy argument in our normalization. For radial nondegeneracy, Sturm theory shows that a radial zero
eigenfunction would have one node, since the quadratic form has at
most one negative direction and is negative on $Q$. The identities
\[
 (-\Delta+1-3Q^2)Q=-2Q^3,\qquad
 (-\Delta+1-3Q^2)\Lambda Q=-2Q
\]
give zero $Q^3$ and $Q$ moments for that eigenfunction. If $r_0$
is its node, the eigenfunction and $Q^2-Q(r_0)^2$ change sign at
the same radius, so their product has a fixed nonzero sign. Its
pairing with $Q>0$ cannot vanish, contradicting those two moments.
For clarity, if $\mathcal B_Q{f}=0$, set $u=R_Q(Q{f})$.
The definition of $\mathcal B_Q$ and the reduced inverse give
\[
 {f}=2Qu,\quad (u,Q)=0,\quad
 L_+u=-\frac{({f},Q^2)}{\|Q\|_2^2}Q.
\]
Radial invertibility of $L_+$ and $L_+\Lambda Q=-2Q$ then imply
\[
 u=\frac{({f},Q^2)}{2\|Q\|_2^2}\Lambda Q,\qquad
 {f}\in\operatorname{span}\{Q\Lambda Q\}.
\]
Equivalently, the null direction of the radial constrained
$L_+$ form is $\Lambda Q$; this is distinct from an unconstrained
zero eigenfunction of $L_+$.
Finally,
\[
 (Q\Lambda Q,Q^2)
 =\int Q^4+\tfrac14\int y\cdot\nabla(Q^4)
 =\tfrac12\int Q^4\ne0.
\]
Thus $({f},Q^2)=0$ removes the kernel.
If the strict constrained gap failed, a normalized constrained
sequence with $(\mathcal B_Qf^{(k)},f^{(k)})\to0$ would satisfy
$\|\mathcal B_Qf^{(k)}\|_2\to0$, since $\mathcal B_Q\ge0$ and is bounded.
The compactness of $1-\mathcal B_Q$ then makes
$f^{(k)}=(1-\mathcal B_Q)f^{(k)}+o_{L^2}(1)$ converge strongly along a
subsequence to a unit kernel vector satisfying that moment, a
contradiction. This proves the strict density gap.
\paragraph{Continuation to the small positive branch.}
We supply the continuity estimates used in the main proof. Put
\[
 p_a=P_a/\sqrt{M(a)},\qquad
 H_a=L_a+p_a\otimes p_a,\qquad
 p_0=Q/\|Q\|_2,\qquad H_0=L_-+p_0\otimes p_0 .
\]
At the limit $H_0$ is invertible. The profile bounds give
\[
 \|N_a+Q^2\|_\infty\to0,\qquad
 \|y\cdot\nabla(N_a+Q^2)\|_\infty\to0,\qquad
 p_a\to p_0\quad\hbox{in }H^2_{d'} .
\]
For the potential limits, first work on a fixed ball and then
use the uniform polynomial tails on its complement.
Consequently
\[
 \begin{aligned}
 \|H_a-H_0\|_{H^2\to L^2}&\to0,\\
 H_a^{-1}-H_0^{-1}
     &=H_a^{-1}(H_0-H_a)H_0^{-1},\\
 R_a&=H_a^{-1}-p_a\otimes p_a .
 \end{aligned}
\]
The Neumann series is uniform once
$\|(H_a-H_0)H_0^{-1}\|<1/2$. The same argument in
$H^2_d\to L^2_d$, with the weighted elliptic inverse and $d<d'<1$,
applies to localized sources: the multiplication-potential
difference has the same supremum bound, and the rank-one
difference converges in the weighted norm. It follows that
\[
 \|R_aP_a-R_QQ\|_{L^2\to H^2_d}\longrightarrow0 .
\]
In particular $1-2P_aR_aP_a\to\mathcal B_Q$ in $L^2$ operator norm.
Applying
$A_a=-2\Delta-y\cdot\nabla N_a$ to $R_aP_a$ proves the analogous
operator-norm convergence of the density derivative form.
For the Schr\"odinger potentials in the energy and dissipation forms, local convergence
and the bound
\[
 \int_{r>R}\langle r\rangle^{-3}|h|^2
          \le C_pR^{p-3}\|h\|_{\mathrm{loc}}^2,\qquad 2<p<3,
\]
give convergence in form norm on $X$.

To transfer the radial density estimate, take $f\in L^2_{\mathrm{rad}}(\R^2;\R)$
with $(f,P_a^2)=0$ and set
\[
 f^{\,0}=f-
          \frac{(f,Q^2)}{\|Q^2\|_2^2}Q^2 .
\]
Then
\[
 (f^{\,0},Q^2)=0,\qquad
 \|f-f^{\,0}\|_2
 \le\frac{\|P_a^2-Q^2\|_2}{\|Q^2\|_2}\|f\|_2
                      =o_a(1)\|f\|_2 .
\]
For the Schr\"odinger estimate, take radial $h\in X$ with $\pair h{P_a}=0$ and set
\[
 \widehat h=h-\frac{\pair hQ}{\|Q\|_2^2}Q,
 \qquad \pair{\widehat h}Q=0.
\]
Localized profile convergence and the weighted testing estimate give
\[
 |\pair hQ|=|\pair h{Q-P_a}|
             \le o_a(1)\|h\|_{\mathrm{loc}},\qquad
 \|h-\widehat h\|_{\mathrm{loc}}
             \le o_a(1)\|h\|_{\mathrm{loc}}.
\]
Boundedness of the forms then transfers both strict limiting
gaps. Thus, after fixing $a_*>0$ sufficiently small,
\[
 \begin{aligned}
 ((L_a-1)h,h)_\R&\ge c\|h\|_{\mathrm{loc}}^2
                                  &&(\pair h{P_a}=0),\\
 (f,(1-2P_aR_aP_a)f)&\ge c\|f\|_2^2
                                  &&((f,P_a^2)=0).
 \end{aligned}
\]
The derivative gaps transfer in exactly these same norms from
Proposition~\ref{ass:spectral}; this uses the convergence just
proved for the derivative forms.
For radial $h$ without the constraint, write
$h=h^\perp+\pair h{P_a}P_a/M$, where $\pair{h^\perp}{P_a}=0$.
Boundedness of the forms and Young's inequality imply
\[
 \begin{aligned}
 (A_ah,h)_\R&\ge c_I\|h\|_{\mathrm{loc}}^2
                                      -C_I|\pair h{P_a}|^2,\\
 ((L_a-1)h,h)_\R&\ge c_I\|h\|_{\mathrm{loc}}^2
                                      -C_I|\pair h{P_a}|^2.
 \end{aligned}
\]
These are the defect estimates used when the quadratic phase is
nonzero.

\end{proof}

\subsection{Radial velocity coercivity}\label{app:velocity-forms}
In \eqref{eq:velocity-block} the radial components of $\nabla Q^2$ and
$\Dtwo(-yQ^2)$ are, respectively,
\[
 g_1=2QQ',\qquad g_2=-r(3Q^2+2rQQ'),\qquad
 \kappa_0=\int_0^\infty rQ^4\dd r.
\]
The common angular factor $2\pi$ cancels in the rank-two form.

The least spectral value of the rank-two modification in
\eqref{eq:velocity-block} is
\begin{equation}\label{eq:velocity-exact-gap}
 1+\frac{\int r g_1g_2\dd r-
           \sqrt{\int r g_1^2\dd r\int r g_2^2\dd r}}{2\kappa_0}.
\end{equation}
We bound its two moments using only the ground-state equation and
the sharp interpolation inequality.

\begin{lemma}[Radial velocity coercivity]\label{lem:velocity-analytic}
The velocity estimate \eqref{eq:spectral-velocity} holds with $\epsv=w$
for every $w\in L^2_\rad(\R^2;\R^2)$.
\end{lemma}

\begin{proof}
\paragraph{Shape and moment bounds.}
For $\ell=-rQ'/Q$, the proof of
Lemma~\ref{lem:limiting-energy} gives $\ell^2-\ell-r^2<0$.
The radial ODE yields
$r{\ell}'={\ell}^2-r^2+r^2Q^2$. Differentiation and elimination of $Q^2$ yield
\[
 r(r\ell'-2\ell)'=2\ell(\ell^2-\ell-r^2)<0.
\]
The even origin expansion gives $r{\ell}'-2{\ell}\to0$ as $r\downarrow0$.
Consequently
\begin{equation}\label{eq:m-dilation-inequality}
 r{\ell}'(r)<2{\ell}(r)\qquad(r>0).
\end{equation}

Write $M_2=\int_0^\infty rQ^2\dd r$ in this proof.
The Pohozaev identities give $\kappa_0=2M_2$.
The variational characterization in the sharp interpolation
inequality \cite{Weinstein} gives
\[
 \frac{\norm Q_2^2}{2}
 =\inf_{0\ne g\in H^1(\R^2)}
          \frac{\norm{\nabla g}_2^2\norm g_2^2}{\norm g_4^4}
 \le\frac{28\pi}{15}.
\]
Indeed, for $g(y)=(1+|y|^2)^{-2}$, direct integration gives
\[
 \norm g_2^2=\frac\pi3,\qquad
 \norm{\nabla g}_2^2=\frac{4\pi}5,\qquad
 \norm g_4^4=\frac\pi7.
\]
Consequently
\begin{equation}\label{eq:analytic-mass-upper}
 M_2\le\frac{28}{15}<2,\qquad \kappa_0<4.
\end{equation}

For the following one-dimensional calculation set
$p(x)=Q(\sqrt x)^2$, $x\ge0$. By \eqref{eq:m-dilation-inequality},
\[
 (\log p)''(x)=\frac{2{\ell}(r)-r{\ell}'(r)}{2r^4}>0,
 \qquad r=\sqrt x>0,
 \qquad (\log p)'(0)=-\frac{p(0)-1}{2}.
\]
The tangent-line inequality for this strictly convex function gives
$p(x)>p(0)\e^{-(p(0)-1)x/2}$ for $x>0$.
Here $p(0)>1$, since $Q$ has a strict maximum at zero and solves
the ground-state equation. Integration yields
\[
 \kappa_0=\frac12\int_0^\infty p(x)^2\dd x
       >\frac{p(0)^2}{2(p(0)-1)}.
\]
Since $\kappa_0<4$, one has $p(0)^2<8(p(0)-1)$. Solving this
quadratic inequality proves the analytic peak bound below.
The same strict convexity gives
$p((x+y)/2)^2<p(x)p(y)$ for $x\ne y$. Integrating over the
positive quadrant and using $x+y=2u$ gives
\[
 4\int_0^\infty u p(u)^2\dd u
 <\left(\int_0^\infty p(x)\dd x\right)^2=4M_2^2.
\]
We have therefore proved
\begin{equation}\label{eq:analytic-Townes-bounds}
 M_2<2,\qquad Q(0)^2<4+2\sqrt2,\qquad
 \int_0^\infty r^3Q^4\dd r<\frac{M_2^2}{2}
                                      =\frac{\kappa_0^2}{8}.
\end{equation}

\paragraph{The velocity moments.}
Integration by parts gives the three identities
\begin{align}
 \int_0^\infty r g_1g_2\dd r-\kappa_0
 &=\int_0^\infty rQ^4(2{\ell}-r{\ell}')\dd r>0,
 \label{eq:velocity-cross-positive}\\
 \frac{Q(0)^4}{2}-\int_0^\infty r g_1^2\dd r
 &=\int_0^\infty \frac{Q^4}{r}(2{\ell}-r{\ell}')\dd r>0,
 \label{eq:velocity-first-moment}\\
 3\int_0^\infty r^3Q^4\dd r-\int_0^\infty r g_2^2\dd r
 &=\int_0^\infty r^3Q^4(2{\ell}-r{\ell}')\dd r>0.
 \label{eq:velocity-second-moment}
\end{align}
Here are the integrations giving all three identities:
\[
 \begin{aligned}
 \int r{\ell}Q^4&=\kappa_0/2,&
 \int r^2{\ell}'Q^4&=4\int r{\ell}^2Q^4-\kappa_0,\\
 \int\frac{{\ell}Q^4}{r}&=Q(0)^4/4,&
 \int {\ell}'Q^4&=4\int\frac{{\ell}^2Q^4}{r},\\
 \int r^3{\ell}Q^4&=\int r^3Q^4,&
 \int r^4{\ell}'Q^4&=4\int r^3({\ell}^2-{\ell})Q^4.
 \end{aligned}
\]
They follow respectively from integrating derivatives of
$r^2Q^4$, $r^2{\ell}Q^4$, $Q^4$, ${\ell}Q^4$, $r^4Q^4$ and
$r^4{\ell}Q^4$, using $(Q^4)'=-4{\ell}Q^4/r$.
Substitution of $g_1=-2{\ell}Q^2/r$ and $g_2=-r(3-2{\ell})Q^2$
gives \eqref{eq:velocity-cross-positive}--\eqref{eq:velocity-second-moment}.
Only the integral of $(Q^4)'$ retains an origin value; this
produces $Q(0)^4/2$ in the second identity. The other boundary
terms vanish by ${\ell}(r)=O(r^2)$ and exponential decay.
In particular,
\begin{equation}\label{eq:velocity-cross-analytic}
 \int_0^\infty r g_1g_2\dd r>\kappa_0.
\end{equation}
Combining the last two identities with
\eqref{eq:analytic-Townes-bounds} gives
\begin{align}
 \left(\int r g_1^2\dd r\right)
 \left(\int r g_2^2\dd r\right)
 &<\frac{3Q(0)^4}{16}\kappa_0^2
  <\left(\frac92+3\sqrt2\right)\kappa_0^2,\nonumber\\
 \sqrt{\int r g_1^2\dd r\int r g_2^2\dd r}
 &<(\sqrt3+\sqrt{3/2})\kappa_0.
 \label{eq:velocity-product-analytic}
\end{align}
Thus \eqref{eq:velocity-exact-gap} is greater than
$\tfrac12(3-\sqrt3-\sqrt{3/2})$.
The rational comparisons $\sqrt3<26/15$ and
$\sqrt{3/2}<49/40$ give $\sqrt3+\sqrt{3/2}<71/24$;
the gap is therefore greater than $1/48>1/50$.
This proves \eqref{eq:spectral-velocity}.
\end{proof}

\subsection{Velocity coercivity with scale and translation modulation}\label{app:center-velocity-proof}

\begin{proof}[Proof of Proposition~\ref{prop:center-velocity-gap}]
All one-dimensional integrals in this proof use the radial measure
$r\dd r$. Set, only in this proof,
\[
 J_1=\int_0^\infty rQ^2(Q')^2\dd r,\qquad
 J_2=\int_0^\infty rQ^4(Q')^2\dd r.
\]
Recall $\kappa_0=\int rQ^4\dd r=2\int rQ^2\dd r$.
We first control $J_2$ in terms of $J_1$ and $\kappa_0$, without
using a numerical bound on the peak of $Q$.

Multiplying the ground-state equation by $rQ^3$ and $rQ^5$,
respectively, and integrating gives
\begin{equation}\label{eq:center-moment-identities}
 \int_0^\infty rQ^6\dd r=3J_1+\kappa_0,\qquad
 \int_0^\infty rQ^8\dd r=5J_2+3J_1+\kappa_0.
\end{equation}
The boundary terms vanish by smoothness at zero and exponential
decay. Apply the sharp interpolation inequality \cite{Weinstein}
to $Q^2\in H^1(\R^2)$:
\[
 \norm{Q^2}_4^4
 \le\frac{2}{\norm Q_2^2}
          \norm{\nabla(Q^2)}_2^2\norm{Q^2}_2^2.
\]
Since $\norm Q_2^2=\pi\kappa_0$ and
$\norm{\nabla(Q^2)}_2^2=8\pi J_1$, cancellation of $2\pi$
gives $\int rQ^8\dd r\le16J_1$. Thus
\begin{equation}\label{eq:center-moment-bounds}
 5J_2\le13J_1-\kappa_0,\qquad 0<\frac{\kappa_0}{J_1}\le3.
\end{equation}
For the second bound use Cauchy--Schwarz in the first identity of
\eqref{eq:center-moment-identities}:
\[
 \kappa_0^2
 \le\left(\int rQ^2\dd r\right)
       \left(\int rQ^6\dd r\right)
 =\frac{\kappa_0}{2}(3J_1+\kappa_0).
\]

We now estimate $\nabla\partial_j(Q^2)$ and $\partial_j(yQ^2)$. Put $\ell=-rQ'/Q$ in the next
calculation. The analytic shape inequalities
$\ell^2-\ell-r^2<0$ and $r\ell'<2\ell$ imply
\[
 \frac{\ell(r)}{r^2}\le\frac{Q(0)^2-1}{2},\qquad
 \left|\frac{Q'}Q\right|^2
 <1+\frac{\ell}{r^2}\le\frac{Q(0)^2+1}{2}<4.
\]
Here the first bound follows by differentiating $\ell/r^2$ and
using its origin limit. The last bound follows from
$Q(0)^2<4+2\sqrt2<7$, proved in
Lemma~\ref{lem:velocity-analytic}. Consequently
$\int r(Q')^4\dd r<4J_1$.

For each $j$, $\|\partial_j(Q^2)\|_2^2=4\pi J_1$. The ground-state transform for
$L_+$ in angular order one, whose positive radial factor is $-Q'$,
gives the identity
\begin{align}
 \norm{\nabla\partial_j(Q^2)}_2^2
 &=4\pi\int_0^\infty r
       [(Q')^4-Q^2(Q')^2+3Q^4(Q')^2]\dd r\nonumber\\
 &<4\pi(3J_1+3J_2)
 \le \|\partial_j(Q^2)\|_2^2\left(\frac{54}{5}-\frac{3\kappa_0}{5J_1}\right).
 \label{eq:center-A-bound}
\end{align}
To verify it directly, the identity $L_+\partial_jQ=0$ gives
\[
 (L_+(Q\partial_jQ),Q\partial_jQ)
       =\int|\nabla Q|^2|\partial_jQ|^2
       =\pi\int_0^\infty r(Q')^4\dd r.
\]
Expanding $L_+$ on the left and moving the potential term yields
\[
 \|\nabla(Q\partial_jQ)\|_2^2
 =\pi\int_0^\infty r\{(Q')^4-(1-3Q^2)Q^2(Q')^2\}\dd r.
\]
Since $\nabla\partial_j(Q^2)=2\nabla(Q\partial_jQ)$, multiplication by four gives
\eqref{eq:center-A-bound}.

A radial integration by parts also gives
\begin{align}
 \norm{\partial_j(yQ^2)}_2^2
 &=\pi\int_0^\infty r^3[(Q^2)']^2\dd r
 =4\pi\int_0^\infty r\ell^2Q^4\dd r\nonumber\\
 &<2\pi\kappa_0
 =\|\partial_j(Q^2)\|_2^2\frac{\kappa_0}{2J_1}.
 \label{eq:center-B-bound}
\end{align}
The strict inequality is the mixed-moment inequality
\eqref{eq:velocity-cross-positive}. Combining
\eqref{eq:center-A-bound}--\eqref{eq:center-B-bound} gives
\[
 \frac{\norm{\nabla\partial_j(Q^2)}_2^2\norm{\partial_j(yQ^2)}_2^2}{\|\partial_j(Q^2)\|_2^4}
 <\frac{\kappa_0}{2J_1}
       \left(\frac{54}{5}-\frac{3\kappa_0}{5J_1}\right)
 \le\frac{27}{2}.
\]
Indeed the polynomial $x(54-3x)/10$ is increasing on $[0,3]$,
and \eqref{eq:center-moment-bounds} places $\kappa_0/J_1$ in
that interval. Since $(\nabla\partial_j(Q^2),\partial_j(yQ^2))=-2\|\partial_j(Q^2)\|_2^2$, the least spectral value
on each translation plane is therefore greater than
\[
 2-\sqrt{\frac{27}{8}}>\frac3{20}.
\]
The last rational comparison follows by squaring
$\sqrt{27/8}<37/20$. Lemma~\ref{lem:velocity-analytic} gives
the $1/50$ gap on the radial scale-modulation plane. All orthogonal
complementary vectors have coefficient one. The orthogonality of
these planes proves the full $1/50$ bound.
\end{proof}

\subsection{Schr\"odinger estimates in nonzero angular modes}\label{app:angular-comparison}

\begin{lemma}[Angular comparison]\label{lem:angular-bound}
The Townes profile satisfies $\norm Q_2^2<4\pi$. If $h$ has zero
angular mean, then
\begin{equation}\label{eq:angular-potentials}
 \int Q^2|h|^2\le\frac{\norm Q_2^2}{4\pi}\norm{\nabla h}_2^2,
 \qquad
 \int(-rQQ')|h|^2\le\frac{\norm Q_2^2}{4\pi}\norm{\nabla h}_2^2.
\end{equation}
More precisely, on an angular mode of order $m\ge1$, the factor
$\norm Q_2^2/(4\pi)$ can be divided by $m$.
\end{lemma}

\begin{proof}
The mass bound follows from the variational comparison in the proof of
Lemma~\ref{lem:velocity-analytic}.

For a mode $u(r)\e^{\pm\ii m\theta}$ and an interval
$J\subset(0,\infty)$, write
\[
 E_m(u;J)=\int_J (|u'|^2+m^2r^{-2}|u|^2)r\dd r,
 \qquad E_m(u)=E_m(u;(0,\infty)).
\]
For a smooth compactly supported mode, $u(0)=0$. The two
square completions are
\[
 \begin{aligned}
 \int_0^r\left|u'-\frac mt u\right|^2t\dd t
       &=E_m(u;(0,r))-m|u(r)|^2,\\
 \int_r^\infty\left|u'+\frac mt u\right|^2t\dd t
       &=E_m(u;(r,\infty))-m|u(r)|^2.
 \end{aligned}
\]
Adding the nonnegative left sides gives
$E_m(u)\ge2m|u(r)|^2$, hence
\begin{equation}\label{eq:angular-point-evaluation}
 |u(r)|^2\le\frac1{2m}E_m(u).
\end{equation}
Start with smooth compactly supported modes and pass by density.
Integrating this bound against a nonnegative radial potential $V$
gives $\int rV|u|^2\le (2m)^{-1}(\int rV)E_m(u)$.
For the two potentials at hand,
\[
 \int_0^\infty rQ^2\dd r=\frac{\norm Q_2^2}{2\pi},\qquad
 \int_0^\infty r(-rQQ')\dd r
                         =\int_0^\infty rQ^2\dd r.
\]
The second identity follows by integration by parts. Summing the
nonzero Fourier modes proves \eqref{eq:angular-potentials}.
\end{proof}

\begin{proof}[Proof of Proposition~\ref{prop:angular-scalar-velocity}]
Decompose $h$ into its angular mean and its mean-zero part. The
radial operators and radial weights make these components orthogonal
in all the forms. The mean remains orthogonal to $Q$ and satisfies
the Schr\"odinger energy inequality of Lemma~\ref{lem:limiting-energy}
and the Schr\"odinger coercivity estimate of Proposition~\ref{ass:spectral}.
On the other modes,
Lemma~\ref{lem:angular-bound} gives
\[
 \int(|\nabla h|^2-Q^2|h|^2)
 \ge\left(1-\frac{\norm Q_2^2}{4\pi}\right)\norm{\nabla h}_2^2,
\]
and the same argument gives twice this lower bound for $(A_Qh,h)_\R$.
For a mean-zero angular function the local mass obeys
\[
 \int e^{-r}|h|^2
 \le\sup_{r>0}(r^2e^{-r})
       \int r^{-2}|\partial_\theta h|^2
 \le C\|\nabla h\|_2^2.
\]
Combining this with the radial gap and the strict factor
$1-\|Q\|_2^2/(4\pi)>0$ proves both scalar claims.

For the vector form, let ${w}^{\mathrm{rad}}$ be the orthogonal projection of ${w}$
onto the closed subspace of radial no-swirl fields. Both
$\nabla Q^2$ and $\Dtwo(-yQ^2)$ belong to that subspace. The cross
term therefore depends only on ${w}^{\mathrm{rad}}$, while
$\norm {w}_2^2=\norm{{w}^{\mathrm{rad}}}_2^2+\norm{{w}-{w}^{\mathrm{rad}}}_2^2$.
The radial inequality plus this last positive square gives the full
vector estimate with the same constant $1/50$.
\end{proof}

\subsection{Density estimates in nonradial angular modes}\label{app:angular-density}

\begin{proof}[Proof of Proposition~\ref{prop:centered-density-energy}]
The variational inequality $(L_+u,u)\ge0$ for real $u\perp Q$
holds without radial symmetry
\cite[Proposition~2.7]{WeinsteinModulation}. The Cauchy--Schwarz argument preceding
\eqref{eq:B0-density} therefore gives $\mathcal B_Q\ge0$ on the
full real space. This operator is identity minus a compact operator.

The radial kernel of $L_+$ is zero by the argument in
Appendix~\ref{app:energy-coercivity}. In angular order one, the regular
zero solutions are spanned by $\partial_1Q,\partial_2Q$.
The positive factor $-Q'$ is a zero solution in angular
order one, so its ground-state transform makes that form nonnegative.
In higher angular order, comparison with order one adds
$(m^2-1)/r^2$ to a nonnegative form and excludes a zero solution.
Consequently $\ker L_+=\operatorname{span}\{\partial_1Q,\partial_2Q\}$.
For $\mathcal B_Q{f}=0$ set $u=R_Q(Q{f})$. If also $({f},Q^2)=0$,
\[
 {f}=2Qu,\qquad L_-u=Q{f}=2Q^2u,\qquad
 u=\sum_{j=1}^2\alpha_j\partial_jQ.
\]
The translation conditions therefore give
\[
 0=({f},Q\partial_kQ)
   =2\sum_{j=1}^2\alpha_j(Q\partial_jQ,Q\partial_kQ)
   =2\alpha_k\|Q\partial_kQ\|_2^2,\qquad k=1,2.
\]
Thus $u={f}=0$. For a normalized constrained sequence with
$(\mathcal B_Q{f}^{(k)},{f}^{(k)})\to0$, positivity implies
$\|\mathcal B_Q{f}^{(k)}\|_2^2
\le\|\mathcal B_Q\|(\mathcal B_Q{f}^{(k)},{f}^{(k)})\to0$.
Compactness of $1-\mathcal B_Q$ then yields a strongly convergent
subsequence, whose limit has norm one and is a kernel vector
satisfying all three moments. The preceding calculation excludes
such a limit and proves \eqref{eq:centered-density-energy}.
Combining this density bound with
Proposition~\ref{prop:angular-scalar-velocity} and the decomposition
of a field as $h-R_Q(Qf)$ gives coercivity of the full limiting energy.
\end{proof}

\begin{proposition}[Density coercivity in higher angular modes]
\label{prop:high-angular-density}
For every real density ${f}\in L^2(\R^2)$ supported in angular
orders at least two,
\begin{equation}\label{eq:high-angular-density}
 (\mathcal D_Q{f},{f})\ge\frac1{25}\norm {f}_2^2.
\end{equation}
\end{proposition}

\begin{proof}
Young's inequality will reduce the mixed density term to one
radial potential. We bound that potential using weighted moments
of $Q$ and a Hilbert--Schmidt estimate in logarithmic radius,
then absorb it into the angular gradient energy.

\paragraph{Weighted moments.}
Put $\ell=-rQ'/Q$ and, only in this proof,
\[
 J_k=\int_0^\infty r^3Q^4\ell^k\dd r,\qquad 0\le k\le4.
\]
The moment bound in \eqref{eq:analytic-Townes-bounds} and the
improved mass bound \eqref{eq:analytic-mass-upper} give
\[
 J_0<\tfrac12\left(\int_0^\infty rQ^2\dd r\right)^2
 \le\frac{392}{225}.
\]
Integration of $(r^4Q^4\ell^k)'$, using $r\ell'<2\ell$, yields
\[
 J_1=J_0,\qquad
 J_{k+1}=J_k+\frac{k}{4}
       \int_0^\infty r^4Q^4\ell^{k-1}\ell'\dd r
 <\left(1+\frac k2\right)J_k,\quad 1\le k\le3.
\]
All boundary terms vanish by the even origin expansion and
exponential decay. Hence
$J_2<3J_0/2$, $J_3<3J_0$ and $J_4<15J_0/2$.
Since $Q-rQ'=(1+\ell)Q$, expansion gives
\begin{equation}\label{eq:density-potential-moment}
\begin{aligned}
 \int_0^\infty r^3[(Q-rQ')^2-2rQQ']^2\dd r
 &=J_0+8J_1+18J_2+8J_3+J_4\\
 &<\frac{135}{2}J_0<\frac{588}{5}.
\end{aligned}
\end{equation}

\paragraph{The potential bound.}
Use the logarithmic radius $x=\log r$ and define, in this proof,
\[
 \Omega(x)=\left.r^2[(Q-rQ')^2-2rQQ']\right|_{r=\e^x}\ge0.
\]
Then \eqref{eq:density-potential-moment} says
$\int_{\R}\Omega^2\dd x<588/5$.
The Green kernel of $-2\partial_x^2+8$ on $L^2(\R)$ is
$\e^{-2|x-y|}/8$. Therefore the integral operator with kernel
\[
 \frac18\sqrt{\Omega(x)}\e^{-2|x-y|}\sqrt{\Omega(y)}
\]
has squared Hilbert--Schmidt norm bounded by
\begin{align*}
 \frac1{64}\iint\Omega(x)\Omega(y)\e^{-4|x-y|}\dd x\dd y
 &\le\frac1{64}\norm{\e^{-4|\cdot|}}_{L^1(\R)}
                      \norm\Omega_{L^2(\R)}^2\\
 &<\frac{147}{160}<\left(\frac{24}{25}\right)^2.
\end{align*}
The first inequality is Young's convolution inequality. The operator
norm is therefore less than $24/25$. Equivalently,
\[
 \norm{\sqrt\Omega(-2\partial_x^2+8)^{-1/2}}_{L^2\to L^2}^2
 <\frac{24}{25}.
\]
Indeed, the displayed operator times its adjoint is exactly the
integral operator above. Applying it to
$(-2\partial_x^2+8)^{1/2}h$ proves
\[
 \int\Omega|h|^2\dd x
 \le\frac{24}{25}\int(2|h'|^2+8|h|^2)\dd x,
 \qquad h\in H^1(\R).
\]
For a radial amplitude $u$, take $h(x)=u(\e^x)$.
For every angular order $m\ge2$ this gives
\begin{equation}\label{eq:angular-potential-form}
 \int_0^\infty [(Q-rQ')^2-2rQQ']|u|^2r\dd r
 \le\frac{24}{25}\,2\int_0^\infty
               (|u'|^2+m^2r^{-2}|u|^2)r\dd r.
\end{equation}
For $m=2$ this is exactly the change of variables; for $m>2$
the right-hand side only increases. Smooth amplitudes supported
away from zero are dense in this angular energy space, so the
identity extends to $h\in H^1(\R)$.

\paragraph{Absorbing the cross term.}
Let $z=R_Q(Q{f})$ in angular order $m\ge2$. Young's inequality gives
\[
 2|({f},(Q-rQ')z)|
 \le\frac{24}{25}\norm {f}_2^2
                  +\frac{25}{24}\norm{(Q-rQ')z}_2^2.
\]
As $-rQQ'>0$, inequality \eqref{eq:angular-potential-form} implies
\[
 2\int(-rQQ')|z|^2+
       \frac{25}{24}\norm{(Q-rQ')z}_2^2
 \le\frac{25}{24}\int[(Q-rQ')^2-2rQQ']|z|^2
 \le2\norm{\nabla z}_2^2.
\]
The same angular factor occurs on both sides.
Insertion in \eqref{eq:full-density-form} gives the full chain
\[
 \begin{aligned}
 (\mathcal D_Q{f},{f})
 &\ge\tfrac1{25}\|{f}\|_2^2+2\|\nabla z\|_2^2
        -2\int(-rQQ')|z|^2
        -\tfrac{25}{24}\|(Q-rQ')z\|_2^2\\
 &\ge\tfrac1{25}\|{f}\|_2^2.
 \end{aligned}
\]
This proves \eqref{eq:high-angular-density}. Rotation invariance and Parseval's
identity give the same bound on sums of these modes.
\end{proof}

\section{Technical estimates and approximation}
\label{app:energy-details}

The energy identities and their principal cancellations are
derived in Sections~\ref{sec:linear-energy}, \ref{sec:energy}
and~\ref{sec:centered-nonlinear}. This appendix collects the
reduced-inverse commutators, correction-test bounds and
energy-space approximation used there.

\subsection{Linear commutator and correction-test estimates}
\label{app:linear-details}

Fix $a\in I$ and the reference path of Section~\ref{sec:linear-energy}.
Use the corrected fields and real test of
Sections~\ref{sec:linear-response}--\ref{sec:linear-correction}.
All errors are initially smooth and satisfy \eqref{eq:old-slice}.

\paragraph{Reduced-inverse and phase estimates.}
The two projections of the commutator in \eqref{eq:R-commutator}
give
\[
 \begin{aligned}
 \Pi_a[\Lambda,R_a]
       &=R_aA_aR_a-\frac{\chi_a\otimes P_a}{M},\\
 (1-\Pi_a)[\Lambda,R_a]&=-\frac{P_a\otimes\chi_a}{M}.
 \end{aligned}
\]
The second projection follows by testing against $P_a$ and using
$\Lambda^*=-\Lambda$. For the first, use
$L_a[\Lambda,R_a]=A_aR_a+[\Lambda,\Pi_a]$ and apply $R_a$.
The identities are first computed on smooth tests. Their
right-hand sides extend by the localized reduced-inverse bounds.

For general $\beta$ the same product rule gives
\[
 K_\beta\Dtwo-\Lambda K_\beta
 =K_\beta+[\Lambda,R_a]Q_{a,\beta}
                              +R_a(y\cdot\nabla Q_{a,\beta}).
\]
The bracketed profiles in these operator products are multipliers.
Since $Q_{a,\beta}-P_a$ and its finitely many weighted derivatives
are $O_I(\beta)$, these formulas and the reduced-inverse estimates give
\[
 \begin{aligned}
 &\norm{K_\beta+R_aP_a}_{L^2\to H^2_d}\\
 &\quad+\norm{(K_\beta\Dtwo-\Lambda K_\beta)
       -[\Lambda,R_a]P_a+R_a[(2-\Lambda)P_a]}_{L^2\to H^2_d}
                                                  \le C_I\beta,\\
 &\norm{(K_\beta)_s}_{L^2\to H^2_d}\le C_I\beta^2.
 \end{aligned}
\]
The last bound uses $(K_\beta)_s=-\beta^2\partial_\beta K_\beta$
at fixed $a$. The phase difference estimate follows also from
$Q_{a,\beta}-P_a=\int_0^\beta\partial_\theta Q_{a,\theta}\dd\theta$.

\paragraph{The correction derivative.}
The product rule in Section~\ref{sec:linear-correction} leaves
the following projection and test-derivative estimates.
The projection in $L_aR_a=\Pi_a$ has been retained here. Its
omitted kernel contribution is explicitly bounded by
\[
 \begin{aligned}
 &2\beta\left|(\Rea z,(1-\Pi_a)
        [\mathcal M_a\epsn-(c_0a)^{-1}R_a(P_a\diver \epsv)])\right|\\
 &\qquad\le C_I\beta\,|\pair z{P_a}|\,\norm\veps_Y
       \le C_I\beta^2\norm\veps_Y^2.
 \end{aligned}
\]
In fact this contribution is zero on \eqref{eq:old-slice}.
This gives the cancellation against
\eqref{eq:mass-subtracted-identity}--\eqref{eq:diagonal-expansion}.
The correction itself is $O_I(\beta\norm\veps_Y^2)$.
The last three derivative terms are estimated below. First,
\[
 \norm p_{H^2_d}\le C_I\norm\veps_Y.
\]
Along the fixed-orbit linear equation, its derivative is exactly
\begin{align*}
 p_s={}&R_a\mathcal M_a(-\beta\Dtwo\epsn-\lambda\diver\epsv)\\
 &-(c_0a)^{-1}R_a^2\bigl(P_a\diver
 [-\beta\Dtwo\epsv-c_0^2\lambda\nabla
             (\epsn+2\Rea(\overline{Q_{a,\beta}}\eps))]\bigr).
\end{align*}
The localized inverse compositions in Lemma~\ref{lem:linear-tested}
of Appendix~\ref{app:tests} act on the dilation and divergence
terms by duality. Together with $Q_{a,\beta}\eps\in H^1_d$
(Lemma~\ref{lem:anchoring}), they give
\[
 \norm{p_s}_{H^1_d}
 \le C_I(\beta+\lambda)\norm{\veps}_Y
 \le C_I\left(1+\frac1{c_0a_-}\right)\beta\norm{\veps}_Y.
\]
The last factor is absorbed into $C_I$ only after fixing $I$;
this estimate is not uniform as $a_-\downarrow0$.
Using \eqref{eq:good-unknown-linear}, $(p,P_a)=0$ and
$\beta_s=-\beta^2$,
\[
 \begin{aligned}
 |\pair{y\cdot\nabla z}{p}|&\le\norm{\nabla z}_2\norm{yp}_2,\\
 |(\Ima(z_s+\ii L_az),p)|&\le C_I\beta\norm\veps_Y^2,\\
 |\beta_s(\Ima z,p)|+\beta|(\Ima z,p_s)|
                                      &\le C_I\beta^2\norm\veps_Y^2.
 \end{aligned}
\]

Together they prove \eqref{eq:unprojected-diagonal}.

\subsection{Time derivative of the nonlinear correction test}
\label{app:correction-details}

Use the notation of Section~\ref{sec:nonlinear-cancellation}
throughout the bootstrap interval.

\begin{proof}[Proof of Lemma~\ref{lem:correction-test-time}]
Use the definition \eqref{eq:direct-test}, with $a\in I$ and
the actual error equations \eqref{eq:actual-nonlinear-error}.
To verify the second bound, differentiate the displayed formula
for $p$, with $\partial_ap$ holding the two wave errors fixed:
\[
\begin{aligned}
 p_s={}&a_s\partial_ap
       +R_a\mathcal M_a(\epsn)_s
       -(c_0a)^{-1}R_a^2\bigl(P_a\diver(\epsv)_s\bigr)\\
 &+\frac{\chi_a}{M'}\bigl[((\epsn)_s,\partial_a(N_a+P_a^2))
                      +c_0^{-2}((\epsv)_s,\partial_aV_a)\bigr].
\end{aligned}
\]
The coefficient derivative obeys
$\|\partial_ap\|_{H^1_d}\le C_I\|\veps\|_Y$.
Inserting \eqref{eq:actual-nonlinear-error}, every differentiated
wave source carries $\lambda$ or one of the parameter defects.
Lemma~\ref{lem:linear-tested} of Appendix~\ref{app:tests} treats
the linear dilation and divergence sources by duality.
For the nonlinear wave source,
\[
 \|R_a^2(P_a\Delta|\eps|^2)\|_{H^1_d}
       \le C_I\|\eps\|_{\mathrm{loc}}^2.
\]
Indeed, integration by parts gives
\[
 \begin{aligned}
 |\pair{P_a\Delta|\eps|^2}\eta|
 &=|\pair{|\eps|^2}{\Delta(P_a\eta)}|
       \le C_I\|\eps\|_{\mathrm{loc}}^2\|\eta\|_{H^2},\\
 P_a\Delta|\eps|^2&\in H^{-2}_d
       \xrightarrow{R_a}L^2_d
       \xrightarrow{R_a}H^2_d\hookrightarrow H^1_d.
 \end{aligned}
\]
The first inequality uses the weighted local $L^4$ estimate
and the exponential decay of $P_a$ and its derivatives. The
same calculation with the weight inserted proves the indicated
$H^{-2}_d$ bound. The two inverse bounds give the claimed estimate.
For the branch part of $p$, differentiate its scalar coefficient:
\[
 \begin{aligned}
 \partial_s\mathcal J_a(\epsn,\epsv)
 ={}&((\epsn)_s,\partial_a(N_a+P_a^2))
          +c_0^{-2}((\epsv)_s,\partial_aV_a)\\
 &+a_s\bigl[(\epsn,\partial_a^2(N_a+P_a^2))
          +c_0^{-2}(\epsv,\partial_a^2V_a)\bigr].
 \end{aligned}
\]
Insert \eqref{eq:actual-nonlinear-error} and transfer the
spatial derivatives to the profile tests. All linear terms
are bounded by $C_I\lambda\|\veps\|_Y$.
For the quadratic wave source the bound is
\[
 \left|\lambda\int|\eps|^2\diver(\partial_aV_a)\right|
       \le C_I\lambda\int\langle r\rangle^{-3}|\eps|^2
       \le C_I\lambda\|\veps\|_Y^2.
\]
The parameter terms have the same bound by
$|a_s|\le C_I\lambda\|\veps\|_Y$. Hence
\[
 |\partial_s\mathcal J_a(\epsn,\epsv)|
       \le C_I\lambda(\|\veps\|_Y+\|\veps\|_Y^2).
\]
Differentiating the smooth coefficient $\chi_a/M'$ contributes
only $O_I(\lambda\|\veps\|_Y^2)$. Together with the preceding
inverse estimates this proves \eqref{eq:direct-test-time}.\qedhere
\end{proof}

\subsection{Approximation in the energy space}
\label{app:energy-approximation}

\begin{proof}[Proof of Lemma~\ref{lem:energy-solutions}]
\paragraph{Approximation on a common physical interval.}
Fix $s_0<S<S_1<s_*$ and put $t_1=t(S_1)$. Continuity and
positivity give constants $0<\ell<L<\infty$ with
$\ell\le\lambda\le L$ on $[s_0,S_1]$.
Choose smooth initial data converging in
$H^1\times L^2\times L^2$ to the given datum, preserving radial
symmetry and the no-swirl condition in the radial case. Let
$(\psi^k,n^k,v^k)$ be their exact solutions of \eqref{eq:physical}.
The scalar wave datum $n_t(0)=-\diver v(0)$ belongs to
$\dot H^{-1}(\R^2)$, denoted $\widehat H^{-1}$ in
\cite[footnote~1]{BC}. With the norm
$\|h\|_{\dot H^{-1}}=\|(-\Delta)^{-1/2}h\|_2$, the map
$h\mapsto-\nabla\Delta^{-1}h$ is an isometry onto the gradient
subspace of $L^2$. The divergence-free component of
$v$ is its initial value for all times. The same decomposition
applies to the approximants and preserves their higher regularity.
After the fixed sound-speed rescaling, the complex wave variable
$n-\ii(-\Delta)^{-1/2}n_t$ has $L^2$ norm squared
$\|n\|_2^2+\|n_t\|_{\dot H^{-1}}^2$. It obeys the half-wave
equation in \cite[Section~2.1]{SanwalLWP}. At $(s,l)=(1,0)$
in the Sobolev indices of that paper, Theorem~1.1 and its proof
in Section~5 give a lifespan bounded below on bounded sets and
continuous dependence in $C_t(H^1\times L^2)$ for the Schr\"odinger
and complex wave variables. Recovering the real density and
velocity gives continuous dependence in our full energy norm.
For smooth data we use the persistence estimate
\cite[Section~9, eq.~(9.20)]{BC}, whose continuation time is
controlled by the energy norm. Smooth solutions agree with the
energy flow by local uniqueness in the contraction class.

The energy norm of the limiting solution is bounded on $[0,t_1]$.
Take a finite partition of this interval with step size smaller
than a common local lifespan for an energy bound strictly larger
than this bound. At each partition point continuous dependence
gives convergence of the approximating data for the next step.
The larger energy bound provides the same lower lifespan there.
Iteration and local uniqueness give exact approximants on
$[0,t_1]$, with higher regularity, and
\[
 \sup_{0\le t\le t_1}
 \bigl(\|\psi^k-\psi\|_{H^1}+\|n^k-n\|_2
                                      +\|v^k-v\|_2\bigr)\to0.
\]
This construction precedes any use of the modified-energy
inequality. It also passes the smooth mass and Hamiltonian
conservation laws to the energy solution by endpoint continuity.
For the interaction in the Hamiltonian use
$H^1(\R^2)\hookrightarrow L^4(\R^2)$ and $n\in L^2$.

\paragraph{Modulation and rescaled time.}
Apply the physical constraint map \eqref{eq:physical-chart-tests}
with the profiles evaluated at $(x-x_c)/\lambda$ in the nonradial
case and append the two translated tests of \eqref{eq:center-conditions}.
In the complex constraint the full-state term
$\e^{-\ii\gamma}\pair\psi{\phi(n)}$ is independent of $\lambda$ and
$x_c$. Thus $\partial_{\log\lambda}$ and $\lambda\partial_{x_{c,j}}$
act only on the profile term, yielding the pairings with
$\Lambda Q_{a,\beta}-\beta\partial_\beta Q_{a,\beta}$ and
$\partial_jQ_{a,\beta}$, respectively. For the density tests the
same derivatives act on the smooth test functions, yielding
$(\Dtwo N,P_a^2)$ and $(\partial_jN,P_a^2)$ for the scale and
translation columns of the density row, with the derivatives
understood by transposition. The other density rows are obtained
by replacing $P_a^2$ by $\partial_jN_a$. On the constraint slice
these are precisely the coefficients of the tested rate system
after reversing the sign of the complex row. The shift by $\beta$
in the scale rate only affects the forcing. The inverses established
in Sections~\ref{sec:slow-modulation} and~\ref{sec:nonradial}
therefore bound the physical chart inverse uniformly on the compact
interval with $\lambda\ge\ell>0$.
A finite cover by implicit-function neighborhoods therefore
gives continuous approximating parameters, uniformly converging
to the given ones on $[0,t_1]$. Uniqueness near the limiting
parameters joins the local choices; the real phases are chosen
near the same continuous phase. The argument with translation modulation also
gives $x_c^k\to x_c$.

For each approximant define its own rescaled clock by
\[
 s^k(t)=s_0+\int_0^t(\lambda^k(\tau))^{-2}\dd\tau.
\]
For large $k$, $\ell/2\le\lambda^k\le2L$. Consequently
\[
 \sup_{[0,t_1]}|s^k-s|
 \le C_{\ell,t_1}\sup_{[0,t_1]}|\lambda^k-\lambda|\to0,
 \qquad \frac1{4L^2}\le\partial_ts^k\le\frac4{\ell^2}.
\]
Since $s^k(t_1)\to S_1>S$, their inverse clocks $t^k(s)$ are
all defined on $[s_0,S]$ and converge uniformly to $t(s)$.
All rescaled comparisons below are made after composition
with $t^k(s)$ and $t(s)$, respectively.

The rescaling and phase multiplication are strongly continuous
on $H^1$ and $L^2$, with bounded operator norms for scales in
$[\ell/2,2L]$. Their action is uniformly continuous on the compact
image of the limiting trajectory. Combining this fact with
physical convergence and the inverse-clock convergence gives
\[
 \sup_{s_0\le s\le S}
 \left\{|a^k-a|+|\lambda^k-\lambda|+|\gamma^k-\gamma|
       +\|\veps^k-\veps\|_{H^1\times L^2\times L^2}\right\}\to0.
\]
Here the parameters with superscript $k$ are evaluated at
$t^k(s)$. For the decomposition with translation modulation, uniform convergence
of $x_c^k$ and strong continuity of translations give the same
conclusion. The limiting compact interval has a strict
bootstrap margin, so the approximants satisfy the same
bootstrap assumptions for large $k$, with common constants.

\paragraph{Cutoffs in the smooth identities.}
Higher Sobolev regularity suffices for the smooth calculations
even when the approximants have no finite spatial moments.
Use a smooth cutoff $\chi(y/R)$, equal to one on $|y|\le R$ and
zero on $|y|\ge2R$, in the interaction energy and corrected mass.
Compute their integrated identities before sending $R$ to infinity.
The fast Schr\"odinger pairing against $L_N\eps+Q_{a,\beta}\epsn$
still has zero real part under this cutoff. In the wave exchange
the two terms cancel up to the flux with the cutoff gradient:
its factors are $\epsv$ and
$\epsn+2\Rea(\overline{Q_{a,\beta}}\eps)+|\eps|^2$, both in $L^2$.
The kinetic integration by parts produces a cutoff-gradient
term with $\eps_s$ and $\nabla\eps$. For the smooth solution,
$\eps_s-(\lambda_s/\lambda)\Lambda\eps$ is in $L^2$ by the
error equation, including in the presence of translations.
The part containing $y\cdot\nabla\eps$ is bounded by the
annular integral of $|\nabla\eps|^2$, since
$|y\cdot\nabla[\chi(y/R)]|$ is bounded uniformly in $R$.
The other dilation boundary terms are bounded by annular
integrals of the absolute values of the quadratic and cubic
energy densities. For the corrected mass, the analogous terms
use $|z|^2$ and $|z||\nabla z|$. These integrands are in $L^1$.
Translation boundary terms have the cutoff gradient in place
of its dilation derivative. All these errors tend to zero in
the integrated identity by the finite Sobolev bounds on the
compact interval. The correction tests are localized or have
the controlled polynomial tails of \eqref{eq:polynomial-tests};
their boundary terms vanish by the same tested pairings used
below.

\paragraph{Convergence of the energy and tested derivatives.}
The local terms in $\E-H_3$ are continuous in the displayed
norms and parameters. The cubic term is continuous because
\[
 \begin{aligned}
 \left|\int\epsn^k|\eps^k|^2-\int\epsn|\eps|^2\right|
 &\le\|\epsn^k-\epsn\|_2\|\eps^k\|_4^2\\
 &\quad+\|\epsn\|_2(\|\eps^k\|_4+\|\eps\|_4)
                                      \|\eps^k-\eps\|_4,
 \end{aligned}
\]
and $H^1(\R^2)\hookrightarrow L^4(\R^2)$.
The density-dependent inverse and its localized derivatives
are continuous by Lemma~\ref{lem:current-inverse}. It follows
that $\E^k\to\E$ uniformly on $[s_0,S]$.

For the differentiated identities, take the cancellations
before passing to the limit. Their scalar right-hand sides
contain gradient quadratic forms, localized profile tests,
products against localized tests, and the differentiated
inverse terms in \eqref{eq:good-inverse-time}.
For the terms involving the density and velocity derivatives,
use the following continuous interpretations:
\[
 \begin{aligned}
 (\diver\epsv,p)&=-(\epsv,\nabla p),\\
 (\Dtwo\epsn,p)&=-(\epsn,y\cdot\nabla p),\\
 \pair{-\Delta\eps}p&=\pair{\nabla\eps}{\nabla p}.
 \end{aligned}
\]
Lemma~\ref{lem:form-gradients} gives the required weighted
regularity of these tests.
For example, convergence of both arguments gives
\[
 \begin{aligned}
 |(\diver\epsv^k,p^k)-(\diver\epsv,p)|
 &\le\|\epsv^k-\epsv\|_2\|\nabla p^k\|_2
      +\|\epsv\|_2\|\nabla(p^k-p)\|_2\to0.
 \end{aligned}
\]
The same estimate with $y\cdot\nabla p$ treats dilation.
When $p$ is a correction test, convergence in its reserve
weighted norm controls this additional factor of $y$.
For its finite polynomial profile components, the required
weighted gradients converge in $L^2$ by the profile bounds.

The inverse terms require continuity in their density argument
as well. The next density derivative in
\eqref{eq:inverse-mixed-tests} and the mean-value formula give,
for a fixed localized $g$,
\[
 \begin{aligned}
 &\|(R'_{N^k}-R'_N)[\diver V^k]g\|_{H^1_d}
       \le C_I\|N^k-N\|_2\|V^k\|_2\|g\|_{L^2_{d'}},\\
 &\|R'_N[\diver(V^k-V)]g\|_{H^1_d}
       \le C_I\|V^k-V\|_2\|g\|_{L^2_{d'}}.
 \end{aligned}
\]
Here a single common reference $N_a$ may be used locally;
parameter variation is continuous on $I$. If $g=g^k$ also
varies, linearity and the same bound control $g^k-g$.
This applies to $g=Q_{a,\beta}\epsn$. The source $Q_{a,\beta}\diver\epsv$
is treated in $H^{-1}_{d'}$ by
\[
 Q_{a,\beta}\diver\epsv=\diver(Q_{a,\beta}\epsv)-\nabla Q_{a,\beta}\cdot\epsv,
 \qquad \|Q_{a,\beta}\diver\epsv\|_{H^{-1}_{d'}}\le C_I\|\epsv\|_2.
\]
The corresponding eigenfunction pairings are continuous by
\eqref{eq:mixed-spectral}. Thus every term in
\eqref{eq:good-inverse-time}, after pairing with $z$, converges.
In particular, this passage never treats $N_s$ or
$(\epsv)_s$ as an $L^2$ function.

Finally, for the nonlinear source pairings the localized
product estimates in Corollary~\ref{cor:nonlinear-pairings}
are multilinear bounds. Subtracting two such expressions one
argument at a time proves their convergence; for example,
\[
 \begin{aligned}
 |\pair{\epsn\eps}p|
   &\le C_I\|\epsn\|_2\|\eps\|_{\mathrm{loc}}\|p\|_{H^1_d},\\
 |\pair{\nabla|\eps|^2}{p_v}|
   &\le C_I\|\eps\|_{\mathrm{loc}}^2\|p_v\|_{H^1_d}.
 \end{aligned}
\]
The polynomial profile tests use instead
$\int\langle r\rangle^{-3}|\eps|^2
\le C\|\eps\|_{\mathrm{loc}}^2$.
These estimates also cover the differentiated correction and,
after Section~\ref{sec:centered-nonlinear}, its transport work.

\paragraph{Integrated identities and inequalities.}
First apply the preceding continuity to the tested modulation
system. Its matrix is uniformly invertible, so the smooth
outputs $(a_s^k,\lambda_s^k,\gamma_s^k)$ converge uniformly to
continuous functions of the limiting parameters and errors.
Integrating the smooth parameter equations identifies these
functions with the derivatives of $(a,\lambda,\gamma)$.
Thus these parameters are $C^1$. In the nonradial case the
same matrix convergence gives uniform convergence of $\nu^k$,
and integration of $x_{c,s}^k=\lambda^k\nu^k$ identifies the
center derivative. This gives the parameter derivatives before
they enter the limiting energy identity.

We now have strong convergence of the errors and uniform
convergence of the rates, including
$\beta_s=c_0(\lambda a_s+a\lambda_s)$ and, in the nonradial
case, $\nu=x_{c,s}/\lambda$. For $E_s$, use \eqref{eq:action-law}
with $2(\lambda_s/\lambda)H_3$ restored. Its terms are continuous:
the dilation form uses the strong convergence of gradients;
the profile work uses the localized or polynomial tests; and the
only differentiated inverse in $\mathcal F$ has the continuous
representation \eqref{eq:good-inverse-time}.

For the correction, at fixed errors the partial derivatives
$\partial_a B_\beta$ and $\partial_\beta B_\beta$ are continuous
quadratic forms by Lemma~\ref{lem:form-gradients}. In the smooth
identity the remaining derivative is
$2B_\beta(\veps,\veps_s)$, using \eqref{eq:actual-nonlinear-error}
before adding the center transport. Apply the
gradient representatives \eqref{eq:gradient-class} and move
the spatial derivatives onto these tests as in
Corollary~\ref{cor:nonlinear-pairings}. This accounts for each
term of the product rule, including the differentiated wave
source $\nabla|\eps|^2$. For the full-space calculation add the
separate translation contribution \eqref{eq:exact-center-transport-work}.
In it $b_j$ is continuous by the reduced-inverse bounds, and the
wave derivatives in the last term act on the correction
gradients as in \eqref{eq:B-transport-tests}.

Thus the exact scalar right-hand sides $F^k$ for $(\E^k)_s$
converge uniformly on $[s_0,S]$ to a continuous function $F$.
This follows from continuity of the canceled expressions on
the compact limiting trajectory and the uniform convergence
of their arguments. For $s_0\le u\le v\le S$,
\[
 \begin{aligned}
 \E^k(v)-\E^k(u)&=\int_u^v F^k(s)\dd s,\\
 \E(v)-\E(u)&=\int_u^v F(s)\dd s.
 \end{aligned}
\]
The second equality follows by uniform convergence, not by
differentiating a rough energy-space solution. It gives local
absolute continuity of $\E$ and the asserted integrated
identity. The same convergence applied to the smooth
inequality gives
\[
 \E(v)-\E(u)
       \le-\kappa_I\int_u^v\lambda(s)\E(s)\dd s.
\]
All constants are common to the approximants on this compact
interval. Since $S<s_*$ is arbitrary, the conclusions hold on
every compact subinterval of $[s_0,s_*)$.
\end{proof}

\section*{Data and code availability}

The arXiv ancillary files contain the programs, recorded interval
bounds, arithmetic assumptions and reproduction instructions for
the radial Schr\"odinger estimate and the radial and first angular
density estimates. The integral formulas and the origin and exterior
bounds are given in Appendix~\ref{app:certificate}.

\section*{Use of artificial intelligence}

The research questions and main mathematical ideas originate with
the authors. OpenAI Codex was used to assist with revising and
polishing the manuscript, working through and checking intermediate
mathematical derivations, and developing and checking the code for
the computer-assisted spectral verification in
Appendix~\ref{app:certificate}. The authors retain full
responsibility for the mathematical arguments, the computational
verification and the final text.

\end{document}